\documentclass[final,3p,sort&compress]{elsarticle} 

\usepackage{etoolbox}
\usepackage[utf8]{inputenc}
\usepackage{amsmath}
\usepackage{amsfonts}
\mathversion{normal}
\usepackage{graphicx}
\usepackage{tikz}
\usepackage{float}
\usepackage{xcolor}
\usepackage{multirow}

\makeatletter
\@ifclassloaded{elsarticle}{
	\usepackage{amssymb}
	\usepackage{amsthm}
	\usepackage[light]{kpfonts}
	\usepackage{listings}
	\usepackage{caption}
	\usepackage{pdfpages}
	\usepackage[colorlinks = false]{hyperref}
	\usepackage[hyperref,verbose]{backref}
	\renewcommand*{\backrefalt}[4]{%
		\ifcase #1%
			No citations.%
		\or%
			Cited in section #2.%
		\else%
			Cited in sections #2.%
		\fi%
	}
	\journal{{\tt arXiv.org}}
	\newbox\amsbox
	\newenvironment{AMS}{\global\setbox\amsbox=\vbox\bgroup\par}{\egroup} 
	
	\title{Asymptotic preserving $P_N$ methods for haptotaxis equations
		\thanks{This work was financially supported by BMBF in the project GlioMaTh.}}
	
	\author{Gregor Corbin}
	\address{Department of Mathematics, University of Kaiserslautern, \\ 
		P.O. Box 3049, 67653 Kaiserslautern, Germany\\
		{\it corbin@mathematik.uni-kl.de}}
	
	\newcommand{\addreferencestotoc}{\addcontentsline{toc}{section}{References}}
	\newcommand{\externaldir}{./Externals/elsarticle/}
}
{
	\usepackage{textcomp}
	
	\newsiamremark{remark}{Remark}
	\newsiamremark{hypothesis}{Hypothesis}
	\crefname{hypothesis}{Hypothesis}{Hypotheses}
	\newsiamthm{claim}{Claim}
	
	\headers{AP method for haptotaxis}{G. Corbin}
	
	\title{Asymptotic preserving $P_N$ methods for haptotaxis equations\thanks{Submitted to the editors DATE.
			\funding{This work was financially supported by BMBF in the project GlioMaTh.}}}
	
	\author{Gregor Corbin\thanks{AG Technomathematik, Technische Universität Kaiserslautern, Germany
			(\email{corbin@mathematik.uni-kl.de}, \url{https://www.mathematik.uni-kl.de/techno/personen/mitglieder/}).}}
	
	\newcommand{\addreferencestotoc}{} 
	\newcommand{\externaldir}{./Externals/siamart/}
}
\makeatother

\usepackage{standalone}
\usepackage{pgfplots}
\pgfplotsset{compat=1.12}
\usepgfplotslibrary{groupplots} 
\usepgfplotslibrary{colorbrewer}
\newcommand{\externaltikz}[2]{\includegraphics{\externaldir#1}}

\ifcsmacro{thm}{}{}
\ifcsmacro{claim}{}{}   
\ifcsmacro{remark}{}{\newtheorem{remark}{Remark}}
\ifcsmacro{prop}{}{}
\ifcsmacro{ex}{}{}
\ifcsmacro{assumption}{}{}
\ifcsmacro{lemma}{}{\newtheorem{lemma}{Lemma}}
\ifcsmacro{fact}{}{}
\ifcsmacro{definition}{}{\newtheorem*{definition}{Definition}}

\newcommand{\secref}[1]{Section~\ref{#1}}

\newcommand{\lemref}[1]{Lemma~\ref{#1}}

\newcommand{\figref}[1]{Figure~\ref{#1}}
\newcommand{\tabref}[1]{Table~\ref{#1}}

\usepackage[colorinlistoftodos,obeyFinal]{todonotes}

\newlength{\figurehorizontalsep}
\newlength{\figurethreecol}
\newlength{\figuretwocol}
\newlength{\figureonecol}
\newlength{\figuretwocollegend}
\newlength{\figurelegendhorizontalsep}
\newlength{\figureheight}
\newlength{\figurewidth}
\newlength{\figureheightsave}
\newlength{\figurewidthsave}

\newcommand{\withfiguresize}[3]{%
	\setlength{\figureheightsave}{\figureheight}%
	\setlength{\figurewidthsave}{\figurewidth}%
	\setlength{\figureheight}{#1}%
	\setlength{\figurewidth}{#2}%
	#3%
	\setlength{\figureheight}{\figureheightsave}%
	\setlength{\figurewidth}{\figurewidthsave}%
}

\newcommand{\figureticklabelfont}[1]{%
	\ifdim #1 < 3cm%
		\tiny
	\else%
		\ifdim #1 < 5cm%
			\scriptsize
		\else%
			\footnotesize
		\fi%
	\fi%
}

\newcommand{\figurelabelfont}[1]{%
	\ifdim #1 < 3cm%
		\footnotesize
	\else%
		\ifdim #1 < 5cm%
			\small
		\else%
			\normalsize
		\fi%
	\fi%
}
\newcounter{tikzsubfigcounter}[figure]
\renewcommand{\thetikzsubfigcounter}{\the\numexpr\value{figure}+1\relax\alph{tikzsubfigcounter}}

\newcommand{\tikztitle}[1]{ %
	\refstepcounter{tikzsubfigcounter}
	\textbf{(\alph{tikzsubfigcounter})}\space\space #1 
}

\newcounter{tikzsubfigcounterinvisible}[figure]
\makeatletter
\@ifclassloaded{siamart190516}{
	\renewcommand{\thetikzsubfigcounterinvisible}{\thesection.\the\numexpr\value{figure}+1\relax\alph{tikzsubfigcounterinvisible}}
}{
	\renewcommand{\thetikzsubfigcounterinvisible}{\the\numexpr\value{figure}+1\relax\alph{tikzsubfigcounterinvisible}}}
\makeatother

\newcommand{\settikzlabel}[1]{ %
	\refstepcounter{tikzsubfigcounterinvisible} \label{#1}
} 

\usepackage{tikz}
\usepackage{pgfplots}
\pgfplotsset{compat=1.13}
\usetikzlibrary{calc}
\usepgfplotslibrary{colorbrewer}
\usepgfplotslibrary{groupplots}

\pgfplotscreateplotcyclelist{marklist}{%
	mark=x\\
	mark=o\\
	mark=triangle\\
	mark=diamond\\
	mark=square\\
}
\pgfplotsset{
	cycle list/Set1-5,
	cycle multiindex* list={
		marklist\nextlist
		Set1-5\nextlist},
}

\pgfplotsset{
	colormap={cmviridis}{
rgb(0pt)=(0.267004,0.004874,0.329415);
rgb(1pt)=(0.268510,0.009605,0.335427);
rgb(2pt)=(0.269944,0.014625,0.341379);
rgb(3pt)=(0.271305,0.019942,0.347269);
rgb(4pt)=(0.272594,0.025563,0.353093);
rgb(5pt)=(0.273809,0.031497,0.358853);
rgb(6pt)=(0.274952,0.037752,0.364543);
rgb(7pt)=(0.276022,0.044167,0.370164);
rgb(8pt)=(0.277018,0.050344,0.375715);
rgb(9pt)=(0.277941,0.056324,0.381191);
rgb(10pt)=(0.278791,0.062145,0.386592);
rgb(11pt)=(0.279566,0.067836,0.391917);
rgb(12pt)=(0.280267,0.073417,0.397163);
rgb(13pt)=(0.280894,0.078907,0.402329);
rgb(14pt)=(0.281446,0.084320,0.407414);
rgb(15pt)=(0.281924,0.089666,0.412415);
rgb(16pt)=(0.282327,0.094955,0.417331);
rgb(17pt)=(0.282656,0.100196,0.422160);
rgb(18pt)=(0.282910,0.105393,0.426902);
rgb(19pt)=(0.283091,0.110553,0.431554);
rgb(20pt)=(0.283197,0.115680,0.436115);
rgb(21pt)=(0.283229,0.120777,0.440584);
rgb(22pt)=(0.283187,0.125848,0.444960);
rgb(23pt)=(0.283072,0.130895,0.449241);
rgb(24pt)=(0.282884,0.135920,0.453427);
rgb(25pt)=(0.282623,0.140926,0.457517);
rgb(26pt)=(0.282290,0.145912,0.461510);
rgb(27pt)=(0.281887,0.150881,0.465405);
rgb(28pt)=(0.281412,0.155834,0.469201);
rgb(29pt)=(0.280868,0.160771,0.472899);
rgb(30pt)=(0.280255,0.165693,0.476498);
rgb(31pt)=(0.279574,0.170599,0.479997);
rgb(32pt)=(0.278826,0.175490,0.483397);
rgb(33pt)=(0.278012,0.180367,0.486697);
rgb(34pt)=(0.277134,0.185228,0.489898);
rgb(35pt)=(0.276194,0.190074,0.493001);
rgb(36pt)=(0.275191,0.194905,0.496005);
rgb(37pt)=(0.274128,0.199721,0.498911);
rgb(38pt)=(0.273006,0.204520,0.501721);
rgb(39pt)=(0.271828,0.209303,0.504434);
rgb(40pt)=(0.270595,0.214069,0.507052);
rgb(41pt)=(0.269308,0.218818,0.509577);
rgb(42pt)=(0.267968,0.223549,0.512008);
rgb(43pt)=(0.266580,0.228262,0.514349);
rgb(44pt)=(0.265145,0.232956,0.516599);
rgb(45pt)=(0.263663,0.237631,0.518762);
rgb(46pt)=(0.262138,0.242286,0.520837);
rgb(47pt)=(0.260571,0.246922,0.522828);
rgb(48pt)=(0.258965,0.251537,0.524736);
rgb(49pt)=(0.257322,0.256130,0.526563);
rgb(50pt)=(0.255645,0.260703,0.528312);
rgb(51pt)=(0.253935,0.265254,0.529983);
rgb(52pt)=(0.252194,0.269783,0.531579);
rgb(53pt)=(0.250425,0.274290,0.533103);
rgb(54pt)=(0.248629,0.278775,0.534556);
rgb(55pt)=(0.246811,0.283237,0.535941);
rgb(56pt)=(0.244972,0.287675,0.537260);
rgb(57pt)=(0.243113,0.292092,0.538516);
rgb(58pt)=(0.241237,0.296485,0.539709);
rgb(59pt)=(0.239346,0.300855,0.540844);
rgb(60pt)=(0.237441,0.305202,0.541921);
rgb(61pt)=(0.235526,0.309527,0.542944);
rgb(62pt)=(0.233603,0.313828,0.543914);
rgb(63pt)=(0.231674,0.318106,0.544834);
rgb(64pt)=(0.229739,0.322361,0.545706);
rgb(65pt)=(0.227802,0.326594,0.546532);
rgb(66pt)=(0.225863,0.330805,0.547314);
rgb(67pt)=(0.223925,0.334994,0.548053);
rgb(68pt)=(0.221989,0.339161,0.548752);
rgb(69pt)=(0.220057,0.343307,0.549413);
rgb(70pt)=(0.218130,0.347432,0.550038);
rgb(71pt)=(0.216210,0.351535,0.550627);
rgb(72pt)=(0.214298,0.355619,0.551184);
rgb(73pt)=(0.212395,0.359683,0.551710);
rgb(74pt)=(0.210503,0.363727,0.552206);
rgb(75pt)=(0.208623,0.367752,0.552675);
rgb(76pt)=(0.206756,0.371758,0.553117);
rgb(77pt)=(0.204903,0.375746,0.553533);
rgb(78pt)=(0.203063,0.379716,0.553925);
rgb(79pt)=(0.201239,0.383670,0.554294);
rgb(80pt)=(0.199430,0.387607,0.554642);
rgb(81pt)=(0.197636,0.391528,0.554969);
rgb(82pt)=(0.195860,0.395433,0.555276);
rgb(83pt)=(0.194100,0.399323,0.555565);
rgb(84pt)=(0.192357,0.403199,0.555836);
rgb(85pt)=(0.190631,0.407061,0.556089);
rgb(86pt)=(0.188923,0.410910,0.556326);
rgb(87pt)=(0.187231,0.414746,0.556547);
rgb(88pt)=(0.185556,0.418570,0.556753);
rgb(89pt)=(0.183898,0.422383,0.556944);
rgb(90pt)=(0.182256,0.426184,0.557120);
rgb(91pt)=(0.180629,0.429975,0.557282);
rgb(92pt)=(0.179019,0.433756,0.557430);
rgb(93pt)=(0.177423,0.437527,0.557565);
rgb(94pt)=(0.175841,0.441290,0.557685);
rgb(95pt)=(0.174274,0.445044,0.557792);
rgb(96pt)=(0.172719,0.448791,0.557885);
rgb(97pt)=(0.171176,0.452530,0.557965);
rgb(98pt)=(0.169646,0.456262,0.558030);
rgb(99pt)=(0.168126,0.459988,0.558082);
rgb(100pt)=(0.166617,0.463708,0.558119);
rgb(101pt)=(0.165117,0.467423,0.558141);
rgb(102pt)=(0.163625,0.471133,0.558148);
rgb(103pt)=(0.162142,0.474838,0.558140);
rgb(104pt)=(0.160665,0.478540,0.558115);
rgb(105pt)=(0.159194,0.482237,0.558073);
rgb(106pt)=(0.157729,0.485932,0.558013);
rgb(107pt)=(0.156270,0.489624,0.557936);
rgb(108pt)=(0.154815,0.493313,0.557840);
rgb(109pt)=(0.153364,0.497000,0.557724);
rgb(110pt)=(0.151918,0.500685,0.557587);
rgb(111pt)=(0.150476,0.504369,0.557430);
rgb(112pt)=(0.149039,0.508051,0.557250);
rgb(113pt)=(0.147607,0.511733,0.557049);
rgb(114pt)=(0.146180,0.515413,0.556823);
rgb(115pt)=(0.144759,0.519093,0.556572);
rgb(116pt)=(0.143343,0.522773,0.556295);
rgb(117pt)=(0.141935,0.526453,0.555991);
rgb(118pt)=(0.140536,0.530132,0.555659);
rgb(119pt)=(0.139147,0.533812,0.555298);
rgb(120pt)=(0.137770,0.537492,0.554906);
rgb(121pt)=(0.136408,0.541173,0.554483);
rgb(122pt)=(0.135066,0.544853,0.554029);
rgb(123pt)=(0.133743,0.548535,0.553541);
rgb(124pt)=(0.132444,0.552216,0.553018);
rgb(125pt)=(0.131172,0.555899,0.552459);
rgb(126pt)=(0.129933,0.559582,0.551864);
rgb(127pt)=(0.128729,0.563265,0.551229);
rgb(128pt)=(0.127568,0.566949,0.550556);
rgb(129pt)=(0.126453,0.570633,0.549841);
rgb(130pt)=(0.125394,0.574318,0.549086);
rgb(131pt)=(0.124395,0.578002,0.548287);
rgb(132pt)=(0.123463,0.581687,0.547445);
rgb(133pt)=(0.122606,0.585371,0.546557);
rgb(134pt)=(0.121831,0.589055,0.545623);
rgb(135pt)=(0.121148,0.592739,0.544641);
rgb(136pt)=(0.120565,0.596422,0.543611);
rgb(137pt)=(0.120092,0.600104,0.542530);
rgb(138pt)=(0.119738,0.603785,0.541400);
rgb(139pt)=(0.119512,0.607464,0.540218);
rgb(140pt)=(0.119423,0.611141,0.538982);
rgb(141pt)=(0.119483,0.614817,0.537692);
rgb(142pt)=(0.119699,0.618490,0.536347);
rgb(143pt)=(0.120081,0.622161,0.534946);
rgb(144pt)=(0.120638,0.625828,0.533488);
rgb(145pt)=(0.121380,0.629492,0.531973);
rgb(146pt)=(0.122312,0.633153,0.530398);
rgb(147pt)=(0.123444,0.636809,0.528763);
rgb(148pt)=(0.124780,0.640461,0.527068);
rgb(149pt)=(0.126326,0.644107,0.525311);
rgb(150pt)=(0.128087,0.647749,0.523491);
rgb(151pt)=(0.130067,0.651384,0.521608);
rgb(152pt)=(0.132268,0.655014,0.519661);
rgb(153pt)=(0.134692,0.658636,0.517649);
rgb(154pt)=(0.137339,0.662252,0.515571);
rgb(155pt)=(0.140210,0.665859,0.513427);
rgb(156pt)=(0.143303,0.669459,0.511215);
rgb(157pt)=(0.146616,0.673050,0.508936);
rgb(158pt)=(0.150148,0.676631,0.506589);
rgb(159pt)=(0.153894,0.680203,0.504172);
rgb(160pt)=(0.157851,0.683765,0.501686);
rgb(161pt)=(0.162016,0.687316,0.499129);
rgb(162pt)=(0.166383,0.690856,0.496502);
rgb(163pt)=(0.170948,0.694384,0.493803);
rgb(164pt)=(0.175707,0.697900,0.491033);
rgb(165pt)=(0.180653,0.701402,0.488189);
rgb(166pt)=(0.185783,0.704891,0.485273);
rgb(167pt)=(0.191090,0.708366,0.482284);
rgb(168pt)=(0.196571,0.711827,0.479221);
rgb(169pt)=(0.202219,0.715272,0.476084);
rgb(170pt)=(0.208030,0.718701,0.472873);
rgb(171pt)=(0.214000,0.722114,0.469588);
rgb(172pt)=(0.220124,0.725509,0.466226);
rgb(173pt)=(0.226397,0.728888,0.462789);
rgb(174pt)=(0.232815,0.732247,0.459277);
rgb(175pt)=(0.239374,0.735588,0.455688);
rgb(176pt)=(0.246070,0.738910,0.452024);
rgb(177pt)=(0.252899,0.742211,0.448284);
rgb(178pt)=(0.259857,0.745492,0.444467);
rgb(179pt)=(0.266941,0.748751,0.440573);
rgb(180pt)=(0.274149,0.751988,0.436601);
rgb(181pt)=(0.281477,0.755203,0.432552);
rgb(182pt)=(0.288921,0.758394,0.428426);
rgb(183pt)=(0.296479,0.761561,0.424223);
rgb(184pt)=(0.304148,0.764704,0.419943);
rgb(185pt)=(0.311925,0.767822,0.415586);
rgb(186pt)=(0.319809,0.770914,0.411152);
rgb(187pt)=(0.327796,0.773980,0.406640);
rgb(188pt)=(0.335885,0.777018,0.402049);
rgb(189pt)=(0.344074,0.780029,0.397381);
rgb(190pt)=(0.352360,0.783011,0.392636);
rgb(191pt)=(0.360741,0.785964,0.387814);
rgb(192pt)=(0.369214,0.788888,0.382914);
rgb(193pt)=(0.377779,0.791781,0.377939);
rgb(194pt)=(0.386433,0.794644,0.372886);
rgb(195pt)=(0.395174,0.797475,0.367757);
rgb(196pt)=(0.404001,0.800275,0.362552);
rgb(197pt)=(0.412913,0.803041,0.357269);
rgb(198pt)=(0.421908,0.805774,0.351910);
rgb(199pt)=(0.430983,0.808473,0.346476);
rgb(200pt)=(0.440137,0.811138,0.340967);
rgb(201pt)=(0.449368,0.813768,0.335384);
rgb(202pt)=(0.458674,0.816363,0.329727);
rgb(203pt)=(0.468053,0.818921,0.323998);
rgb(204pt)=(0.477504,0.821444,0.318195);
rgb(205pt)=(0.487026,0.823929,0.312321);
rgb(206pt)=(0.496615,0.826376,0.306377);
rgb(207pt)=(0.506271,0.828786,0.300362);
rgb(208pt)=(0.515992,0.831158,0.294279);
rgb(209pt)=(0.525776,0.833491,0.288127);
rgb(210pt)=(0.535621,0.835785,0.281908);
rgb(211pt)=(0.545524,0.838039,0.275626);
rgb(212pt)=(0.555484,0.840254,0.269281);
rgb(213pt)=(0.565498,0.842430,0.262877);
rgb(214pt)=(0.575563,0.844566,0.256415);
rgb(215pt)=(0.585678,0.846661,0.249897);
rgb(216pt)=(0.595839,0.848717,0.243329);
rgb(217pt)=(0.606045,0.850733,0.236712);
rgb(218pt)=(0.616293,0.852709,0.230052);
rgb(219pt)=(0.626579,0.854645,0.223353);
rgb(220pt)=(0.636902,0.856542,0.216620);
rgb(221pt)=(0.647257,0.858400,0.209861);
rgb(222pt)=(0.657642,0.860219,0.203082);
rgb(223pt)=(0.668054,0.861999,0.196293);
rgb(224pt)=(0.678489,0.863742,0.189503);
rgb(225pt)=(0.688944,0.865448,0.182725);
rgb(226pt)=(0.699415,0.867117,0.175971);
rgb(227pt)=(0.709898,0.868751,0.169257);
rgb(228pt)=(0.720391,0.870350,0.162603);
rgb(229pt)=(0.730889,0.871916,0.156029);
rgb(230pt)=(0.741388,0.873449,0.149561);
rgb(231pt)=(0.751884,0.874951,0.143228);
rgb(232pt)=(0.762373,0.876424,0.137064);
rgb(233pt)=(0.772852,0.877868,0.131109);
rgb(234pt)=(0.783315,0.879285,0.125405);
rgb(235pt)=(0.793760,0.880678,0.120005);
rgb(236pt)=(0.804182,0.882046,0.114965);
rgb(237pt)=(0.814576,0.883393,0.110347);
rgb(238pt)=(0.824940,0.884720,0.106217);
rgb(239pt)=(0.835270,0.886029,0.102646);
rgb(240pt)=(0.845561,0.887322,0.099702);
rgb(241pt)=(0.855810,0.888601,0.097452);
rgb(242pt)=(0.866013,0.889868,0.095953);
rgb(243pt)=(0.876168,0.891125,0.095250);
rgb(244pt)=(0.886271,0.892374,0.095374);
rgb(245pt)=(0.896320,0.893616,0.096335);
rgb(246pt)=(0.906311,0.894855,0.098125);
rgb(247pt)=(0.916242,0.896091,0.100717);
rgb(248pt)=(0.926106,0.897330,0.104071);
rgb(249pt)=(0.935904,0.898570,0.108131);
rgb(250pt)=(0.945636,0.899815,0.112838);
rgb(251pt)=(0.955300,0.901065,0.118128);
rgb(252pt)=(0.964894,0.902323,0.123941);
rgb(253pt)=(0.974417,0.903590,0.130215);
rgb(254pt)=(0.983868,0.904867,0.136897);
rgb(255pt)=(0.993248,0.906157,0.143936);
}
}
\pgfplotsset{
colormap={cmRdBu_r}{
rgb(000pt)=(0.019608,0.188235,0.380392);
rgb(001pt)=(0.023914,0.196540,0.391926);
rgb(002pt)=(0.028220,0.204844,0.403460);
rgb(003pt)=(0.032526,0.213149,0.414994);
rgb(004pt)=(0.036832,0.221453,0.426528);
rgb(005pt)=(0.041138,0.229758,0.438062);
rgb(006pt)=(0.045444,0.238062,0.449596);
rgb(007pt)=(0.049750,0.246367,0.461130);
rgb(008pt)=(0.054056,0.254671,0.472664);
rgb(009pt)=(0.058362,0.262976,0.484198);
rgb(010pt)=(0.062668,0.271280,0.495732);
rgb(011pt)=(0.066974,0.279585,0.507266);
rgb(012pt)=(0.071280,0.287889,0.518800);
rgb(013pt)=(0.075586,0.296194,0.530335);
rgb(014pt)=(0.079892,0.304498,0.541869);
rgb(015pt)=(0.084198,0.312803,0.553403);
rgb(016pt)=(0.088504,0.321107,0.564937);
rgb(017pt)=(0.092810,0.329412,0.576471);
rgb(018pt)=(0.097116,0.337716,0.588005);
rgb(019pt)=(0.101423,0.346021,0.599539);
rgb(020pt)=(0.105729,0.354325,0.611073);
rgb(021pt)=(0.110035,0.362630,0.622607);
rgb(022pt)=(0.114341,0.370934,0.634141);
rgb(023pt)=(0.118647,0.379239,0.645675);
rgb(024pt)=(0.122953,0.387543,0.657209);
rgb(025pt)=(0.127259,0.395848,0.668743);
rgb(026pt)=(0.132026,0.403460,0.676278);
rgb(027pt)=(0.137255,0.410381,0.679815);
rgb(028pt)=(0.142484,0.417301,0.683353);
rgb(029pt)=(0.147712,0.424221,0.686890);
rgb(030pt)=(0.152941,0.431142,0.690427);
rgb(031pt)=(0.158170,0.438062,0.693964);
rgb(032pt)=(0.163399,0.444983,0.697501);
rgb(033pt)=(0.168627,0.451903,0.701038);
rgb(034pt)=(0.173856,0.458824,0.704575);
rgb(035pt)=(0.179085,0.465744,0.708112);
rgb(036pt)=(0.184314,0.472664,0.711649);
rgb(037pt)=(0.189542,0.479585,0.715186);
rgb(038pt)=(0.194771,0.486505,0.718724);
rgb(039pt)=(0.200000,0.493426,0.722261);
rgb(040pt)=(0.205229,0.500346,0.725798);
rgb(041pt)=(0.210458,0.507266,0.729335);
rgb(042pt)=(0.215686,0.514187,0.732872);
rgb(043pt)=(0.220915,0.521107,0.736409);
rgb(044pt)=(0.226144,0.528028,0.739946);
rgb(045pt)=(0.231373,0.534948,0.743483);
rgb(046pt)=(0.236601,0.541869,0.747020);
rgb(047pt)=(0.241830,0.548789,0.750557);
rgb(048pt)=(0.247059,0.555709,0.754095);
rgb(049pt)=(0.252288,0.562630,0.757632);
rgb(050pt)=(0.257516,0.569550,0.761169);
rgb(051pt)=(0.262745,0.576471,0.764706);
rgb(052pt)=(0.274894,0.584160,0.768858);
rgb(053pt)=(0.287043,0.591849,0.773010);
rgb(054pt)=(0.299193,0.599539,0.777163);
rgb(055pt)=(0.311342,0.607228,0.781315);
rgb(056pt)=(0.323491,0.614917,0.785467);
rgb(057pt)=(0.335640,0.622607,0.789619);
rgb(058pt)=(0.347789,0.630296,0.793772);
rgb(059pt)=(0.359939,0.637985,0.797924);
rgb(060pt)=(0.372088,0.645675,0.802076);
rgb(061pt)=(0.384237,0.653364,0.806228);
rgb(062pt)=(0.396386,0.661053,0.810381);
rgb(063pt)=(0.408535,0.668743,0.814533);
rgb(064pt)=(0.420684,0.676432,0.818685);
rgb(065pt)=(0.432834,0.684122,0.822837);
rgb(066pt)=(0.444983,0.691811,0.826990);
rgb(067pt)=(0.457132,0.699500,0.831142);
rgb(068pt)=(0.469281,0.707190,0.835294);
rgb(069pt)=(0.481430,0.714879,0.839446);
rgb(070pt)=(0.493579,0.722568,0.843599);
rgb(071pt)=(0.505729,0.730258,0.847751);
rgb(072pt)=(0.517878,0.737947,0.851903);
rgb(073pt)=(0.530027,0.745636,0.856055);
rgb(074pt)=(0.542176,0.753326,0.860208);
rgb(075pt)=(0.554325,0.761015,0.864360);
rgb(076pt)=(0.566474,0.768704,0.868512);
rgb(077pt)=(0.577393,0.775010,0.871972);
rgb(078pt)=(0.587082,0.779931,0.874740);
rgb(079pt)=(0.596770,0.784852,0.877509);
rgb(080pt)=(0.606459,0.789773,0.880277);
rgb(081pt)=(0.616148,0.794694,0.883045);
rgb(082pt)=(0.625836,0.799616,0.885813);
rgb(083pt)=(0.635525,0.804537,0.888581);
rgb(084pt)=(0.645213,0.809458,0.891349);
rgb(085pt)=(0.654902,0.814379,0.894118);
rgb(086pt)=(0.664591,0.819300,0.896886);
rgb(087pt)=(0.674279,0.824221,0.899654);
rgb(088pt)=(0.683968,0.829143,0.902422);
rgb(089pt)=(0.693656,0.834064,0.905190);
rgb(090pt)=(0.703345,0.838985,0.907958);
rgb(091pt)=(0.713033,0.843906,0.910727);
rgb(092pt)=(0.722722,0.848827,0.913495);
rgb(093pt)=(0.732411,0.853749,0.916263);
rgb(094pt)=(0.742099,0.858670,0.919031);
rgb(095pt)=(0.751788,0.863591,0.921799);
rgb(096pt)=(0.761476,0.868512,0.924567);
rgb(097pt)=(0.771165,0.873433,0.927336);
rgb(098pt)=(0.780854,0.878354,0.930104);
rgb(099pt)=(0.790542,0.883276,0.932872);
rgb(100pt)=(0.800231,0.888197,0.935640);
rgb(101pt)=(0.809919,0.893118,0.938408);
rgb(102pt)=(0.819608,0.898039,0.941176);
rgb(103pt)=(0.825452,0.900807,0.942253);
rgb(104pt)=(0.831296,0.903576,0.943329);
rgb(105pt)=(0.837140,0.906344,0.944406);
rgb(106pt)=(0.842983,0.909112,0.945483);
rgb(107pt)=(0.848827,0.911880,0.946559);
rgb(108pt)=(0.854671,0.914648,0.947636);
rgb(109pt)=(0.860515,0.917416,0.948712);
rgb(110pt)=(0.866359,0.920185,0.949789);
rgb(111pt)=(0.872203,0.922953,0.950865);
rgb(112pt)=(0.878047,0.925721,0.951942);
rgb(113pt)=(0.883891,0.928489,0.953018);
rgb(114pt)=(0.889735,0.931257,0.954095);
rgb(115pt)=(0.895579,0.934025,0.955171);
rgb(116pt)=(0.901423,0.936794,0.956248);
rgb(117pt)=(0.907266,0.939562,0.957324);
rgb(118pt)=(0.913110,0.942330,0.958401);
rgb(119pt)=(0.918954,0.945098,0.959477);
rgb(120pt)=(0.924798,0.947866,0.960554);
rgb(121pt)=(0.930642,0.950634,0.961630);
rgb(122pt)=(0.936486,0.953403,0.962707);
rgb(123pt)=(0.942330,0.956171,0.963783);
rgb(124pt)=(0.948174,0.958939,0.964860);
rgb(125pt)=(0.954018,0.961707,0.965936);
rgb(126pt)=(0.959862,0.964475,0.967013);
rgb(127pt)=(0.965705,0.967243,0.968089);
rgb(128pt)=(0.969089,0.966474,0.964937);
rgb(129pt)=(0.970012,0.962168,0.957555);
rgb(130pt)=(0.970934,0.957862,0.950173);
rgb(131pt)=(0.971857,0.953556,0.942791);
rgb(132pt)=(0.972780,0.949250,0.935409);
rgb(133pt)=(0.973702,0.944944,0.928028);
rgb(134pt)=(0.974625,0.940638,0.920646);
rgb(135pt)=(0.975548,0.936332,0.913264);
rgb(136pt)=(0.976471,0.932026,0.905882);
rgb(137pt)=(0.977393,0.927720,0.898501);
rgb(138pt)=(0.978316,0.923414,0.891119);
rgb(139pt)=(0.979239,0.919108,0.883737);
rgb(140pt)=(0.980161,0.914802,0.876355);
rgb(141pt)=(0.981084,0.910496,0.868973);
rgb(142pt)=(0.982007,0.906190,0.861592);
rgb(143pt)=(0.982930,0.901884,0.854210);
rgb(144pt)=(0.983852,0.897578,0.846828);
rgb(145pt)=(0.984775,0.893272,0.839446);
rgb(146pt)=(0.985698,0.888966,0.832065);
rgb(147pt)=(0.986621,0.884660,0.824683);
rgb(148pt)=(0.987543,0.880354,0.817301);
rgb(149pt)=(0.988466,0.876048,0.809919);
rgb(150pt)=(0.989389,0.871742,0.802537);
rgb(151pt)=(0.990311,0.867436,0.795156);
rgb(152pt)=(0.991234,0.863130,0.787774);
rgb(153pt)=(0.992157,0.858824,0.780392);
rgb(154pt)=(0.990773,0.850519,0.769781);
rgb(155pt)=(0.989389,0.842215,0.759170);
rgb(156pt)=(0.988005,0.833910,0.748558);
rgb(157pt)=(0.986621,0.825606,0.737947);
rgb(158pt)=(0.985236,0.817301,0.727336);
rgb(159pt)=(0.983852,0.808997,0.716724);
rgb(160pt)=(0.982468,0.800692,0.706113);
rgb(161pt)=(0.981084,0.792388,0.695502);
rgb(162pt)=(0.979700,0.784083,0.684890);
rgb(163pt)=(0.978316,0.775779,0.674279);
rgb(164pt)=(0.976932,0.767474,0.663668);
rgb(165pt)=(0.975548,0.759170,0.653057);
rgb(166pt)=(0.974164,0.750865,0.642445);
rgb(167pt)=(0.972780,0.742561,0.631834);
rgb(168pt)=(0.971396,0.734256,0.621223);
rgb(169pt)=(0.970012,0.725952,0.610611);
rgb(170pt)=(0.968627,0.717647,0.600000);
rgb(171pt)=(0.967243,0.709343,0.589389);
rgb(172pt)=(0.965859,0.701038,0.578777);
rgb(173pt)=(0.964475,0.692734,0.568166);
rgb(174pt)=(0.963091,0.684429,0.557555);
rgb(175pt)=(0.961707,0.676125,0.546944);
rgb(176pt)=(0.960323,0.667820,0.536332);
rgb(177pt)=(0.958939,0.659516,0.525721);
rgb(178pt)=(0.957555,0.651211,0.515110);
rgb(179pt)=(0.954556,0.641753,0.505729);
rgb(180pt)=(0.949942,0.631142,0.497578);
rgb(181pt)=(0.945329,0.620531,0.489427);
rgb(182pt)=(0.940715,0.609919,0.481276);
rgb(183pt)=(0.936102,0.599308,0.473126);
rgb(184pt)=(0.931488,0.588697,0.464975);
rgb(185pt)=(0.926874,0.578085,0.456824);
rgb(186pt)=(0.922261,0.567474,0.448674);
rgb(187pt)=(0.917647,0.556863,0.440523);
rgb(188pt)=(0.913033,0.546251,0.432372);
rgb(189pt)=(0.908420,0.535640,0.424221);
rgb(190pt)=(0.903806,0.525029,0.416071);
rgb(191pt)=(0.899193,0.514418,0.407920);
rgb(192pt)=(0.894579,0.503806,0.399769);
rgb(193pt)=(0.889965,0.493195,0.391619);
rgb(194pt)=(0.885352,0.482584,0.383468);
rgb(195pt)=(0.880738,0.471972,0.375317);
rgb(196pt)=(0.876125,0.461361,0.367166);
rgb(197pt)=(0.871511,0.450750,0.359016);
rgb(198pt)=(0.866897,0.440138,0.350865);
rgb(199pt)=(0.862284,0.429527,0.342714);
rgb(200pt)=(0.857670,0.418916,0.334564);
rgb(201pt)=(0.853057,0.408305,0.326413);
rgb(202pt)=(0.848443,0.397693,0.318262);
rgb(203pt)=(0.843829,0.387082,0.310112);
rgb(204pt)=(0.839216,0.376471,0.301961);
rgb(205pt)=(0.833679,0.365398,0.296732);
rgb(206pt)=(0.828143,0.354325,0.291503);
rgb(207pt)=(0.822607,0.343253,0.286275);
rgb(208pt)=(0.817070,0.332180,0.281046);
rgb(209pt)=(0.811534,0.321107,0.275817);
rgb(210pt)=(0.805998,0.310035,0.270588);
rgb(211pt)=(0.800461,0.298962,0.265359);
rgb(212pt)=(0.794925,0.287889,0.260131);
rgb(213pt)=(0.789389,0.276817,0.254902);
rgb(214pt)=(0.783852,0.265744,0.249673);
rgb(215pt)=(0.778316,0.254671,0.244444);
rgb(216pt)=(0.772780,0.243599,0.239216);
rgb(217pt)=(0.767243,0.232526,0.233987);
rgb(218pt)=(0.761707,0.221453,0.228758);
rgb(219pt)=(0.756171,0.210381,0.223529);
rgb(220pt)=(0.750634,0.199308,0.218301);
rgb(221pt)=(0.745098,0.188235,0.213072);
rgb(222pt)=(0.739562,0.177163,0.207843);
rgb(223pt)=(0.734025,0.166090,0.202614);
rgb(224pt)=(0.728489,0.155017,0.197386);
rgb(225pt)=(0.722953,0.143945,0.192157);
rgb(226pt)=(0.717416,0.132872,0.186928);
rgb(227pt)=(0.711880,0.121799,0.181699);
rgb(228pt)=(0.706344,0.110727,0.176471);
rgb(229pt)=(0.700807,0.099654,0.171242);
rgb(230pt)=(0.692272,0.092272,0.167705);
rgb(231pt)=(0.680738,0.088581,0.165859);
rgb(232pt)=(0.669204,0.084890,0.164014);
rgb(233pt)=(0.657670,0.081200,0.162168);
rgb(234pt)=(0.646136,0.077509,0.160323);
rgb(235pt)=(0.634602,0.073818,0.158478);
rgb(236pt)=(0.623068,0.070127,0.156632);
rgb(237pt)=(0.611534,0.066436,0.154787);
rgb(238pt)=(0.600000,0.062745,0.152941);
rgb(239pt)=(0.588466,0.059054,0.151096);
rgb(240pt)=(0.576932,0.055363,0.149250);
rgb(241pt)=(0.565398,0.051672,0.147405);
rgb(242pt)=(0.553864,0.047982,0.145559);
rgb(243pt)=(0.542330,0.044291,0.143714);
rgb(244pt)=(0.530796,0.040600,0.141869);
rgb(245pt)=(0.519262,0.036909,0.140023);
rgb(246pt)=(0.507728,0.033218,0.138178);
rgb(247pt)=(0.496194,0.029527,0.136332);
rgb(248pt)=(0.484660,0.025836,0.134487);
rgb(249pt)=(0.473126,0.022145,0.132641);
rgb(250pt)=(0.461592,0.018454,0.130796);
rgb(251pt)=(0.450058,0.014764,0.128950);
rgb(252pt)=(0.438524,0.011073,0.127105);
rgb(253pt)=(0.426990,0.007382,0.125260);
rgb(254pt)=(0.415456,0.003691,0.123414);
rgb(255pt)=(0.403922,0.000000,0.121569);
}
}

\pgfplotsset{
colormap={cmRdBu_r_plus}{
rgb(000pt)=(0.969089,0.966474,0.964937);
rgb(001pt)=(0.969089,0.966474,0.964937);
rgb(002pt)=(0.970012,0.962168,0.957555);
rgb(003pt)=(0.970012,0.962168,0.957555);
rgb(004pt)=(0.970934,0.957862,0.950173);
rgb(005pt)=(0.970934,0.957862,0.950173);
rgb(006pt)=(0.971857,0.953556,0.942791);
rgb(007pt)=(0.971857,0.953556,0.942791);
rgb(008pt)=(0.972780,0.949250,0.935409);
rgb(009pt)=(0.972780,0.949250,0.935409);
rgb(010pt)=(0.973702,0.944944,0.928028);
rgb(011pt)=(0.973702,0.944944,0.928028);
rgb(012pt)=(0.974625,0.940638,0.920646);
rgb(013pt)=(0.974625,0.940638,0.920646);
rgb(014pt)=(0.975548,0.936332,0.913264);
rgb(015pt)=(0.975548,0.936332,0.913264);
rgb(016pt)=(0.976471,0.932026,0.905882);
rgb(017pt)=(0.976471,0.932026,0.905882);
rgb(018pt)=(0.977393,0.927720,0.898501);
rgb(019pt)=(0.977393,0.927720,0.898501);
rgb(020pt)=(0.978316,0.923414,0.891119);
rgb(021pt)=(0.978316,0.923414,0.891119);
rgb(022pt)=(0.979239,0.919108,0.883737);
rgb(023pt)=(0.979239,0.919108,0.883737);
rgb(024pt)=(0.980161,0.914802,0.876355);
rgb(025pt)=(0.980161,0.914802,0.876355);
rgb(026pt)=(0.981084,0.910496,0.868973);
rgb(027pt)=(0.981084,0.910496,0.868973);
rgb(028pt)=(0.982007,0.906190,0.861592);
rgb(029pt)=(0.982007,0.906190,0.861592);
rgb(030pt)=(0.982930,0.901884,0.854210);
rgb(031pt)=(0.982930,0.901884,0.854210);
rgb(032pt)=(0.983852,0.897578,0.846828);
rgb(033pt)=(0.983852,0.897578,0.846828);
rgb(034pt)=(0.984775,0.893272,0.839446);
rgb(035pt)=(0.984775,0.893272,0.839446);
rgb(036pt)=(0.985698,0.888966,0.832065);
rgb(037pt)=(0.985698,0.888966,0.832065);
rgb(038pt)=(0.986621,0.884660,0.824683);
rgb(039pt)=(0.986621,0.884660,0.824683);
rgb(040pt)=(0.987543,0.880354,0.817301);
rgb(041pt)=(0.987543,0.880354,0.817301);
rgb(042pt)=(0.988466,0.876048,0.809919);
rgb(043pt)=(0.988466,0.876048,0.809919);
rgb(044pt)=(0.989389,0.871742,0.802537);
rgb(045pt)=(0.989389,0.871742,0.802537);
rgb(046pt)=(0.990311,0.867436,0.795156);
rgb(047pt)=(0.990311,0.867436,0.795156);
rgb(048pt)=(0.991234,0.863130,0.787774);
rgb(049pt)=(0.991234,0.863130,0.787774);
rgb(050pt)=(0.992157,0.858824,0.780392);
rgb(051pt)=(0.992157,0.858824,0.780392);
rgb(052pt)=(0.990773,0.850519,0.769781);
rgb(053pt)=(0.990773,0.850519,0.769781);
rgb(054pt)=(0.989389,0.842215,0.759170);
rgb(055pt)=(0.989389,0.842215,0.759170);
rgb(056pt)=(0.988005,0.833910,0.748558);
rgb(057pt)=(0.988005,0.833910,0.748558);
rgb(058pt)=(0.986621,0.825606,0.737947);
rgb(059pt)=(0.986621,0.825606,0.737947);
rgb(060pt)=(0.985236,0.817301,0.727336);
rgb(061pt)=(0.985236,0.817301,0.727336);
rgb(062pt)=(0.983852,0.808997,0.716724);
rgb(063pt)=(0.983852,0.808997,0.716724);
rgb(064pt)=(0.982468,0.800692,0.706113);
rgb(065pt)=(0.982468,0.800692,0.706113);
rgb(066pt)=(0.981084,0.792388,0.695502);
rgb(067pt)=(0.981084,0.792388,0.695502);
rgb(068pt)=(0.979700,0.784083,0.684890);
rgb(069pt)=(0.979700,0.784083,0.684890);
rgb(070pt)=(0.978316,0.775779,0.674279);
rgb(071pt)=(0.978316,0.775779,0.674279);
rgb(072pt)=(0.976932,0.767474,0.663668);
rgb(073pt)=(0.976932,0.767474,0.663668);
rgb(074pt)=(0.975548,0.759170,0.653057);
rgb(075pt)=(0.975548,0.759170,0.653057);
rgb(076pt)=(0.974164,0.750865,0.642445);
rgb(077pt)=(0.974164,0.750865,0.642445);
rgb(078pt)=(0.972780,0.742561,0.631834);
rgb(079pt)=(0.972780,0.742561,0.631834);
rgb(080pt)=(0.971396,0.734256,0.621223);
rgb(081pt)=(0.971396,0.734256,0.621223);
rgb(082pt)=(0.970012,0.725952,0.610611);
rgb(083pt)=(0.970012,0.725952,0.610611);
rgb(084pt)=(0.968627,0.717647,0.600000);
rgb(085pt)=(0.968627,0.717647,0.600000);
rgb(086pt)=(0.967243,0.709343,0.589389);
rgb(087pt)=(0.967243,0.709343,0.589389);
rgb(088pt)=(0.965859,0.701038,0.578777);
rgb(089pt)=(0.965859,0.701038,0.578777);
rgb(090pt)=(0.964475,0.692734,0.568166);
rgb(091pt)=(0.964475,0.692734,0.568166);
rgb(092pt)=(0.963091,0.684429,0.557555);
rgb(093pt)=(0.963091,0.684429,0.557555);
rgb(094pt)=(0.961707,0.676125,0.546944);
rgb(095pt)=(0.961707,0.676125,0.546944);
rgb(096pt)=(0.960323,0.667820,0.536332);
rgb(097pt)=(0.960323,0.667820,0.536332);
rgb(098pt)=(0.958939,0.659516,0.525721);
rgb(099pt)=(0.958939,0.659516,0.525721);
rgb(100pt)=(0.957555,0.651211,0.515110);
rgb(101pt)=(0.957555,0.651211,0.515110);
rgb(102pt)=(0.954556,0.641753,0.505729);
rgb(103pt)=(0.954556,0.641753,0.505729);
rgb(104pt)=(0.949942,0.631142,0.497578);
rgb(105pt)=(0.949942,0.631142,0.497578);
rgb(106pt)=(0.945329,0.620531,0.489427);
rgb(107pt)=(0.945329,0.620531,0.489427);
rgb(108pt)=(0.940715,0.609919,0.481276);
rgb(109pt)=(0.940715,0.609919,0.481276);
rgb(110pt)=(0.936102,0.599308,0.473126);
rgb(111pt)=(0.936102,0.599308,0.473126);
rgb(112pt)=(0.931488,0.588697,0.464975);
rgb(113pt)=(0.931488,0.588697,0.464975);
rgb(114pt)=(0.926874,0.578085,0.456824);
rgb(115pt)=(0.926874,0.578085,0.456824);
rgb(116pt)=(0.922261,0.567474,0.448674);
rgb(117pt)=(0.922261,0.567474,0.448674);
rgb(118pt)=(0.917647,0.556863,0.440523);
rgb(119pt)=(0.917647,0.556863,0.440523);
rgb(120pt)=(0.913033,0.546251,0.432372);
rgb(121pt)=(0.913033,0.546251,0.432372);
rgb(122pt)=(0.908420,0.535640,0.424221);
rgb(123pt)=(0.908420,0.535640,0.424221);
rgb(124pt)=(0.903806,0.525029,0.416071);
rgb(125pt)=(0.903806,0.525029,0.416071);
rgb(126pt)=(0.899193,0.514418,0.407920);
rgb(127pt)=(0.899193,0.514418,0.407920);
rgb(128pt)=(0.894579,0.503806,0.399769);
rgb(129pt)=(0.894579,0.503806,0.399769);
rgb(130pt)=(0.889965,0.493195,0.391619);
rgb(131pt)=(0.889965,0.493195,0.391619);
rgb(132pt)=(0.885352,0.482584,0.383468);
rgb(133pt)=(0.885352,0.482584,0.383468);
rgb(134pt)=(0.880738,0.471972,0.375317);
rgb(135pt)=(0.880738,0.471972,0.375317);
rgb(136pt)=(0.876125,0.461361,0.367166);
rgb(137pt)=(0.876125,0.461361,0.367166);
rgb(138pt)=(0.871511,0.450750,0.359016);
rgb(139pt)=(0.871511,0.450750,0.359016);
rgb(140pt)=(0.866897,0.440138,0.350865);
rgb(141pt)=(0.866897,0.440138,0.350865);
rgb(142pt)=(0.862284,0.429527,0.342714);
rgb(143pt)=(0.862284,0.429527,0.342714);
rgb(144pt)=(0.857670,0.418916,0.334564);
rgb(145pt)=(0.857670,0.418916,0.334564);
rgb(146pt)=(0.853057,0.408305,0.326413);
rgb(147pt)=(0.853057,0.408305,0.326413);
rgb(148pt)=(0.848443,0.397693,0.318262);
rgb(149pt)=(0.848443,0.397693,0.318262);
rgb(150pt)=(0.843829,0.387082,0.310112);
rgb(151pt)=(0.843829,0.387082,0.310112);
rgb(152pt)=(0.839216,0.376471,0.301961);
rgb(153pt)=(0.839216,0.376471,0.301961);
rgb(154pt)=(0.833679,0.365398,0.296732);
rgb(155pt)=(0.833679,0.365398,0.296732);
rgb(156pt)=(0.828143,0.354325,0.291503);
rgb(157pt)=(0.828143,0.354325,0.291503);
rgb(158pt)=(0.822607,0.343253,0.286275);
rgb(159pt)=(0.822607,0.343253,0.286275);
rgb(160pt)=(0.817070,0.332180,0.281046);
rgb(161pt)=(0.817070,0.332180,0.281046);
rgb(162pt)=(0.811534,0.321107,0.275817);
rgb(163pt)=(0.811534,0.321107,0.275817);
rgb(164pt)=(0.805998,0.310035,0.270588);
rgb(165pt)=(0.805998,0.310035,0.270588);
rgb(166pt)=(0.800461,0.298962,0.265359);
rgb(167pt)=(0.800461,0.298962,0.265359);
rgb(168pt)=(0.794925,0.287889,0.260131);
rgb(169pt)=(0.794925,0.287889,0.260131);
rgb(170pt)=(0.789389,0.276817,0.254902);
rgb(171pt)=(0.789389,0.276817,0.254902);
rgb(172pt)=(0.783852,0.265744,0.249673);
rgb(173pt)=(0.783852,0.265744,0.249673);
rgb(174pt)=(0.778316,0.254671,0.244444);
rgb(175pt)=(0.778316,0.254671,0.244444);
rgb(176pt)=(0.772780,0.243599,0.239216);
rgb(177pt)=(0.772780,0.243599,0.239216);
rgb(178pt)=(0.767243,0.232526,0.233987);
rgb(179pt)=(0.767243,0.232526,0.233987);
rgb(180pt)=(0.761707,0.221453,0.228758);
rgb(181pt)=(0.761707,0.221453,0.228758);
rgb(182pt)=(0.756171,0.210381,0.223529);
rgb(183pt)=(0.756171,0.210381,0.223529);
rgb(184pt)=(0.750634,0.199308,0.218301);
rgb(185pt)=(0.750634,0.199308,0.218301);
rgb(186pt)=(0.745098,0.188235,0.213072);
rgb(187pt)=(0.745098,0.188235,0.213072);
rgb(188pt)=(0.739562,0.177163,0.207843);
rgb(189pt)=(0.739562,0.177163,0.207843);
rgb(190pt)=(0.734025,0.166090,0.202614);
rgb(191pt)=(0.734025,0.166090,0.202614);
rgb(192pt)=(0.728489,0.155017,0.197386);
rgb(193pt)=(0.728489,0.155017,0.197386);
rgb(194pt)=(0.722953,0.143945,0.192157);
rgb(195pt)=(0.722953,0.143945,0.192157);
rgb(196pt)=(0.717416,0.132872,0.186928);
rgb(197pt)=(0.717416,0.132872,0.186928);
rgb(198pt)=(0.711880,0.121799,0.181699);
rgb(199pt)=(0.711880,0.121799,0.181699);
rgb(200pt)=(0.706344,0.110727,0.176471);
rgb(201pt)=(0.706344,0.110727,0.176471);
rgb(202pt)=(0.700807,0.099654,0.171242);
rgb(203pt)=(0.700807,0.099654,0.171242);
rgb(204pt)=(0.692272,0.092272,0.167705);
rgb(205pt)=(0.692272,0.092272,0.167705);
rgb(206pt)=(0.680738,0.088581,0.165859);
rgb(207pt)=(0.680738,0.088581,0.165859);
rgb(208pt)=(0.669204,0.084890,0.164014);
rgb(209pt)=(0.669204,0.084890,0.164014);
rgb(210pt)=(0.657670,0.081200,0.162168);
rgb(211pt)=(0.657670,0.081200,0.162168);
rgb(212pt)=(0.646136,0.077509,0.160323);
rgb(213pt)=(0.646136,0.077509,0.160323);
rgb(214pt)=(0.634602,0.073818,0.158478);
rgb(215pt)=(0.634602,0.073818,0.158478);
rgb(216pt)=(0.623068,0.070127,0.156632);
rgb(217pt)=(0.623068,0.070127,0.156632);
rgb(218pt)=(0.611534,0.066436,0.154787);
rgb(219pt)=(0.611534,0.066436,0.154787);
rgb(220pt)=(0.600000,0.062745,0.152941);
rgb(221pt)=(0.600000,0.062745,0.152941);
rgb(222pt)=(0.588466,0.059054,0.151096);
rgb(223pt)=(0.588466,0.059054,0.151096);
rgb(224pt)=(0.576932,0.055363,0.149250);
rgb(225pt)=(0.576932,0.055363,0.149250);
rgb(226pt)=(0.565398,0.051672,0.147405);
rgb(227pt)=(0.565398,0.051672,0.147405);
rgb(228pt)=(0.553864,0.047982,0.145559);
rgb(229pt)=(0.553864,0.047982,0.145559);
rgb(230pt)=(0.542330,0.044291,0.143714);
rgb(231pt)=(0.542330,0.044291,0.143714);
rgb(232pt)=(0.530796,0.040600,0.141869);
rgb(233pt)=(0.530796,0.040600,0.141869);
rgb(234pt)=(0.519262,0.036909,0.140023);
rgb(235pt)=(0.519262,0.036909,0.140023);
rgb(236pt)=(0.507728,0.033218,0.138178);
rgb(237pt)=(0.507728,0.033218,0.138178);
rgb(238pt)=(0.496194,0.029527,0.136332);
rgb(239pt)=(0.496194,0.029527,0.136332);
rgb(240pt)=(0.484660,0.025836,0.134487);
rgb(241pt)=(0.484660,0.025836,0.134487);
rgb(242pt)=(0.473126,0.022145,0.132641);
rgb(243pt)=(0.473126,0.022145,0.132641);
rgb(244pt)=(0.461592,0.018454,0.130796);
rgb(245pt)=(0.461592,0.018454,0.130796);
rgb(246pt)=(0.450058,0.014764,0.128950);
rgb(247pt)=(0.450058,0.014764,0.128950);
rgb(248pt)=(0.438524,0.011073,0.127105);
rgb(249pt)=(0.438524,0.011073,0.127105);
rgb(250pt)=(0.426990,0.007382,0.125260);
rgb(251pt)=(0.426990,0.007382,0.125260);
rgb(252pt)=(0.415456,0.003691,0.123414);
rgb(253pt)=(0.415456,0.003691,0.123414);
rgb(254pt)=(0.403922,0.000000,0.121569);
rgb(255pt)=(0.403922,0.000000,0.121569);
}
}

\pgfplotsset{
colormap={cmRdBu_r_minus}{
rgb(000pt)=(0.019608,0.188235,0.380392);
rgb(001pt)=(0.019608,0.188235,0.380392);
rgb(002pt)=(0.023914,0.196540,0.391926);
rgb(003pt)=(0.023914,0.196540,0.391926);
rgb(004pt)=(0.028220,0.204844,0.403460);
rgb(005pt)=(0.028220,0.204844,0.403460);
rgb(006pt)=(0.032526,0.213149,0.414994);
rgb(007pt)=(0.032526,0.213149,0.414994);
rgb(008pt)=(0.036832,0.221453,0.426528);
rgb(009pt)=(0.036832,0.221453,0.426528);
rgb(010pt)=(0.041138,0.229758,0.438062);
rgb(011pt)=(0.041138,0.229758,0.438062);
rgb(012pt)=(0.045444,0.238062,0.449596);
rgb(013pt)=(0.045444,0.238062,0.449596);
rgb(014pt)=(0.049750,0.246367,0.461130);
rgb(015pt)=(0.049750,0.246367,0.461130);
rgb(016pt)=(0.054056,0.254671,0.472664);
rgb(017pt)=(0.054056,0.254671,0.472664);
rgb(018pt)=(0.058362,0.262976,0.484198);
rgb(019pt)=(0.058362,0.262976,0.484198);
rgb(020pt)=(0.062668,0.271280,0.495732);
rgb(021pt)=(0.062668,0.271280,0.495732);
rgb(022pt)=(0.066974,0.279585,0.507266);
rgb(023pt)=(0.066974,0.279585,0.507266);
rgb(024pt)=(0.071280,0.287889,0.518800);
rgb(025pt)=(0.071280,0.287889,0.518800);
rgb(026pt)=(0.075586,0.296194,0.530335);
rgb(027pt)=(0.075586,0.296194,0.530335);
rgb(028pt)=(0.079892,0.304498,0.541869);
rgb(029pt)=(0.079892,0.304498,0.541869);
rgb(030pt)=(0.084198,0.312803,0.553403);
rgb(031pt)=(0.084198,0.312803,0.553403);
rgb(032pt)=(0.088504,0.321107,0.564937);
rgb(033pt)=(0.088504,0.321107,0.564937);
rgb(034pt)=(0.092810,0.329412,0.576471);
rgb(035pt)=(0.092810,0.329412,0.576471);
rgb(036pt)=(0.097116,0.337716,0.588005);
rgb(037pt)=(0.097116,0.337716,0.588005);
rgb(038pt)=(0.101423,0.346021,0.599539);
rgb(039pt)=(0.101423,0.346021,0.599539);
rgb(040pt)=(0.105729,0.354325,0.611073);
rgb(041pt)=(0.105729,0.354325,0.611073);
rgb(042pt)=(0.110035,0.362630,0.622607);
rgb(043pt)=(0.110035,0.362630,0.622607);
rgb(044pt)=(0.114341,0.370934,0.634141);
rgb(045pt)=(0.114341,0.370934,0.634141);
rgb(046pt)=(0.118647,0.379239,0.645675);
rgb(047pt)=(0.118647,0.379239,0.645675);
rgb(048pt)=(0.122953,0.387543,0.657209);
rgb(049pt)=(0.122953,0.387543,0.657209);
rgb(050pt)=(0.127259,0.395848,0.668743);
rgb(051pt)=(0.127259,0.395848,0.668743);
rgb(052pt)=(0.132026,0.403460,0.676278);
rgb(053pt)=(0.132026,0.403460,0.676278);
rgb(054pt)=(0.137255,0.410381,0.679815);
rgb(055pt)=(0.137255,0.410381,0.679815);
rgb(056pt)=(0.142484,0.417301,0.683353);
rgb(057pt)=(0.142484,0.417301,0.683353);
rgb(058pt)=(0.147712,0.424221,0.686890);
rgb(059pt)=(0.147712,0.424221,0.686890);
rgb(060pt)=(0.152941,0.431142,0.690427);
rgb(061pt)=(0.152941,0.431142,0.690427);
rgb(062pt)=(0.158170,0.438062,0.693964);
rgb(063pt)=(0.158170,0.438062,0.693964);
rgb(064pt)=(0.163399,0.444983,0.697501);
rgb(065pt)=(0.163399,0.444983,0.697501);
rgb(066pt)=(0.168627,0.451903,0.701038);
rgb(067pt)=(0.168627,0.451903,0.701038);
rgb(068pt)=(0.173856,0.458824,0.704575);
rgb(069pt)=(0.173856,0.458824,0.704575);
rgb(070pt)=(0.179085,0.465744,0.708112);
rgb(071pt)=(0.179085,0.465744,0.708112);
rgb(072pt)=(0.184314,0.472664,0.711649);
rgb(073pt)=(0.184314,0.472664,0.711649);
rgb(074pt)=(0.189542,0.479585,0.715186);
rgb(075pt)=(0.189542,0.479585,0.715186);
rgb(076pt)=(0.194771,0.486505,0.718724);
rgb(077pt)=(0.194771,0.486505,0.718724);
rgb(078pt)=(0.200000,0.493426,0.722261);
rgb(079pt)=(0.200000,0.493426,0.722261);
rgb(080pt)=(0.205229,0.500346,0.725798);
rgb(081pt)=(0.205229,0.500346,0.725798);
rgb(082pt)=(0.210458,0.507266,0.729335);
rgb(083pt)=(0.210458,0.507266,0.729335);
rgb(084pt)=(0.215686,0.514187,0.732872);
rgb(085pt)=(0.215686,0.514187,0.732872);
rgb(086pt)=(0.220915,0.521107,0.736409);
rgb(087pt)=(0.220915,0.521107,0.736409);
rgb(088pt)=(0.226144,0.528028,0.739946);
rgb(089pt)=(0.226144,0.528028,0.739946);
rgb(090pt)=(0.231373,0.534948,0.743483);
rgb(091pt)=(0.231373,0.534948,0.743483);
rgb(092pt)=(0.236601,0.541869,0.747020);
rgb(093pt)=(0.236601,0.541869,0.747020);
rgb(094pt)=(0.241830,0.548789,0.750557);
rgb(095pt)=(0.241830,0.548789,0.750557);
rgb(096pt)=(0.247059,0.555709,0.754095);
rgb(097pt)=(0.247059,0.555709,0.754095);
rgb(098pt)=(0.252288,0.562630,0.757632);
rgb(099pt)=(0.252288,0.562630,0.757632);
rgb(100pt)=(0.257516,0.569550,0.761169);
rgb(101pt)=(0.257516,0.569550,0.761169);
rgb(102pt)=(0.262745,0.576471,0.764706);
rgb(103pt)=(0.262745,0.576471,0.764706);
rgb(104pt)=(0.274894,0.584160,0.768858);
rgb(105pt)=(0.274894,0.584160,0.768858);
rgb(106pt)=(0.287043,0.591849,0.773010);
rgb(107pt)=(0.287043,0.591849,0.773010);
rgb(108pt)=(0.299193,0.599539,0.777163);
rgb(109pt)=(0.299193,0.599539,0.777163);
rgb(110pt)=(0.311342,0.607228,0.781315);
rgb(111pt)=(0.311342,0.607228,0.781315);
rgb(112pt)=(0.323491,0.614917,0.785467);
rgb(113pt)=(0.323491,0.614917,0.785467);
rgb(114pt)=(0.335640,0.622607,0.789619);
rgb(115pt)=(0.335640,0.622607,0.789619);
rgb(116pt)=(0.347789,0.630296,0.793772);
rgb(117pt)=(0.347789,0.630296,0.793772);
rgb(118pt)=(0.359939,0.637985,0.797924);
rgb(119pt)=(0.359939,0.637985,0.797924);
rgb(120pt)=(0.372088,0.645675,0.802076);
rgb(121pt)=(0.372088,0.645675,0.802076);
rgb(122pt)=(0.384237,0.653364,0.806228);
rgb(123pt)=(0.384237,0.653364,0.806228);
rgb(124pt)=(0.396386,0.661053,0.810381);
rgb(125pt)=(0.396386,0.661053,0.810381);
rgb(126pt)=(0.408535,0.668743,0.814533);
rgb(127pt)=(0.408535,0.668743,0.814533);
rgb(128pt)=(0.420684,0.676432,0.818685);
rgb(129pt)=(0.420684,0.676432,0.818685);
rgb(130pt)=(0.432834,0.684122,0.822837);
rgb(131pt)=(0.432834,0.684122,0.822837);
rgb(132pt)=(0.444983,0.691811,0.826990);
rgb(133pt)=(0.444983,0.691811,0.826990);
rgb(134pt)=(0.457132,0.699500,0.831142);
rgb(135pt)=(0.457132,0.699500,0.831142);
rgb(136pt)=(0.469281,0.707190,0.835294);
rgb(137pt)=(0.469281,0.707190,0.835294);
rgb(138pt)=(0.481430,0.714879,0.839446);
rgb(139pt)=(0.481430,0.714879,0.839446);
rgb(140pt)=(0.493579,0.722568,0.843599);
rgb(141pt)=(0.493579,0.722568,0.843599);
rgb(142pt)=(0.505729,0.730258,0.847751);
rgb(143pt)=(0.505729,0.730258,0.847751);
rgb(144pt)=(0.517878,0.737947,0.851903);
rgb(145pt)=(0.517878,0.737947,0.851903);
rgb(146pt)=(0.530027,0.745636,0.856055);
rgb(147pt)=(0.530027,0.745636,0.856055);
rgb(148pt)=(0.542176,0.753326,0.860208);
rgb(149pt)=(0.542176,0.753326,0.860208);
rgb(150pt)=(0.554325,0.761015,0.864360);
rgb(151pt)=(0.554325,0.761015,0.864360);
rgb(152pt)=(0.566474,0.768704,0.868512);
rgb(153pt)=(0.566474,0.768704,0.868512);
rgb(154pt)=(0.577393,0.775010,0.871972);
rgb(155pt)=(0.577393,0.775010,0.871972);
rgb(156pt)=(0.587082,0.779931,0.874740);
rgb(157pt)=(0.587082,0.779931,0.874740);
rgb(158pt)=(0.596770,0.784852,0.877509);
rgb(159pt)=(0.596770,0.784852,0.877509);
rgb(160pt)=(0.606459,0.789773,0.880277);
rgb(161pt)=(0.606459,0.789773,0.880277);
rgb(162pt)=(0.616148,0.794694,0.883045);
rgb(163pt)=(0.616148,0.794694,0.883045);
rgb(164pt)=(0.625836,0.799616,0.885813);
rgb(165pt)=(0.625836,0.799616,0.885813);
rgb(166pt)=(0.635525,0.804537,0.888581);
rgb(167pt)=(0.635525,0.804537,0.888581);
rgb(168pt)=(0.645213,0.809458,0.891349);
rgb(169pt)=(0.645213,0.809458,0.891349);
rgb(170pt)=(0.654902,0.814379,0.894118);
rgb(171pt)=(0.654902,0.814379,0.894118);
rgb(172pt)=(0.664591,0.819300,0.896886);
rgb(173pt)=(0.664591,0.819300,0.896886);
rgb(174pt)=(0.674279,0.824221,0.899654);
rgb(175pt)=(0.674279,0.824221,0.899654);
rgb(176pt)=(0.683968,0.829143,0.902422);
rgb(177pt)=(0.683968,0.829143,0.902422);
rgb(178pt)=(0.693656,0.834064,0.905190);
rgb(179pt)=(0.693656,0.834064,0.905190);
rgb(180pt)=(0.703345,0.838985,0.907958);
rgb(181pt)=(0.703345,0.838985,0.907958);
rgb(182pt)=(0.713033,0.843906,0.910727);
rgb(183pt)=(0.713033,0.843906,0.910727);
rgb(184pt)=(0.722722,0.848827,0.913495);
rgb(185pt)=(0.722722,0.848827,0.913495);
rgb(186pt)=(0.732411,0.853749,0.916263);
rgb(187pt)=(0.732411,0.853749,0.916263);
rgb(188pt)=(0.742099,0.858670,0.919031);
rgb(189pt)=(0.742099,0.858670,0.919031);
rgb(190pt)=(0.751788,0.863591,0.921799);
rgb(191pt)=(0.751788,0.863591,0.921799);
rgb(192pt)=(0.761476,0.868512,0.924567);
rgb(193pt)=(0.761476,0.868512,0.924567);
rgb(194pt)=(0.771165,0.873433,0.927336);
rgb(195pt)=(0.771165,0.873433,0.927336);
rgb(196pt)=(0.780854,0.878354,0.930104);
rgb(197pt)=(0.780854,0.878354,0.930104);
rgb(198pt)=(0.790542,0.883276,0.932872);
rgb(199pt)=(0.790542,0.883276,0.932872);
rgb(200pt)=(0.800231,0.888197,0.935640);
rgb(201pt)=(0.800231,0.888197,0.935640);
rgb(202pt)=(0.809919,0.893118,0.938408);
rgb(203pt)=(0.809919,0.893118,0.938408);
rgb(204pt)=(0.819608,0.898039,0.941176);
rgb(205pt)=(0.819608,0.898039,0.941176);
rgb(206pt)=(0.825452,0.900807,0.942253);
rgb(207pt)=(0.825452,0.900807,0.942253);
rgb(208pt)=(0.831296,0.903576,0.943329);
rgb(209pt)=(0.831296,0.903576,0.943329);
rgb(210pt)=(0.837140,0.906344,0.944406);
rgb(211pt)=(0.837140,0.906344,0.944406);
rgb(212pt)=(0.842983,0.909112,0.945483);
rgb(213pt)=(0.842983,0.909112,0.945483);
rgb(214pt)=(0.848827,0.911880,0.946559);
rgb(215pt)=(0.848827,0.911880,0.946559);
rgb(216pt)=(0.854671,0.914648,0.947636);
rgb(217pt)=(0.854671,0.914648,0.947636);
rgb(218pt)=(0.860515,0.917416,0.948712);
rgb(219pt)=(0.860515,0.917416,0.948712);
rgb(220pt)=(0.866359,0.920185,0.949789);
rgb(221pt)=(0.866359,0.920185,0.949789);
rgb(222pt)=(0.872203,0.922953,0.950865);
rgb(223pt)=(0.872203,0.922953,0.950865);
rgb(224pt)=(0.878047,0.925721,0.951942);
rgb(225pt)=(0.878047,0.925721,0.951942);
rgb(226pt)=(0.883891,0.928489,0.953018);
rgb(227pt)=(0.883891,0.928489,0.953018);
rgb(228pt)=(0.889735,0.931257,0.954095);
rgb(229pt)=(0.889735,0.931257,0.954095);
rgb(230pt)=(0.895579,0.934025,0.955171);
rgb(231pt)=(0.895579,0.934025,0.955171);
rgb(232pt)=(0.901423,0.936794,0.956248);
rgb(233pt)=(0.901423,0.936794,0.956248);
rgb(234pt)=(0.907266,0.939562,0.957324);
rgb(235pt)=(0.907266,0.939562,0.957324);
rgb(236pt)=(0.913110,0.942330,0.958401);
rgb(237pt)=(0.913110,0.942330,0.958401);
rgb(238pt)=(0.918954,0.945098,0.959477);
rgb(239pt)=(0.918954,0.945098,0.959477);
rgb(240pt)=(0.924798,0.947866,0.960554);
rgb(241pt)=(0.924798,0.947866,0.960554);
rgb(242pt)=(0.930642,0.950634,0.961630);
rgb(243pt)=(0.930642,0.950634,0.961630);
rgb(244pt)=(0.936486,0.953403,0.962707);
rgb(245pt)=(0.936486,0.953403,0.962707);
rgb(246pt)=(0.942330,0.956171,0.963783);
rgb(247pt)=(0.942330,0.956171,0.963783);
rgb(248pt)=(0.948174,0.958939,0.964860);
rgb(249pt)=(0.948174,0.958939,0.964860);
rgb(250pt)=(0.954018,0.961707,0.965936);
rgb(251pt)=(0.954018,0.961707,0.965936);
rgb(252pt)=(0.959862,0.964475,0.967013);
rgb(253pt)=(0.959862,0.964475,0.967013);
rgb(254pt)=(0.965705,0.967243,0.968089);
rgb(255pt)=(0.969089,0.966474,0.964937);
}
}

\usepackage{amsmath}
\usepackage{amsfonts}
\usepackage{ifthen}
\newcommand{\bs}[1]{\boldsymbol{#1}}

\newcommand{\R}{\mathbb{R}}

\DeclareMathOperator{\trace}{tr}
\newcommand{\vecspan}[1]{\text{span} \lbrace #1 \rbrace}

\newcommand{\abs}[1]{\left| #1 \right|}
\newcommand{\norm}[1]{\ensuremath{\Vert #1 \Vert}}
\newcommand{\trans}[1][]{^{#1\top}}
\newcommand{\placeholder}{\cdot}

\DeclareMathOperator{\Eigenvalue}{EV}

\newcommand{\ltwospace}{L^2}
\newcommand{\weightedltwospace}{L^2_{\FD}}

\newcommand{\dirac}[2]{\delta_{#1}(#2)}

\newcommand{\de}{\partial}
\newcommand{\dt}{\de_t}

\newcommand{\divx}{\ensuremath{\nabla_x \cdot}}
\newcommand{\gradx}{\nabla_x}
\newcommand{\pfrac}[2]{\ensuremath{\frac{\partial #1}{\partial #2}}}

\newcommand{\xcoord}{x}
\newcommand{\xcx}{\xi}
\newcommand{\xcy}{\eta}
\newcommand{\xcz}{\zeta}
\newcommand{\tcoord}{t}
\newcommand{\vcoord}{v}
\newcommand{\nxcoord}{\normal{\xcoord}}
\newcommand{\ntcoord}{\normal{\tcoord}}
\newcommand{\nvcoord}{\normal{\vcoord}}
\newcommand{\speed}{c}
\newcommand{\xdomain}{\Omega_{\xcoord}}
\newcommand{\tdomain}{\Omega_{\tcoord}}
\newcommand{\vdomain}{\Omega_{\vcoord}}
\newcommand{\nxdomain}{\normal{\Omega}_{\xcoord}}
\newcommand{\ntdomain}{[0,1]}
\newcommand{\nvdomain}{\US}
\newcommand{\domain}{\Omega_{\tcoord\xcoord\vcoord}}

\newcommand{\xref}{X}
\newcommand{\tref}{T}
\newcommand{\vref}{c}

\newcommand{\polradius}{r}
\newcommand{\polangle}{\theta}
\newcommand{\quadrant}[1]{Q_{#1}}

\newcommand{\normal}[1]{\ensuremath{\hat{#1}}}
\newcommand{\US}{\mathbb{S}^2}

\newcommand{\intV}[1]{\int_{\vdomain} #1 d\vcoord}

\newcommand{\intVn}[1]{\int_{\nvdomain} #1 d\normal\vcoord}
\newcommand{\intVnprime}[1]{\int_{\nvdomain} #1 d\normal\vcoord'}
\newcommand{\ints}[1]{\ensuremath{\left\langle #1 \right\rangle}}

\newcommand{\outernormal}{n}
\newcommand{\boundarykernel}{B}
\newcommand{\Vplus}[1]{#1 \cdot \outernormal > 0}
\newcommand{\Vminus}[1]{#1 \cdot \outernormal < 0}
\newcommand{\intBplus}[2]{\int_{\Vplus{#1}} #2 d#1}
\newcommand{\intBminus}[2]{\int_{\Vminus{#1}} #2 d#1}
\DeclareMathOperator{\St}{St}
\DeclareMathOperator{\Kn}{Kn}
\newcommand{\Knt}{\Kn_t}
\newcommand{\Knp}{\Kn_p}
\newcommand{\pareps}{\varepsilon}
\newcommand{\pardel}{\delta}
\newcommand{\parnu}{\nu}
\newcommand{\parthet}{\theta}

\newcommand{\PD}{\ensuremath{f}}
\newcommand{\FD}{\ensuremath{E}}
\newcommand{\VF}{\ensuremath{Q}}
\newcommand{\PP}{\ensuremath{g}}
\newcommand{\PM}{\rho}
\newcommand{\density}{{\rho_0}}

\newcommand{\dif}{D}
\newcommand{\adv}{a}

\newcommand{\LCO}[1][]{\mathcal{L}_{#1}}
\newcommand{\nLCO}[1][]{\normal{\LCO}_{#1}}
\newcommand{\Source}{\mathcal{S}}
\newcommand{\nSource}{\normal{\Source}}

\newcommand{\Nsp}{\mathcal{N}}
\newcommand{\Nsporth}{\Nsp^\bot}
\newcommand{\Range}{\mathcal{R}}
\newcommand{\Identity}{I}
\newcommand{\Nspproj}{\Pi}
\newcommand{\Ptbproj}{(\Identity - \Nspproj)}

\newcommand{\factorcol}[1][]{\kappa_{#1}}
\newcommand{\nfactorcol}[1][]{\normal{\factorcol}_{#1}}
\newcommand{\factorprol}{\mu}
\newcommand{\nfactorprol}{\normal \mu}

\newcommand{\factordifref}{K_{\dif}}
\newcommand{\factoradvref}{K_{\adv}}
\newcommand{\factorprolref}{M}
\newcommand{\carryingcapacity}{\PM_{\text{cc}}}

\newcommand{\kernelcol}[1][]{k_{#1}}
\newcommand{\nkernelcol}[1][]{\normal{\kernelcol}_{#1}}

\newcommand{\Difftens}{D}
\newcommand{\DW}{D_W}
\newcommand{\DT}{D_T}
\newcommand{\drift}{a}
\newcommand{\driftT}{a_{T}}

\newcommand{\basisPD}{\bs{m}}
\newcommand{\basisPDord}[1]{\basisPD^{(#1)}}
\newcommand{\basisPDcomp}[1]{m_{#1}}
\newcommand{\basisPP}{\bs{a}}
\newcommand{\basisPPord}[1]{\basisPP^{(#1)}}
\newcommand{\basisPPcomp}[1]{a_{#1}}

\newcommand{\momPD}{\bs{w}}

\newcommand{\momPDcomp}[1]{w_{#1}}
\newcommand{\momPP}{\bs{u}}
\newcommand{\momPPord}[1]{\momPP^{(#1)}}
\newcommand{\momPPcomp}[1]{u_{#1}}

\newcommand{\recPD}{\mathfrak{f}}

\newcommand{\recPP}{\mathfrak{g}}

\newcommand{\mischeme}[2][]{\ensuremath{MM#2_{#1}}} 
\newcommand{\ivscheme}[2][]{\ensuremath{MM#2i_{#1}}} 

\newcommand{\orderof}[1]{\mathcal{O}\left(#1\right)}
\newcommand{\avg}[1]{\ensuremath{\left\lbrace #1 \right\rbrace}}
\newcommand{\spquad}{\mathcal{Q}}

\newcommand{\reldiff}[2]{\ensuremath{\Delta_{rel}\left(#1,#2\right)}}
\newcommand{\scinum}[2]{%
	\ifthenelse{\equal{#1}{1}}{%
		\ifthenelse{\equal{#2}{0}}{%
			1%
		}{%
			10^{#2}%
		}%
	}{%
		\ifthenelse{\equal{#2}{0}}{%
			#1%
		}{%
			#1 \times 10^{#2}%
		}%
	}%
}

\newcommand{\timestep}{\Delta t}
\newcommand{\timeidx}{n}
\newcommand{\timetwostepfraction}{\tau}
\newcommand{\timetwostepweight}{\sigma}
\newcommand{\gridsize}{\Delta x}
\newcommand{\evalat}[2]{\left. #1 \right|_{#2}}
\newcommand{\spacedim}{{S}}
\newcommand{\dimidx}{d}

\newcommand{\macroexpl}{\Phi^{\PM}}
\newcommand{\macroimpl}{\Gamma^{\PM}}
\newcommand{\microexpl}{\Phi^{\PP}}
\newcommand{\microexplfibres}{{\Phi^{\PP}_{\FD}}}
\newcommand{\microexplfibressk}[2]{\Phi^{\PP, (#2)}_{\FD, #1}}
\newcommand{\microimpl}{\Gamma^{\PP}}

\newcommand{\macroivexpl}{\tilde \Phi^{\PM}}
\newcommand{\macroivimpl}{\tilde \Gamma^{\PM}}
\newcommand{\microivexpl}{\tilde \Phi^{\PP}}
\newcommand{\microivexplfibres}{{\tilde \Phi^{\PP}_{\FD}}}

\newcommand{\microivimpl}{\tilde \Gamma^{\PP}}

\newcommand{\onehalf}{\ensuremath{\frac{1}{2}}}
\newcommand{\threehalf}{\ensuremath{\frac{3}{2}}}
\DeclareMathOperator{\cuboid}{Box}
\newcommand{\multiidx}{\bs{i}}
\newcommand{\xidx}{l}
\newcommand{\yidx}{m}
\newcommand{\pnodeidx}{(\xidx + \onehalf, \yidx + \onehalf)}
\newcommand{\dnodeidx}{(\xidx, \yidx)}
\newcommand{\pnodell}{(\xidx, \yidx)}
\newcommand{\pnodelr}{(\xidx + 1, \yidx)}
\newcommand{\pnodeul}{(\xidx, \yidx + 1)}
\newcommand{\pnodeur}{(\xidx + 1, \yidx + 1)}
\newcommand{\dnodell}{(\xidx - \onehalf, \yidx - \onehalf)}
\newcommand{\dnodelr}{(\xidx + \onehalf, \yidx - \onehalf)}
\newcommand{\dnodeul}{(\xidx - \onehalf, \yidx + \onehalf)}
\newcommand{\dnodeur}{(\xidx + \onehalf, \yidx + \onehalf)}

\newcommand{\mentity}[1]{\mathcal{E}_{#1}}
\newcommand{\midx}{i}
\newcommand{\mneighbors}[1]{\mathfrak{N}_{#1}}
\newcommand{\madjacents}[2]{\mathfrak{A}_{#1, #2}}
\newcommand{\mverticesof}[1]{\mathfrak{V}_{#1}}

\newcommand{\mavg}[2]{\avg{#2}_{#1}}
\newcommand{\mapproxavg}[2]{\avg{#2}_{#1}^{\spquad}}

\newcommand{\mcell}[1]{\Omega_{#1}}
\newcommand{\msubcell}[2]{\omega_{#1}^{#2}}
\newcommand{\mvert}[1]{\xcoord_{#1}}
\newcommand{\mface}[2]{\partial \Omega_{#1, #2}}
\newcommand{\mfacet}[3]{\partial \omega_{#1,#2}^{#3}}
\newcommand{\medge}[2]{\overline{\mvert{#1}\mvert{#2}}}

\newcommand{\isect}{\cap}
\newcommand{\union}{\cup}

\newcommand{\pcellidx}{j}
\newcommand{\pcellidxnb}{k}
\newcommand{\pvertidx}{r}
\newcommand{\pvertidxnb}{s}
\newcommand{\dcellidx}{\pvertidx}
\newcommand{\dcellidxnb}{\pvertidxnb}
\newcommand{\dvertidx}{\pcellidx}

\newcommand{\pcell}{\mcell{\pcellidx}}
\newcommand{\psubcell}{\msubcell{\pcellidx}{\pvertidx}}
\newcommand{\pcellnb}{\mcell{\pcellidxnb}}
\newcommand{\pvert}{\mvert{\pvertidx}}
\newcommand{\pvertnb}{\mvert{\pvertidxnb}}
\newcommand{\dcell}{\mcell{\dcellidx}}

\newcommand{\dvert}{\mvert{\dvertidx}}

\newcommand{\pface}{\mface{\pcellidx}{\pcellidxnb}}
\newcommand{\pfacet}{\mfacet{\pcellidx}{\pcellidxnb}{\dcellidx}}
\newcommand{\pedge}{\medge{\pvertidx}{\pvertidxnb}}
\newcommand{\dface}{\mface{\dcellidx}{\dcellidxnb}}
\newcommand{\dfacet}{\mfacet{\dcellidx}{\dcellidxnb}{\pcellidx}}
\newcommand{\dfacetone}{\mfacet{\dcellidx}{\dcellidxnb}{\pcellidx, 1}}
\newcommand{\dfacettwo}{\mfacet{\dcellidx}{\dcellidxnb}{\pcellidx, 2}}

\newcommand{\dPM}{\PM_{\dcellidx}}
\newcommand{\pPP}{\PP_{\pcellidx}}

\newcommand{\avgvec}[1]{\bar{#1}}
\newcommand{\PMvec}{\avgvec{\PM}}
\newcommand{\PPvec}{\avgvec{\PP}}
\newcommand{\macroexplvec}{\avgvec{\Phi}^{\PM}}
\newcommand{\macroimplvec}{\avgvec{\Gamma}^{\PM}}
\newcommand{\microexplvec}{\avgvec{\Phi}^{\PP}}
\newcommand{\microexplfibresvec}{\avgvec{\Phi}^{\PP}_{\FD}}
\newcommand{\microimplvec}{\avgvec{\Gamma}^{\PP}}

\newcommand{\fluxrg}[1][\dcellidx, \dcellidxnb]{F^{(\PM, \PP)}_{#1}}
\newcommand{\fluxgg}[1][\pcellidx, \pcellidxnb]{F^{(\PP, \PP)}_{#1}}
\newcommand{\fluxgr}[1][\pcellidx, \pcellidxnb]{F^{(\PP, \PM)}_{#1}}

\newcommand{\sprec}[1]{\bs{#1}}

\begin{document}

\begin{abstract}
	The so-called haptotaxis equation is a special class of transport equation that arises from models of biological cell movement along tissue fibers. 
	This equation has an anisotropic advection-diffusion equation as its macroscopic limit.  	
	An up to second-order accurate asymptotic preserving method is developed for the haptotaxis equation in space dimension up to three.
	For this the micro-macro decomposition proposed by Lemou and Mieussens is generalized in the context of finite-volume schemes on staggered grids. 
	The spurious modes that arise from this discretization can be eliminated by combining flux evaluations from different points in the right way. 
	The velocity space is discretized by an arbitrary-order linear moment system ($P_N$). 
\end{abstract}
\begin{keyword}
	Multiscale model, glioma invasion, kinetic transport equation, asymptotic preserving, macroscopic scaling, moment closure, reaction-diffusion-transport equations
\end{keyword}
\begin{AMS}
	35L04 
	35Q92 
	65M08 
	65M12 
	65M20 
	65M70 
	92C17 
	92D25 
\end{AMS}

\maketitle

\listoftodos

\section{Introduction}
\label{sec:intro}
The kinetic theory of active particles(KTAP) \cite{bellomo2009complexity} provides a framework to describe large systems of interacting living particles on multiple scales. 
Prominent examples of phenomena modeled in this setting include bacterial movement, cell migration, animal swarms and pedestrian crowds. 
Viewed at very small length and time scales, one can observe individual particles, each with its own complex internal dynamic and interactions with the environment or other particles. 
When many particles are involved, this level of detail is not practical. 
As a first level of abstraction, the KTAP theory models the microscopic scale with PDEs for the expected distribution of particles in time, physical space and state space; so-called kinetic equations. 

The connection between particle systems and kinetic equations has been established formally for example for neutron transport \cite{golse2012recent} and the movement of a bacterium \cite{stroock1974some}. 
However, in the context of the kinetic theory of active particles,  the models are formulated directly as a kinetic equation \cite{PH13, BBNS10a, bellomo2006onset, hillenM5}. 
Kinetic equations are characterized by a free-streaming transport term resulting from particles movement and a collision operator modeling particle interactions as instantaneous state changes. 
At larger scales only the resulting macroscopic population behavior can be observed, that is, the total number density of particles regardless of their internal microscopic state.
To pass from the microscopic description to a population law, one considers the limit of the kinetic equation when the mean free path of particles tends to zero. 
Analytically, passage to the limiting macroscopic equation has been extensively studied for neutron transport \cite{LarKel75, BSS84} and more recently also in the context of biological cell migration \cite{othmer2000diffusion, burini2017hilbert}.  
When only interactions between particles and the environment are considered and interactions between particles are neglected, the collision operator is linear. 
In this case the resulting macroscopic equations are of diffusion type \cite{burini2017hilbert}. 

A macroscopic equation derived in this manner can of course only be an approximation and one may ask how accurate it is in any given situation. 
From a computational standpoint this means that we would like to compare simulations of the microscopic and macroscopic models. 

However, when the mean free path is small, the collision term is very stiff and a straightforward discretization of the kinetic equation would need infeasible spatial and temporal resolution to resolve the small scales accurately \cite{larsen1987asymptotic}.
Therefore, a variety of so-called asymptotic preserving schemes have been developed \cite{jin1996numerical, klar1998asymptotic, jin2000uniformly, gosse2002asymptotic, lemou2008ap, buet2012design}.
These methods are constructed in such a way that---for a fixed resolution---they converge to a discretization of the limit equation.
A large portion of the work has been done in the context of the telegraph equation and the neutron transport equation, preferably in one space dimension. 
To obtain analytical insights about the method, for instance stability conditions or consistency errors, it is reasonable to simplify the situation as much as possible. 

But in this work the emphasis is on application rather than analysis. 
As a step towards adapting AP methods for more applied situations we consider a kinetic model for glioma invasion in the human brain, developed in \cite{EHKS14, EHS}. 
Malignant glioma are a type of brain tumor arising from mutations of glia cells. 
Tumor recurrence after treatment is very probable because glioma cells migrate far from the original tumor site without being detected by state-of-the-art imaging methods \cite{claes2007}. 
Predictive models could be used to estimate the invisible parts of the tumor and improve treatment success. 
The model takes haptic interactions between glioma cells and white matter tissue into account.  
According to the classification in \cite{dickinson1993stochastic}, this effect can be classified as either klinokinesis or taxis. 
In addition to an anisotropic diffusion, the resulting macroscopic model features a drift towards regions with denser fibers.

We develop an AP method against this prototype model, which introduces some extra real-world complications. 
In clinical praxis, information about the tissue structure of a patient's brain is contained in a diffusion tensor image (DTI) \cite{LeBihan2001DTI} obtained from a MRI scan. 
The three dimensional DTI data comes in the form of a constant tensor per voxel with a spatial resolution of a few millimeters. 
To avoid interpolation artifacts, the discretization should respect the data resolution. 
Also, the scheme must be robust against discontinuities in the data.  

Our scheme is an extension of the method of Lemou and Mieussens \cite{lemou2008ap} who employ a micro-macro decomposition on staggered grids. 

In the following \secref{sec:kinetic-equation}, we introduce the kinetic equation first in a general form and then in the specific form of the glioma invasion model. 
We also introduce a parabolic scaling of this equation. 
Then in \secref{sec:micro-macro}, we briefly introduce the micro-macro decomposition and use this to informally derive the macroscopic limit of the kinetic equation. 
A large part of the paper is dedicated to a detailed description of the AP method.
In \secref{sec:ap-scheme}, we first describe the space discretization on general primal-dual mesh pairs and then also present the scheme for the special situation of a regular grid.
We also discuss the resulting numerical scheme in the parabolic limit and how to overcome some of the problems of this limit scheme. 
Time stepping and boundary conditions will also be described briefly. 
It remains to find a suitable discretization of the velocity. 
The linear spectral method that we use is described in \secref{sec:spectral-method}. 
We do not do much analysis on the developed method but rather assess the method's properties numerically. 
Therefore we present the results of a number of benchmark tests in \secref{sec:results}. 
The emphasis is on situations close to the parabolic limit, also in the presence of discontinuous coefficients. 
Finally we perform a series of computations on the glioma model with measured DTI data and realistic parameters.

\section{Haptotaxis models and their diffusion limit}
\label{sec:kinetic-equation}
First, we recall the general class of kinetic equations from \cite{corbin2018higher}. 
Then we perform a parabolic scaling of this equation and present the resulting diffusion limit from \cite{EHKS14,corbin2018higher} without any derivation.
Finally we introduce a  model for glioma invasion as a special case of the general setting.
 
\subsection{General microscopic setting}
The population is described by a distribution function $\PD(\tcoord,\xcoord,\nvcoord)$ which can be interpreted as the number density of particles with speed $\nvcoord \in \US$ at time $\tcoord \in \R^{+}$ and position $\xcoord = (\xcx, \xcy, \xcz)$. 
The particle distribution is governed by a linear kinetic equation of the form
\begin{align}
	\dt \PD + \speed \divx (\nvcoord \PD) &= (\LCO[\dif] + \LCO[\adv]) \PD + \Source \PD,
	\label{eq:lke-general}
\end{align}
on the domain 
\begin{align*}
	\domain &= \tdomain \times \xdomain \times \vdomain\\
	&= \tref \ntdomain \times \xref \nxdomain \times \nvdomain.
\end{align*}
The left hand side models the free flight of particles with constant speed $\speed$ in arbitrary direction $\nvcoord \in \nvdomain$. 
Changes in velocity happen in so-called collisions, i.e. particles change their velocity instantaneously at certain times. 
This is modeled by the linear turning operator $(\LCO[\dif] + \LCO[\adv])$ on the right hand side of the equation.
Let $\kernelcol(\xcoord, \nvcoord, \nvcoord') := \kernelcol[\dif](\xcoord, \nvcoord, \nvcoord') + \kernelcol[\adv](\xcoord, \nvcoord, \nvcoord')$ be the rate at which particles at position $\xcoord$ with direction $\nvcoord'$ collide and change their direction to $\nvcoord$.
The interpretation as a rate is only meaningful if $\kernelcol$ is strictly positive and bounded from above: 
\begin{alignat}{3}
0 &< \kernelcol[min]  &&\leq \kernelcol[\dif](\xcoord,\nvcoord^\prime,\nvcoord) + \kernelcol[\adv] (\xcoord,\nvcoord^\prime,\nvcoord) &&\leq \kernelcol[max]. 
\label{eq:kernel-bounds2}
\end{alignat}
The turning operator $\LCO$ then maps the distribution $\PD$ onto another distribution $\LCO \PD$ via the kernel integral
\begin{align*}
	\LCO \PD = (\LCO[\dif] + \LCO[\adv]) \PD &= \intVnprime{ \kernelcol (\xcoord, \nvcoord, \nvcoord') \PD(\nvcoord') - \kernelcol(\xcoord,\nvcoord',\nvcoord) \PD(\nvcoord) } . 
\end{align*}
The first summand counts the gain for direction $\nvcoord$ due to particles turning from any direction $\nvcoord'$ to $\nvcoord$. 
Accordingly the second term describes the particle losses for direction $\nvcoord$.
By this construction the operator $\LCO$ (as well as both parts $\LCO[\dif], \LCO[\adv]$ individually) preserves mass: 
\begin{align}
	\intVn { (\LCO \PD)(\nvcoord) } = 0 . 
	\label{eq:lco-mass-conservation}
\end{align}

We need some additional structure for the turning to derive a diffusion limit. 
The first kernel $\kernelcol[\dif]$ on its own is a turning rate, i.e. positive and bounded from above: 
\begin{alignat}{3}
\label{eq:kernel-bounds1}
0 &< \kernelcol[\dif, min] &&\leq \kernelcol[\dif](\xcoord,\nvcoord^\prime,\nvcoord)  &&\leq \kernelcol[\dif, max]. 
\end{alignat}
There is a positive normalization factor 
\begin{align*}
\factorcol[\dif](\xcoord) := \intVn{\kernelcol[\dif](\xcoord, \nvcoord, \nvcoord') }
\end{align*}
that does not depend on the velocity $\vcoord'$.  
The kernel is strictly positive, normalized and first-order symmetric: 
\begin{equation}
\label{eq:equilibrium-assumptions} 
\begin{aligned}
\FD(\xcoord, \nvcoord) &> 0, \\
\intVn{\FD(\xcoord, \nvcoord)} &= 1, \\
\intVn{\nvcoord \FD (\xcoord, \nvcoord)} &= 0.
\end{aligned}
\end{equation}
Additionally it admits a local equilibrium $\FD(\xcoord,\nvcoord)$ that fulfills the detailed balance 
\begin{align}
\label{eq:kernel-equilibrium}
\kernelcol[\dif](\xcoord,\nvcoord,\nvcoord') \FD(\xcoord,\nvcoord') = \kernelcol[\dif](\xcoord,\nvcoord',\nvcoord) \FD(\xcoord,\nvcoord).
\end{align}
This is a slightly more general assumption than the symmetry assumption $\kernelcol(\nvcoord, \nvcoord') = \kernelcol(\nvcoord', \nvcoord)$ in classic linear kinetic theory. 

The kernel $\kernelcol[\adv]$ should be interpreted as a perturbation of the turning rate $\kernelcol[\dif]$. 
It is only restricted by the bounds \eqref{eq:kernel-bounds2} on the full kernel $\kernelcol = \kernelcol[\dif] + \kernelcol[\adv]$. 
The integral 
\begin{align*}
	\tilde\factorcol_\adv(\xcoord, \nvcoord') &= \intVn{\kernelcol[\adv](\xcoord, \nvcoord, \nvcoord') }
\end{align*}  
in general still depends on the direction $\nvcoord'$ and can even be negative. 
We define the normalization factor
\begin{align*}
	\factorcol[\adv](\xcoord) :=  \frac{1}{\speed}\max_{\nvcoord' \in \nvdomain} \lbrace \abs{\tilde\factorcol_\adv(\xcoord, \nvcoord')}\rbrace. 
\end{align*}
Finally, birth and death of particles enters the model via the source term 
\begin{align*}
	\Source \PD &=  \factorprol(\xcoord, \PM) \nSource \PD. 
\end{align*}
The net growth rate $\factorprol(\xcoord, \PM)$ depends on the local particle density $\PM = \intVn{\PD(\nvcoord)}$.
The operator $\nSource$ accounts for direction changes during proliferation. 
We define the growth rate such that the source is normalized, i.e.,  $\int_{\vdomain} \nSource\PD d\vcoord = \PM$. 

\subsection{Parabolic scaling and diffusion limit}
\label{sec:parabolic-scaling}
To derive the diffusion limit of \eqref{eq:lke-general}, it is helpful to write it in a dimensionless form. 
Therefore we introduce non-dimensional coordinates via $\xcoord = \xref \normal \xcoord$, $\tcoord = \tref \normal \tcoord$ together with the non-dimensional particle distribution $\PD(\tcoord,\xcoord,\nvcoord) = \PD_0 \normal \PD(\ntcoord, \nxcoord, \nvcoord)$ and $\factorcol[\dif](\xcoord) = \factordifref \nfactorcol[\dif](\normal \xcoord)$, $\factorcol[\adv](\xcoord) = \frac{\factoradvref}{\xref} \nfactorcol[\adv](\normal \xcoord), \factorprol(\xcoord, \PM) = \factorprolref \nfactorprol(\normal \xcoord, \normal \PM)$. 
With this we can define dimensionless kernels via $\kernelcol[\dif](\xcoord, \nvcoord, \nvcoord') = \factordifref \nfactorcol[\dif](\nxcoord) \nkernelcol[\dif](\nxcoord, \nvcoord, \nvcoord') $ and $\kernelcol[\adv](\xcoord, \nvcoord, \nvcoord') = \frac{\factoradvref}{\xref} \nfactorcol[\adv](\nxcoord) \nkernelcol[\adv](\nxcoord, \nvcoord, \nvcoord') $.
The dimensionless turning operators are 
\begin{align*}
	\nLCO[\dif] \normal\PD &= \intVnprime{ \nkernelcol[\dif](\nxcoord, \nvcoord, \nvcoord') \normal \PD(\nvcoord') - \nkernelcol[\dif](\nxcoord, \nvcoord', \nvcoord) \normal \PD(\nvcoord)}, \\ 
	\nLCO[\adv] \normal\PD &= \intVnprime{ \nkernelcol[\adv](\nxcoord, \nvcoord, \nvcoord') \normal \PD(\nvcoord') - \nkernelcol[\adv](\nxcoord, \nvcoord', \nvcoord) \normal \PD(\nvcoord)}, 
\end{align*}
and finally a non-dimensional form of \eqref{eq:lke-general} is 
\begin{align}
	\label{eq:lke-dimensionless}
	\partial_{\normal \tcoord} \normal \PD + \frac{\tref \vref}{\xref} \nabla_{\normal \xcoord} (\normal \vcoord \normal \PD) &= \tref \factordifref  \nfactorcol[\dif](\normal \xcoord) \nLCO[\dif] \normal \PD + \frac{\tref \vref}{\xref} \factoradvref \nfactorcol[\adv](\normal \xcoord) \nLCO[\adv] \normal \PD + \tref \factorprolref \nfactorprol(\normal \xcoord, \normal \PM) \nSource \normal \PD.
\end{align} 
We recognize the Strouhal number $\St = \frac{\xref}{\vref \tref}$, a Knudsen number for turning events  $\Knt = \frac{1}{\factordifref \tref}$, and a Knudsen number for proliferation events $\Knp = \frac{1}{\factorprolref \tref}$. 
Using these characteristic numbers and dropping the hats everywhere, we write the equation as 
\begin{align}
	\dt \PD + \frac{1}{\St} \divx (\vcoord \PD) = \frac{1}{\Knt} \factorcol[\dif](\xcoord) \LCO[\dif] \PD + \frac{\factoradvref}{\St} \factorcol[\adv](\xcoord) \LCO[\adv] \PD + \frac{1}{\Knp} \factorprol(\xcoord, \PM) \Source \PD
\end{align}
on the unit domain 
\begin{align*}
\normal \Omega_{txv} &= [0,1]\times \normal{\Omega}_x \times \US, \\
\normal \Omega_x & \subseteq [0,1]^\spacedim.
\end{align*}
In accordance with \cite{jin1996numerical}, we take the parabolic scaling parameter
\begin{align*}
\pareps := \frac{\Knt}{\St} = \frac{\vref}{\xref \factordifref}
\end{align*}
as the ratio of mean free path and domain length. 
To make the parabolic scaling apparent, we write \eqref{eq:lke-dimensionless} as 
\begin{align}
\label{eq:lke-scaled}
\dt \PD + \frac{\pardel}{\pareps} \divx (\vcoord \PD) &= \frac{\pardel}{\pareps^2}\factorcol[\dif](\xcoord) \LCO[\dif] \PD + \frac{\pardel \parnu}{\pareps} \factorcol[\adv](\xcoord) \LCO[\adv] \PD + \parthet \factorprol(\xcoord, \PM) \Source \PD,
\end{align}
with the parameters $\pardel = \frac{\Knt}{\St^2}$, $\parnu = \factoradvref$, $\parthet = \frac{1}{\Knp}$. 
In the literature usually $\pardel = \parthet = 1, \parnu = 0$ is assumed (see e.g. \cite{lemou2008ap, jin1996numerical, jin2000uniformly, klar1998asymptotic} ), which is not a problem from a theoretical perspective. 
From the perspective of the application the characteristic numbers are determined by the physical parameters and thus cannot be chosen arbitrarily. 
For fixed characteristic numbers $\pardel, \parnu, \parthet$, equation \eqref{eq:lke-scaled} converges to an advection-diffusion equation for the density $\density(\tcoord,\xcoord)$ as the parabolic scaling parameter approaches zero:
\begin{align}
\dt \density + \pardel \divx \left( \frac{1}{\factorcol[\dif]}  \divx \left(\density \ints{\vcoord \LCO[\dif]^{-1}( \vcoord \FD) \right)} - \frac{\parnu \factorcol[\adv] }{\factorcol[\dif]} \density \ints{\vcoord \LCO[\dif]^{-1} \LCO[\adv] \FD}\right) &= \parthet \factorprol(\xcoord, \density) \density. 
\label{eq:diffusion-limit-general}
\end{align}
Herein we use the angle brackets
\begin{align*}
	\ints{\placeholder} &= \intVn{\placeholder~}
\end{align*}
as shorthand notation for the integral over the unit sphere. 
We identify the symmetric positive definite diffusion tensor
\begin{align}
\label{eq:general-diffusion-tensor}
\Difftens:= -\frac{1}{\factorcol[\dif]}\ints{\LCO[\dif]^{-1}(\vcoord\FD) \vcoord\trans},
\end{align}
and the drift vector
\begin{align}
	\label{eq:general-drift-vector}
	\drift := -\frac{\factorcol[\adv]}{\factorcol[\dif]}\ints{\vcoord \LCO[\dif]^{-1}\LCO[\adv] \FD}.
\end{align}
Modulo hats, the diffusion equation transformed back to physical coordinates is
\begin{align*}
	\dt \density - \frac{\pardel \xref^2}{\tref \Difftens_0} \divx \left(\divx(\density \Difftens) - \frac{\parnu \Difftens_0}{\drift_0 \xref} \drift \density\right) = \frac{\parthet}{\factorprolref \tref} \factorprol(\xcoord, \density) \density,
\end{align*}
with a characteristic diffusion speed $\Difftens_0$ and a characteristic drift speed $\drift_0$ related to the microscopic scales via
\begin{align*}
	\Difftens_0 &= \frac{\pardel\xref^2}{\tref} = \frac{\speed^2}{\factordifref}, \\
	\drift_0 &= \frac{\parnu \Difftens_0}{\xref} = \frac{\speed^2 \factoradvref}{\xref \factordifref}. 
\end{align*}
Then finally the parabolic limit of \eqref{eq:lke-general} in physical coordinates is 
\begin{align}
	\label{eq:diffusion-limit-physical}
	\dt \density - \divx \left(\divx (\Difftens \density)  - \drift \density \right) &= \factorprol(\xcoord, \density) \density. 
\end{align}

A formal proof of the limit via a Hilbert expansion in $\pareps$ can be found in \cite{othmer2000diffusion, EHKS14,corbin2018higher}.
We will not repeat this proof here but rather use the micro-macro decomposition in the next section to compute the limit in a less rigorous way.
In any case, the limit only exists if the operator $\LCO[\dif]$ is invertible on an appropriate space. 
This is guaranteed by the following \lemref{lem:lcol-properties} from \cite{bellomo2006onset}. 

\begin{definition}[Weighted $\ltwospace$ space]
	With  $\weightedltwospace$ we denote the $\ltwospace$-space on $\US$ with respect to the weighted scalar product
	\begin{align*}
	\left( f(\nvcoord),g(\nvcoord) \right)_{\FD} = \ints{\frac{f(\nvcoord) g(\nvcoord)}{\FD(\nvcoord)} }. 
	\end{align*}
\end{definition}

\begin{lemma}[Properties of ${\LCO[\dif]}$]
	\label{lem:lcol-properties}
	Under assumptions  \eqref{eq:kernel-bounds1}, \eqref{eq:kernel-equilibrium}, the turning operator $\LCO[\dif]: \weightedltwospace \mapsto \weightedltwospace$ has the following properties for each $\xcoord \in \Omega_x$:
	\begin{enumerate}
		\item $\LCO[\dif]$ is self-adjoint;	
		\item The one-dimensional nullspace of $\LCO[\dif]$ is $\Nsp(\LCO[\dif]) = \vecspan{\FD}$;
		\item There exists a unique solution to $\LCO[\dif] \PD = \PP$ for every $\PP \in \Nsporth$, i.e. $\PP$ such that $\left( \PP, \FD \right)_{\FD} = \intVn{\PP(\nvcoord)} = 0$.
	\end{enumerate}
\end{lemma}
\subsection{A simple haptotaxis model for glioma}
\label{sec:glioma-model}
For the computations we use a model for haptotaxis induced glioma migration from \cite{EHKS14,corbin2018higher} that can be cast into the general setting. 
Because it would exceed the scope of this paper to discuss the details of its derivation we only give a brief summary. 
First of all, assume that a field of symmetric positive definite tensors $\DW(\xcoord) : \xdomain \mapsto \R^{3\times 3}$ is given. 
In practice, diffusion tensor imaging (DTI) provides piecewise constant measurements of the diffusion of water molecules through the tissue \cite{LeBihan2001DTI}.  
As in \cite{EHS} we use this information to estimate the directional distribution of extracellular matrix (ECM) fibers $\FD[\DW](\xcoord, \vcoord)$ and the fraction of volume $\VF[\DW](\xcoord)$ these fibers occupy. 
One important aspect of the model is that glioma cells use ECM fibers for contact guidance, i.e., they align themselves to the fibers. 
The fiber distribution $\FD$ plays the role of the collision equilibrium and therefore should fulfill assumptions \eqref{eq:equilibrium-assumptions} and \eqref{eq:kernel-equilibrium}. 
A simple estimate for the fiber distribution is the so-called peanut distribution
\begin{align}
\label{eq:peanut}
\FD(\xcoord,\nvcoord) &= \frac{3}{4 \pi \trace (\DW)} (\nvcoord \trans \DW \nvcoord),
\end{align}

The turning rate for the first turning operator is constant, i.e. $\factorcol[\dif] = \lambda_0$, and the turning kernel $\kernelcol[\dif]$ is proportional to the fiber distribution
\begin{align*}
	\kernelcol[\dif](\xcoord,\nvcoord,\nvcoord') = \lambda_0 \FD(\xcoord,\nvcoord),
\end{align*}
such that the turning operator $\LCO[\dif] = \lambda_0 \left(\ints{\PD} \FD - \PD \right)$ is a simple relaxation to local equilibrium. 
For any $\phi \in \Nsporth$, i.e., $\ints{\phi} = 0$, the inverse of $\LCO[\dif]$ is simply
\begin{align}
\LCO[\dif]^{-1} (\phi) = -\frac{1}{\lambda_0}\phi. 
\label{eq:glioma-lco-inverse}
\end{align}
The turning perturbation $\LCO[\adv]$ stems from a subcellular model that includes internal state changes of cells.
In this model cells change their turning behavior according to the ECM concentration.
This results in a collective movement in direction of the fiber gradient:  
\begin{align*}
	\kernelcol[\adv](\xcoord,\nvcoord,\nvcoord') &= -\lambda_H(\xcoord) \speed \left( \nabla_x \VF(\xcoord) \cdot \nvcoord' \right) \FD(\xcoord,\nvcoord), \\
		\factorcol[\adv](\xcoord) &= \lambda_H(\xcoord) \norm{\nabla_x \VF(\xcoord)}. 
\end{align*}

For the source, we consider logistic growth towards the carrying capacity $\carryingcapacity$, thus the growth rate is given by
\begin{align*}
	\factorprol(\xcoord, \PM) = \factorprolref \left(1 - \frac{\PM}{\carryingcapacity}\right).
\end{align*}
We assume that no changes in direction occur during growth, which is expressed by
\begin{align*}
	\nSource \PD = \PD. 
\end{align*} 
For a more detailed discussion the interested reader is referred to \cite{EHKS14,EHS,corbin2018higher}.
With these definitions, the glioma equation in physical coordinates reads
\begin{align*}
	\dt \PD + \speed \divx (\nvcoord \PD) &= \lambda_0 \left( \ints{\PD} \FD(\xcoord,\nvcoord) - \PD \right) - \speed \lambda_H(\xcoord) \nabla_x\VF(\xcoord) \cdot \left( \ints{\nvcoord \PD} \FD(\xcoord,\nvcoord) - \nvcoord\PD \right) + \factorprol(\xcoord, \PM) \PD.
\end{align*}
After applying the parabolic scaling from \secref{sec:parabolic-scaling}, the glioma model in dimensionless form becomes
\begin{align}
	\dt \PD + \frac{\pardel}{\pareps} \divx (\vcoord \PD) &= \frac{\pardel}{\pareps^2} \left(\FD \ints{\PD} - \PD\right) - \frac{\pardel\parnu}{\pareps} \hat \lambda_H \nabla_x \VF \cdot \left(\FD \ints{\PD \vcoord} - \PD \vcoord\right) + \parthet \nfactorprol \PD, 
	\label{eq:lke-glioma-scaled}
\end{align}
with $\lambda_H = \frac{\lambda_1}{\lambda_0} \hat \lambda_H$  and the characteristic numbers
\begin{align*}
	\pareps = \frac{\speed}{\xref \lambda_0}, \quad \pardel = \frac{\speed^2}{\lambda_0} \frac{\tref}{\xref^2}, \quad \parnu = \frac{\lambda_1}{\lambda_0}, \quad \parthet = \factorprolref \tref.
\end{align*}
The diffusion approximation is given by
\begin{align}
\dt \density - \pardel \divx \left( \divx (\density \DT) -  \parnu \driftT \density \right) = \parthet \factorprol \density
\label{eq:DiffusionLimitGlioma}.
\end{align}
Using the inversion formula \eqref{eq:glioma-lco-inverse}, the tumor diffusion tensor and drift resulting from the peanut distribution \eqref{eq:peanut} are given by
\begin{align}
	\label{eq:DT-glioma}
	\DT &= \ints{vv\FD} = \frac{1}{5} \left(\Identity + \frac{2\DW}{\trace \DW}  \right).  \\
	\label{eq:drift-glioma}
	\driftT &= \hat\lambda_H \nabla_x \VF \cdot \DT. 
\end{align}

\section{Micro-Macro decomposition and the diffusion limit}
\label{sec:micro-macro}
In the next section, we follow the work of Lemou and Mieussens \cite{lemou2008ap} quite closely to perform a micro-macro decomposition of equation \eqref{eq:lke-scaled} in the parabolic dimensionless form. 
This serves as the starting point for the numerical discretization scheme.
From \lemref{lem:lcol-properties} we recall the nullspace $\Nsp(\LCO[\dif]) = \vecspan{E}$ and range $\Range(\LCO[\dif]) = \Nsporth(\LCO[\dif])$ of the turning operator. 
Orthogonal projections onto those spaces are 
\begin{align*}
\Nspproj(\phi) &= \ints{\phi} \FD,  \\
\Ptbproj(\phi) &= \phi - \ints{\phi} \FD, 
\end{align*}
respectively. 
Using these projections, we split the particle distribution into an equilibrium part and a perturbation: 
\begin{equation}
	\begin{aligned}
		\PD &= \Nspproj \PD + \Ptbproj \PD \\
		    &= \PM \FD + \pareps \PP. 
		    \label{eq:APsplitPD}
	\end{aligned}
\end{equation}
Here, $\PM(\tcoord, \xcoord) = \ints{\PD}$ is the local particle density. 

Now the kinetic equation is split into a system of two equations\textemdash one for the macroscopic density $\PM$ and one for the microscopic perturbation $\PP$. 
We obtain the $\PM$-equation by inserting the perturbation formula \eqref{eq:APsplitPD} into \eqref{eq:lke-scaled} and applying the projection $\Nspproj$: 
\begin{align}
	\dt \PM + \pardel \divx \ints{\PP \vcoord} &= \parthet \factorprol \PM,
	\label{eq:APrho}
\end{align}
where we use the positivity and symmetry of the equilibrium \eqref{eq:equilibrium-assumptions} and the mass conservation $\ints{\LCO[\dif]} = 0, \ints{\LCO[\adv]} = 0$ of the turning operators \eqref{eq:lco-mass-conservation}. 
Then, applying $\Ptbproj$ to \eqref{eq:lke-scaled} and dividing by $\pareps$ gives 
\begin{align}
	\dt \PP + \frac{\pardel}{\pareps} \Ptbproj \divx (\vcoord \PP) &= -\frac{\pardel}{\pareps^2} \divx (\vcoord \PM \FD) + \frac{\pardel \factorcol[\dif]}{\pareps^2} \LCO[\dif] \PP + \frac{\pardel \parnu \factorcol[\adv]}{\pareps^2} \LCO[\adv] \PD + \frac{\parthet \factorprol}{\pareps} \Ptbproj \Source \PD,
	\label{eq:APremainder}
\end{align}
where we use
\begin{align*}
	\Nspproj \LCO \phi &= \ints{\LCO \phi} \FD \\
	                   &= 0,\\ 
	\Ptbproj \LCO[\dif] \PD &= \LCO[\dif] (\PM \FD + \pareps \PP) \\
	                        &= \pareps \LCO[\dif] \PP, \\
	\Nspproj \divx (\vcoord \PM \FD ) &= \ints{\divx \vcoord \PM \FD} \FD \\
	                                       &= \divx \ints{\vcoord \PM \FD} \FD \\
	                                       &= 0.
\end{align*}
Apart from the new $\LCO[\adv]$ term, this formulation coincides with the decomposition in \cite{lemou2008ap}. The authors of \cite{lemou2008ap} show, that\textemdash for compatible initial and boundary conditions\textemdash the micro-macro decomposition is equivalent to the original kinetic equation \eqref{eq:lke-scaled}.

It is easy to see the diffusion limit from the decomposition in a rather informal way. 
In the limit of $\pareps \rightarrow 0$, only the $\frac{1}{\pareps^2}$ terms remain in \eqref{eq:APremainder} and thus it is reduced to 
\begin{align*}
	\PP_0 &=  \frac{1}{\factorcol[\dif]}  \LCO[\dif]^{-1} \left( \divx (\vcoord \density \FD)  - \parnu \factorcol[\adv] \density \LCO[\adv] \FD \right).
\end{align*}
Since $\ints{\vcoord\FD} = \ints{\LCO[\adv] \PD} = 0$, \lemref{lem:lcol-properties} assures that the inverse of $\LCO[\dif]$ in this expression exists and is unique.
Inserting this into the macro equation \eqref{eq:APrho} immediately gives the diffusion limit \eqref{eq:diffusion-limit-general}. 

The main idea behind the asymptotic preserving scheme is to do something similar in a discrete way. 
First the perturbation $\PP^{n+1}$ on the next time-level is computed using the micro equation, then this is inserted into the macro equation to update the density $\PM^{n+1}$.

\section{The asymptotic preserving method}
\label{sec:ap-scheme}
In general, a numerical scheme is called asymptotic preserving (AP) with respect to a scaling limit, if it converges to a valid scheme for the limit equation as $\pareps \rightarrow 0$ and the spatial discretization is fixed. 
The stability criterion for the time step size $\timestep$ must be bounded from below by a positive value independent of $\pareps$. 
The main objective of this work is to develop such an asymptotic preserving scheme for the kinetic equation \eqref{eq:lke-scaled}. 

We start from the micro-macro decomposition from the previous \secref{sec:micro-macro} and write it as 
\begin{equation}
\label{eq:continuous-system}
\begin{alignedat}{2}
\dt \PM &= \macroexpl(\PM, \PP) &+& \macroimpl(\PM,\PP), \\
\dt \PP      &= \left(\microexplfibres(\PM) + \microexpl(\PM, \PP) \right)&+& \microimpl(\PM,\PP). 
\end{alignedat}
\end{equation}
Here the individual terms are grouped into those that will later be discretized explicitly in time 
\begin{equation}
\label{eq:continuous-system-explicit-terms}
\begin{aligned}
\macroexpl(\PM,\PP)   &= -\pardel \divx \ints{\PP \vcoord} + \parthet \factorprol \PM,\\
\microexplfibres(\PM) &= -\frac{\pardel}{\pareps^2} \divx (\vcoord \PM \FD), \\
\microexpl(\PM, \PP)  &=  -\frac{\pardel}{\pareps} \Ptbproj (\divx (vg)) + \frac{\pardel \parnu \factorcol[\adv]}{\pareps^2} \LCO[\adv] \PD + \frac{\parthet \factorprol}{\pareps} \Ptbproj (\Source \PD), 
\end{aligned}
\end{equation}
and those that will be discretized partially implicit
\begin{equation}
\label{eq:continuous-system-implicit-terms}
\begin{aligned}
\macroimpl(\PM,\PP)   &= 0, \\
\microimpl(\PM, \PP)  &= \frac{\pardel \factorcol[\dif]}{\pareps^2} \LCO[\dif] \PP.
\end{aligned}
\end{equation}
In \cite{lemou2008ap} the authors argued that it is enough to treat only the term $\LCO[\dif]$ in an implicit way to get an AP scheme. 
We call the first-order scheme derived from the micro-macro decomposition in the form \eqref{eq:continuous-system}-\eqref{eq:continuous-system-implicit-terms}, in which only $\LCO[\dif]$ is treated implicitly, \mischeme{1}; and the second-order scheme \mischeme{2}.

But it is also possible to solve the source and $\LCO[\adv]$ terms implicitly in time.
That is, we regroup the terms into
\begin{equation}
\label{eq:continuous-system-iv}
\begin{aligned}
\macroivexpl(\PM,\PP)   &= -\pardel \divx \ints{\PP \vcoord}, \\
\microivexplfibres(\PM) &= -\frac{\pardel}{\pareps^2} \divx (\vcoord \PM \FD), \\
\microivexpl(\PM, \PP)  &=  -\frac{\pardel}{\pareps} \Ptbproj (\divx (vg)),  \\
\macroivimpl(\PM,\PP)   &= \parthet \factorprol \PM, \\
\microivimpl(\PM, \PP)  &= \frac{\pardel \factorcol[\dif]}{\pareps^2} \LCO[\dif] \PP + \frac{\pardel \parnu \factorcol[\adv]}{\pareps^2} \LCO[\adv] \PD + \frac{\parthet \factorprol}{\pareps} \Ptbproj (\Source \PD).
\end{aligned}
\end{equation}
and solve $\macroivexpl, \microivexplfibres, \microivexpl$ explicitly and $\macroivimpl, \microivimpl$ implicitly. 
This variant of the scheme will be called \ivscheme{1}, or \ivscheme{2}. 
In the following sections, we will see that the implicit time update for this scheme can still be done on each grid cell separately. 

\subsection{Space discretization}

In \cite{lemou2008ap} the authors discretize the micro and macro equation with finite differences on staggered grids in one space dimension. 
To generalize the method to arbitrary dimension $\spacedim$, we reformulate the method in the context of finite volumes on primal-dual mesh pairs. 

Although the implementation supports only tensor-product grids at the moment, we write the scheme for conforming polyhedral meshes. 
This has several benefits. 
Most aspects of the scheme do not depend on the tensor-product structure, and also the implementation in DUNE (see \cite{dune-web-page}) is grid-agnostic in most parts. 
The general notation is quite close to the implementation, which helps understanding the code and also will make an implementation on unstructured conforming meshes easier. 
We choose a notation that is similar to that in \cite{buet2012design}.
We use the symbol $\mentity{\midx}$ wherever any kind of entity on the grid can be inserted(cell, face, edge, dual cell, \dots). 
The index $\midx$ is used to label these generic entities.

Only considering topology, the dual mesh belonging to a primal mesh is defined as follows: 
Each cell in the original mesh is identified with a vertex in the dual mesh and each primal vertex with a dual cell. 
Wherever two primal cells intersect in a face, two dual vertices are connected with an edge and where two primal vertices are connected, there is a face between two dual cells. 

We always use the indices $\pcellidx, \pcellidxnb \in \mathbb{N}$ to label cells $\pcell, \pcellnb$ in the primal grid and $\pvertidx, \pvertidxnb$ to identify primal vertices $\pvert, \pvertnb$. 
Considering the primal-dual mapping, any primal cell index $\dvertidx$ also identifies a dual vertex $\dvert$ and a primal index $\dcellidx$ corresponds to a dual cell $\dcell$. 
In one mesh two cells $\mcell{\midx}, \mcell{\midx'}$ are neighbors, if they intersect in a face $\mface{\midx}{\midx'} = \mcell{\midx} \isect \mcell{\midx'}$. 
Then the two vertices $\mvert{\midx}, \mvert{\midx'}$ in the other grid are also neighbors, i.e., they are connected with an edge $\medge{\midx}{\midx'}$.
In this sense, the neighbors of an index $\midx$ are those indices $\midx'$ for which in one mesh the corresponding cells are neighbors and thus in the other grid the corresponding vertices are neighbors. 
We write $\mneighbors{\midx}$ for the set of all neighbors of $\midx$. 
A related concept is the adjacency between entities of different dimension. 
If the edge $\pedge$ is part of the cell $\pcell$ we say that $\pcell$ is adjacent to $\pedge$, and denote this by $\pcellidx \in \madjacents{\pvertidx}{\pvertidxnb}$.
The index pair $(\pvertidx, \pvertidxnb)$ also identifies a dual face $\dface$, thus $\madjacents{\dcellidx}{\dcellidxnb}$ equivalently is the set of all dual vertices $\dvert$ that are part of that face. 
Lastly we denote the set of vertices of a cell $\midx$ with $\mverticesof{\midx}$. 
The example mesh in \figref{fig:annotated-mesh-2d-a} is helpful to visualize these definitions. 
\begin{figure}
	\def\localpath{figures/primal_dual_mesh/}
	\centering
	\parbox{\figuretwocol}{%
		\centering
		\tikztitle{Faces}
		\settikzlabel{fig:annotated-mesh-2d-a}
		\withfiguresize{\figuretwocol}{\figuretwocol}{\externaltikz{annotated_mesh_2d_a}{
%
%
	\input{\localpath mesh_colors}
	\input{\localpath mesh_macros}
\begin{tikzpicture}
	\makemeshpair{B/0.38/0.022,C/0.026/0.40,D/0.38/0.38,E/0.95/0.36,F/0.60/0.60,G/0.22/0.81,H/0.80/0.90}{poly2/C/D/B,poly3/B/E/D,poly4/E/F/D,poly5/C/G/D,poly6/G/D/F,poly7/G/H/F,poly8/F/H/E}{}
	\path[primal cell] (pvertD.center) -- (pvertF.center) -- (pvertG.center) -- cycle; 	
	\path[] (pvertD) -- node[midway, coord]{$\Omega_r$} (dvertpoly2);  
	
	\path[dual cell] (dvertpoly2.center) -- (pfacemidBD.center) 
	-- (dvertpoly3.center) -- (pfacemidDE.center)
	-- (dvertpoly4.center) -- (pfacemidDF.center) 
	-- (dvertpoly6.center) -- (pfacemidDG.center) 
	-- (dvertpoly5.center) -- (pfacemidCD.center)
	-- (dvertpoly2.center) -- (pfacemidBD.center)
	--cycle; 
	\path (dvertpoly6) -- node[midway, coord]{$\Omega_j$} (pvertG);

	\path[draw, primal face, primalhigh] (pvertD) -- node[pos=0.8, anchor=north west, coord] {\textcolor{primalcolorhigh}{$\partial \Omega_{j,k}$}} (pvertF); 	
	
	\path[draw, dual face, dualhigh] (dvertpoly6) -- (pfacemidDF) -- node[below=1ex,coord] {\textcolor{dualcolorhigh}{$\partial \Omega_{r,s}$}} (dvertpoly4); 
	
	\path (pvertD) node[primal vertex, opacity=1] {$x_r$};
	\path (pvertF) node[primal vertex, opacity=1] {$x_s$};
	
	\path (dvertpoly6) node[dual vertex, opacity=1] {$x_j$};
	\path (dvertpoly4) node[dual vertex, opacity=1] {$x_k$};
	
\end{tikzpicture}
}}
	}
	\hspace{\figurehorizontalsep}
	\parbox{\figuretwocol}{%
		\centering
		\tikztitle{Facets}
		\settikzlabel{fig:annotated-mesh-2d-b}
		\withfiguresize{\figuretwocol}{\figuretwocol}{\externaltikz{annotated_mesh_2d_b}{	\input{\localpath mesh_colors}
	\input{\localpath mesh_macros}
\begin{tikzpicture}
	\makemeshpair{B/0.38/0.022,C/0.026/0.40,D/0.38/0.38,E/0.95/0.36,F/0.60/0.60,G/0.22/0.81,H/0.80/0.90}{poly2/C/D/B,poly3/B/E/D,poly4/E/F/D,poly5/C/G/D,poly6/G/D/F,poly7/G/H/F,poly8/F/H/E}{}
	\path[primal cell] (pvertD.center) -- (pvertF.center) -- (pvertG.center) -- cycle; 	 
	
	\path[dual cell] (dvertpoly2.center) -- (pfacemidBD.center) 
	-- (dvertpoly3.center) -- (pfacemidDE.center)
	-- (dvertpoly4.center) -- (pfacemidDF.center) 
	-- (dvertpoly6.center) -- (pfacemidDG.center) 
	-- (dvertpoly5.center) -- (pfacemidCD.center)
	-- (dvertpoly2.center) -- (pfacemidBD.center)
	--cycle; 
			
	\path[opacity=0.3, fill=primalcolor!50!dualcolor] (pvertD.center) -- (pfacemidDF.center) -- (dvertpoly6.center) -- (pfacemidDG.center) --cycle; 
	\path (pvertD) -- node[midway, coord]{$\omega_r^j$} (dvertpoly6);
	
	\path[draw, primal face, primalhigh] (pvertD) -- node[anchor=north west, coord] {$\partial \omega_{j,k}^{r}$} (pfacemidDF); 
	\path[draw, dual face, dualhigh] (pfacemidDF) -- node[anchor= south west,coord] {$\partial \omega_{r,s}^{j}$} (dvertpoly6); 
	
	\path (pvertD) node[primal vertex, opacity=1] {$x_r$};
	\path (pvertF) node[primal vertex, opacity=1] {$x_s$};
	
	\path (dvertpoly6) node[dual vertex, opacity=1] {$x_j$};
	\path (dvertpoly4) node[dual vertex, opacity=1] {$x_k$};
	
\end{tikzpicture}}}
	}
	\caption{The primal-dual mesh pair in two dimensions. The primal cell $\pcell$ is marked green and the dual cell $\dcell$ in gray. \ref{fig:annotated-mesh-2d-a}: Highlighted are the primal face $\pface$ and the dual face $\dface$. \ref{fig:annotated-mesh-2d-b}: Highlighted are the subcell $\psubcell = \pcell \isect \dcell$, the primal facet $\pfacet = \pface \isect \dcell$ and the dual facet $\dfacet = \dface \isect \pcell$.}
	\label{fig:annotated-mesh-2d}
\end{figure}

Given a primal mesh, the topological mapping alone does not define the geometry of the dual mesh uniquely. 
For instance the dual vertex $\dvert$ can be anywhere inside the primal cell $\pcell$. 
For the numerical scheme we need to know the geometry of the dual cells and especially their faces. 
First note that a dual face $\dface$, which is the intersection between two dual cells, does not need to be planar.
In two space dimensions it can be constructed, however, from one planar facet $\dfacet = \dface \isect \pcell$ for each intersection with an adjacent primal cell $\pcell; \pcellidx \in \madjacents{\dcellidx}{\dcellidxnb}$. 
The facet $\dfacet$ is just the line $\medge{\pcellidx}{\dcellidx, \dcellidxnb}$ between the primal 'cell center' $\dvert$ and some arbitrary point $\mvert{\dcellidx, \dcellidxnb}$ on the edge $\pedge$(which coincides with a face $\pface$, for some $\pcellidxnb$).
This construction is depicted in \figref{fig:annotated-mesh-2d-b} and is identical to the definition of a control volume in \cite{buet2012design}. 
In three space dimensions the construction is similar but a bit more complicated. 
For a sketch of the construction, see \figref{fig:annotated-mesh-3d}. 
Because the primal mesh is polyhedral and conforming, the facet $\dfacet$ is bounded by line segments connecting the four points $\dvert, \mvert{\pcellidx, \pcellidxnb}, \mvert{\dcellidx, \dcellidxnb}, \mvert{\pcellidx, \pcellidxnb'}$.
The indices $\pcellidxnb, \pcellidxnb' \in \mneighbors{\pcellidx} \isect \madjacents{\dcellidx}{\dcellidxnb}$ label those two neighbors of cell $\pcell$ that have $\pedge$ as an edge. 
With $\mvert{\pcellidx, \pcellidxnb}, \mvert{\pcellidx, \pcellidxnb'}$ we denote arbitrary points on the faces $\pface, \mface{\pcellidx}{\pcellidxnb'}$, for example their barycenters. 
As in the two-dimensional setting, $\mvert{\dcellidx, \dcellidxnb}$ is an arbitrary point on the edge $\pedge$. 
In general, the four points do not have to lie in a plane. 
Thus if we want to have a polyhedral dual mesh, the facet $\dfacet$ must be split into two triangles $\dfacetone \union \dfacettwo = \dfacet$ defined by the triplets $\dvert, \mvert{\pcellidx, \pcellidxnb}, \mvert{\dcellidx, \dcellidxnb}$, and $\dvert, \mvert{\dcellidx, \dcellidxnb}, \mvert{\pcellidx, \pcellidxnb'}$. 
For tensor product grids and tetrahedral meshes (see \cite{weiss2013primal}), the four points lie in a plane if they are chosen as the barycenters of their respective entities, making the split into triangles unnecessary. 
\begin{figure}
	\def\localpath{figures/primal_dual_mesh/}
	\centering
	\withfiguresize{\figureonecol}{\figureonecol}{\externaltikz{annotated_mesh_3d}{	\input{\localpath mesh_colors}
	\input{\localpath mesh_macros}	
\begin{tikzpicture}
	[coordmark/.style={draw, circle, minimum width=2pt, inner sep = 0pt}]
	\newlength\cubewidth
	\setlength{\cubewidth}{0.3\figurewidth}

	\def\orthoprojalpha{-70}
	\def\orthoprojbeta{0}
	\def\orthoprojgamma{10}
	\def\orthoprojx#1#2#3{#1 * cos(\orthoprojbeta)*cos(\orthoprojgamma) - #3 * sin(\orthoprojbeta) - #2 * cos(\orthoprojbeta) * sin(\orthoprojgamma)}
	\def\orthoprojy#1#2#3{#1 * (cos(\orthoprojalpha)*sin(\orthoprojgamma) - cos(\orthoprojgamma)*sin(\orthoprojalpha)*sin(\orthoprojbeta)) + #2 * (cos(\orthoprojalpha)*cos(\orthoprojgamma) + sin(\orthoprojalpha)*sin(\orthoprojbeta)*sin(\orthoprojgamma)) - #3 * cos(\orthoprojbeta)*sin(\orthoprojalpha) }
	\def\orthoprojz#1#2#3{#1 * (sin(\orthoprojalpha*sin(\orthoprojgamma))+ cos(\orthoprojalpha)*cos(\orthoprojgamma)*sin(\orthoprojbeta))+ #2 * (cos(\orthoprojgamma)*sin(\orthoprojalpha) - cos(\orthoprojalpha)*sin(\orthoprojbeta)*sin(\orthoprojgamma)) + #3 * (cos(\orthoprojalpha)*cos(\orthoprojbeta))}
	\def\orthoproj#1#2#3#4{%
		\path let%
		\n1={\orthoprojx{#1}{#2}{#3}},%
		\n2={\orthoprojy{#1}{#2}{#3}},%
		\p1=(\n1, \n2)%
		in (\p1) coordinate(#4);%
	}
	
	\def\makecubecoords#1#2#3#4#5{%
		\foreach \pj/\px/\py/\pz in {A/0/0/0, B/#4/0/0, C/0/#4/0, D/#4/#4/0, E/0/0/#4, F/#4/0/#4, G/0/#4/#4, H/#4/#4/#4} 
		{
			\path let \n1={\px + #1}, \n2={\py + #2}, \n3={\pz + #3}, \n4 = {\orthoprojx{\n1}{\n2}{\n3}}, \n5 = {\orthoprojy{\n1}{\n2}{\n3}}, \p1=(\n4,\n5) in (\p1) coordinate(#5\pj);
		}
		\path let \n1={#1 + 0.5*#4}, \n2={#2 + 0.5*#4}, \n3={#3 + 0.5*#4}, \n4 = {\orthoprojx{\n1}{\n2}{\n3}}, \n5 = {\orthoprojy{\n1}{\n2}{\n3}}, \p1=(\n4,\n5) in (\p1) coordinate(#5Center);
	}
		
	\def\drawcube#1#2{%
		\path[#2] (#1A) -- (#1B) -- (#1D) -- (#1C) -- (#1A); 
		\path[#2] (#1E) -- (#1F) -- (#1H) -- (#1G) -- (#1E); 
		\path[#2] (#1A) -- (#1E); 
		\path[#2] (#1B) -- (#1F); 
		\path[#2] (#1C) -- (#1G); 
		\path[#2] (#1D) -- (#1H);
	}
	
	\makecubecoords{0}{0}{0}{\cubewidth}{cubemid}
	\makecubecoords{\cubewidth}{0}{0}{\cubewidth}{cuberight}
	\makecubecoords{0}{0}{\cubewidth}{\cubewidth}{cubetop}
	
	\drawcube{cubemid}{primal face}
	\drawcube{cuberight}{primal face, dashed}
	\drawcube{cubetop}{primal face, dashed}	
	
	\path (cubemidCenter) node[coordmark]{} node[above]{$x_j$}; 
	\path (cuberightCenter) node[coordmark]{} node[above]{$x_k$}; 
	\path (cubetopCenter) node[coordmark]{} node[above]{$x_{k'}$}; 
	
	\path (cubemidF) node[coordmark]{} node[above]{$x_r$}; 
	\path (cubemidH) node[coordmark]{} node[above]{$x_s$}; 
	
	\makecubecoords{0.5\cubewidth}{0}{0.5\cubewidth}{0.5\cubewidth}{dual}
	\drawcube{dual}{dual face}
	
	\orthoproj{\cubewidth}{0.5\cubewidth}{0.5\cubewidth}{facemidjk}
	\orthoproj{0.5\cubewidth}{0.5\cubewidth}{\cubewidth}{facemidjkp}
	\orthoproj{\cubewidth}{0.5\cubewidth}{\cubewidth}{edgemidrs}
	
	\path (facemidjk) node[coordmark]{} node[above]{$x_{j,k}$}; 
	\path (facemidjkp) node[coordmark]{} node[above]{$x_{j,k'}$}; 
	\path (edgemidrs) node[coordmark]{} node[above]{$x_{r,s}$};
	
	\path [dual face] (cubemidCenter) --  node[midway] {$\partial \omega_{r,s}^{j}$}(edgemidrs); 	
	\path [dual cell] (cubemidCenter) -- (facemidjk) -- (edgemidrs) -- (facemidjkp) -- cycle; 
\end{tikzpicture}}}
	\caption{The primal-dual mesh pair in three dimensions. Shown is the primal cell $\pcell$ (green wireframe) and the facet $\dfacet = \dface \isect \pcell$ (gray solid).}
	\label{fig:annotated-mesh-3d}
\end{figure}

We write the average of some function over the domain $\mentity{\midx}$ as
\begin{align*}
\mavg{\mentity{\midx}}{\placeholder} &:= \frac{1}{\abs{\mentity{\midx}}}\int_{\mentity{\midx}} \placeholder d\xcoord,
\end{align*} 
in which 
\begin{align*}
\abs{\mentity{\midx}} &= \int_{\mentity{\midx}} 1 d\xcoord
\end{align*}
is the volume of entity $\mentity{\midx}$.  
In the following, we derive the minimally implicit variant \mischeme{1} of the scheme. 
All that is required to obtain the variant with implicit volume terms \ivscheme{1} is a reordering of terms, analogously to \eqref{eq:continuous-system-iv}. 
Let $(\PM, \PP)$ be the solution of \eqref{eq:continuous-system}, with the average densities $\mavg{\dcell}{\PM}$ on dual cells, and the average perturbations $\mavg{\pcell}{\PP}$ on the primal cells. 
The projection of equation \eqref{eq:continuous-system} onto the cell averages is a finite system of equations for the values $\dPM \approx \mavg{\dcell}{\PM}$, $\pPP \approx \mavg{\pcell}{\PP}$ which approximate the averages of the exact solution.
We collect these values in the vectors $\PMvec = (\dots, \dPM, \PM_{\dcellidx+1} \dots)\trans$ and $\PPvec = (\dots, \PP_{\pcellidx}, \PP_{\pcellidx + 1}, \dots)\trans$ and write the resulting space-discrete system as
\begin{equation}
\label{eq:semi-discrete-system}
\begin{alignedat}{2}
\dt \PMvec &= \macroexplvec(\PMvec, \PPvec) &+& \macroimplvec(\PMvec, \PPvec) \\
\dt \PPvec      &= (\microexplfibresvec(\PMvec) + \microexplvec(\PMvec, \PPvec)) &+& \microimplvec(\PMvec, \PPvec),
\end{alignedat}
\end{equation}
using the same notation for the approximations of the projected operators. For instance we have $\macroexplvec(\PMvec, \PPvec) = (...,\macroexpl_{\dcellidx}(\PMvec, \PPvec), \macroexpl_{\dcellidx+1}(\PMvec, \PPvec),...)\trans$, where $\macroexpl_{\dcellidx}$ is an approximation to $\mavg{\dcell}{\macroexpl}$.
With second-order accuracy, the average $\mavg{\mentity{\midx}}{\placeholder}$ can be swapped with a product or a chained function, i.e. given functions  $u(\xcoord), w(\xcoord) \in C^2(\mentity{\midx})$, and $z(u) \in C^2(u(\mentity{\midx}))$ we have 
\begin{align*}
	\mavg{\mentity{\midx}}{u(\xcoord) w(\xcoord)} &= \mavg{\mentity{\midx}}{u(\xcoord)} \mavg{\mentity{\midx}}{w(\xcoord)} + \orderof{\gridsize^2} \\
	\mavg{\mentity{\midx}}{z(u(\xcoord))} &= z\left(\mavg{\mentity{\midx}}{u(\xcoord)}\right) + \orderof{\gridsize^2}.
\end{align*}
Up to second-order accurate approximations to the explicit operators on each cell are 
\begin{equation}
	\label{eq:semidiscrete-system-explicit-terms}
	\begin{aligned}
	\macroexpl_{\dcellidx} &= -\pardel \sum_{\dcellidxnb \in \mneighbors{\dcellidx}} \fluxrg  + \parthet \factorprol(\dPM) \dPM \\
	\microexplfibres_{\pcellidx} &= -\frac{\pardel}{\pareps^2} \sum_{\pcellidxnb \in \mneighbors{\pcellidx}} \fluxgr \\
	\microexpl_{\pcellidx}  &= -\frac{\pardel}{\pareps} \sum_{\pcellidxnb \in \mneighbors{\pcellidx}} \fluxgg +\frac{\pardel \parnu \factorcol[\adv, \pcellidx]}{\pareps^2} \LCO[\adv] \left(\tilde\PM_{\pcellidx} \FD_{\pcellidx} + \pareps \PP_{\pcellidx}\right) + \frac{\parthet \factorprol(\tilde \PM_{\pcellidx})}{\pareps} \Ptbproj \Source (\tilde\PM_{\pcellidx} \FD_{\pcellidx} + \pareps \pPP)
	\end{aligned}
\end{equation}
The average density on a primal cell $\tilde{\PM}_{\pcellidx}$ is not a degree of freedom of the scheme and needs to be computed from the averages on contributing dual cells: 
\begin{align}
	\label{eq:density-on-primal-cell}
	\tilde{\PM}_{\pcellidx} &= \frac{1}{\abs{\pcell}} \sum_{\dcellidx \in \mverticesof{\pcellidx}} \abs{\psubcell} \dPM.
\end{align}
The fluxes $\fluxrg$ are obtained by using Gauss' theorem on the term $\mavg{\dcell}{\macroexpl}$ from equation \eqref{eq:continuous-system-explicit-terms}: 
\begin{align*}
	\fluxrg &= \frac{\abs{\dface}}{\abs{\dcell}} \mapproxavg{\dface}{\ints{\vcoord \sprec{\PP}} \cdot \outernormal_{\dcellidx, \dcellidxnb}} \\
	&= \frac{1}{\abs{\dcell}} \sum_{\pcellidx \in \madjacents{\dcellidx}{\dcellidxnb}} \abs{\dfacet} \mapproxavg{\dfacet}{\ints{\vcoord \sprec{\PP}} \cdot \outernormal_{\dcellidx, \dcellidxnb}^{\pcellidx}} \\
	&\overset{(SO_1)}{=} \frac{1}{\abs{\dcell}} \sum_{\pcellidx \in \madjacents{\dcellidx}{\dcellidxnb}} \left( \abs{\dfacetone} \ints{\vcoord \pPP} \cdot \outernormal_{\dcellidx, \dcellidxnb}^{\pcellidx,1}  + \abs{\dfacettwo} \ints{\vcoord \pPP} \cdot \outernormal_{\dcellidx, \dcellidxnb}^{\pcellidx,2} \right)\\
	&\overset{(P)}{=} \frac{1}{\abs{\dcell}} \sum_{\pcellidx \in \madjacents{\dcellidx}{\dcellidxnb}} \abs{\dfacet} \ints{\vcoord \pPP} \cdot \outernormal_{\dcellidx, \dcellidxnb}^{\pcellidx} \\
\end{align*}
together with a quadrature rule $\spquad$. 
The unit outer normal of a facet $\dfacet$ is $\outernormal_{\dcellidx, \dcellidxnb}^{\pcellidx}$. 
The reconstruction $\sprec{\PP}(\xcoord)$ is a function that is piecewise continuous on primal cells and interpolates the averages: $\mavg{\pcell}{\sprec{\PP}} = \pPP$.  
In the first-order scheme the reconstruction is piecewise constant and equal to the cell average: 
\begin{align*}
	\evalat{\sprec{\PP}(\xcoord)}{\pcell} &= \pPP. 
\end{align*}
In the second-order scheme we make a piecewise linear ansatz
\begin{align*}
\evalat{\sprec{\PP(\xcoord)}}{\pcell} = \pPP + b \cdot (\xcoord-x_{\pcellidx}),
\end{align*}
for the reconstruction, where $b$ is a limited estimate of the slope that can be obtained by a minmod or WENO ansatz. 
Because we compute the flux on dual faces which are inside the primal cells where $\sprec{\PP}$ is continuous, we do not need an approximate flux function and only have to approximate the integrals by some quadrature. 
Using a piecewise constant reconstruction, these simplify to a single evaluation of the cell mean. 

Next we consider the fluxes $\fluxgr$ resulting from $\mavg{\pcell}{\microexplfibres}$ in \eqref{eq:continuous-system-explicit-terms}: 
\begin{align*}
\fluxgr &= \frac{\abs{\pface}}{\abs{\pcell}} \mapproxavg{\pface}{(\vcoord \sprec{\PM} \FD)  \cdot \outernormal_{\pcellidx, \pcellidxnb}} \\
&= \frac{1}{\abs{\pcell}}  \left( \sum_{\dcellidx \in \madjacents{\pcellidx}{\pcellidxnb}} \abs{\pfacet}\mapproxavg{\pfacet}{(\vcoord \sprec{\PM} \FD)} \right) \cdot \outernormal_{\pcellidx, \pcellidxnb} \\
&\overset{(SO_1)}{=} \frac{1}{\abs{\pcell}}  \left( \sum_{\dcellidx \in \madjacents{\pcellidx}{\pcellidxnb}} \abs{\pfacet} \vcoord \dPM \FD_{\pcellidx} \right) \cdot \outernormal_{\pcellidx, \pcellidxnb} \\
\end{align*}
This time, the facets $\pfacet$ which are parts of the primal face $\pface$ all share the same constant normal $\outernormal_{\pcellidx, \pcellidxnb}$.
$\sprec{\PM}(\xcoord)$ is a piecewise continuous reconstruction of the density on dual cells. 

Finally, application of the divergence theorem to $\mavg{\pcell}{\microexpl}$ in equation \eqref{eq:continuous-system-explicit-terms}, together with the projection 
\begin{align*}
	\Ptbproj (\divx (\vcoord \PP)) &= \divx (\vcoord \PP) - \divx \ints{\vcoord \PP} \FD 
\end{align*}
gives: 
\begin{align*}
	\fluxgg &= \frac{\abs{\pface}}{\abs{\pcell}} \left( \mapproxavg{\pface}{(\widehat{\vcoord \PP})} - \mapproxavg{\pface}{\widehat{\ints{\vcoord \PP}} \FD_{\pcell} } \right) \cdot \outernormal_{\pcellidx, \pcellidxnb}. 
\end{align*}
Here, $\widehat{\vcoord \PP}$ is an approximate flux function, for example the upwind flux, that depends on the left and right state $\sprec{\PP}_{\pcell}, \sprec{\PP}_{\pcellnb}$ of the face $\pface$. 
The second term of the projection is not in conservation form. In the spirit of wave-propagation for heterogeneous media as proposed by LeVeque (\cite{leveque2002finite}), we simply evaluate the equilibrium function $\FD_{\pcell}$ on the current cell $\pcell$.  

The approximate implicit operators are
\begin{align*}
\macroimpl_{\dcellidx} &= 0 \\
\microimpl_{\pcellidx} &= \frac{\pardel \factorcol[\dif, \pcellidx]}{\pareps^2} \LCO[\dif] \pPP = \mavg{\pcell}{\microimpl} + \orderof{\gridsize^2}.  
\end{align*}
If $\factorcol[\dif](\xcoord)$ is a constant on each cell, this is even exact, because $\LCO[\dif]$ is linear. 
Note that the implicit operator on a cell only depends on the cell mean. 
Thus the implicit part can be solved on each cell separately. 
This is still true for the \ivscheme{1} and \ivscheme{2} variants in which all of the volume terms are treated implicitly.

\subsection{The resulting scheme on a square grid}
\label{sec:ap-scheme-square-grid}

We consider the tensor-product grid defined by a list of nodes $\left(\xcoord_{\dimidx, 1}, \dots, \xcoord_{\dimidx, \midx^{\max}_{\dimidx}} \right)$ for each space dimension $\dimidx \in 1,\dots, \spacedim$.
Let $\multiidx = (\midx_1, ..., \midx_\spacedim)$ be a multi-index. 
The vertices of the tensor-product grid are all the points $\xcoord_{\multiidx} = (\xcoord_{1, \midx_1},\dots, \xcoord_{\spacedim, \midx_{\spacedim}})$ such that $1 \leq \midx_{\dimidx} \leq \midx^{\max}_{\dimidx}$. 
The primal cells $\mcell{\multiidx+\onehalf}$ of this grid are the boxes $\cuboid_{\spacedim} (\xcoord_{\multiidx}, \xcoord_{\multiidx + 1})$ with centers $\xcoord_{\multiidx + \onehalf} := \frac{\xcoord_{\multiidx} + \xcoord_{\multiidx + 1}}{2}$. 
The box spanned by the two points $\xcoord_{low}, \xcoord_{up}$ is defined as 
\begin{align*}
\cuboid_{\spacedim}(\xcoord_{low},\xcoord_{up}) = \left \lbrace \xcoord \in \R^{\spacedim}: \abs{\xcoord_{low}}_{\infty} \leq \abs{\xcoord}_{\infty} \leq \abs{\xcoord_{up}}_{\infty} \right \rbrace.
\end{align*}
With a slight abuse of multi-index notation, the sum of a multi-index and a scalar as in $\multiidx + 1 := (\midx_1 + 1, \dots, \midx_{\spacedim} + 1)$ is applied component-wise. 
The dual cell $\mcell{\multiidx}$ with center $\xcoord_{\multiidx}$ is the box $\cuboid_{\spacedim}(\xcoord_{\multiidx - \onehalf}, \xcoord_{\multiidx + \onehalf})$. 
In the following we show the \mischeme{1} scheme on a two-dimensional square-grid, i.e. a tensor-product grid where all nodes are equally spaced:
\begin{align*}
	\xcoord_{\multiidx} = \dnodeidx \gridsize.
\end{align*}
In the first-order \mischeme{1} scheme, the reconstructions $\sprec{\PM}, \sprec{\PP}$ are piecewise constant and equal to the cell means. 
All occurrences of a quadrature rule $\spquad$ are replaced by the midpoint-rule. 
Then the right-hand side of the macro equation becomes
\begin{align*}
	\macroexpl_{\dnodeidx} = -\pardel \frac{1}{2\gridsize} &\left\langle  -\vcoord_{\xcx}(\PP_{\dnodell} + \PP_{\dnodeul}) -\vcoord_{\xcy} (\PP_{\dnodell} + \PP_{\dnodeur} ) \right. \\
                                                          &\left.  +\vcoord_{\xcx}(\PP_{\dnodelr} + \PP_{\dnodeur}) +\vcoord_{\xcy} (\PP_{\dnodeur} + \PP_{\dnodeul} \right\rangle + \parthet \factorprol(\PM_{\dnodeidx})\PM_{\dnodeidx},
\end{align*}
when we insert the fluxes on all four faces. 
The term $\microexplfibres$ is 
\begin{align*}
\microexplfibres_{\pnodeidx} = -\frac{\pardel}{\pareps^2} \frac{1}{2\gridsize} \FD_{\pnodeidx} 
&\left[ -\vcoord_{\xcx} (\PM_{\pnodell} + \PM_{\pnodeul}) - \vcoord_{\xcy} (\PM_{\pnodell} + \PM_{\pnodelr}) \right. \\ 
&\left. +\vcoord_{\xcx} (\PM_{\pnodelr} + \PM_{\pnodeur}) + \vcoord_{\xcy} (\PM_{\pnodeur} + \PM_{\pnodeul}) \right]
\end{align*}
And finally we have: 
\begin{align*}
	\microexpl_{\pnodeidx} = 
	  &-\frac{\pardel}{\pareps}\frac{1}{2\gridsize} \left[ \widehat{-\vcoord_{\xcx}\PP}(\PP_{\pnodeidx}, \PP_{(\xidx-\onehalf, \yidx+\onehalf)}) 
	                                                      \widehat{-\vcoord_{\xcy} \PP}(\PP_{\pnodeidx}, \PP_{(\xidx+\onehalf, \yidx-\onehalf)}) \right. \\
	                                               &\qquad\qquad\left. + \widehat{\vcoord_{\xcx}\PP}(\PP_{\pnodeidx}, \PP_{(\xidx+\threehalf, \yidx+\onehalf)}) 
	                                                      +  \widehat{\vcoord_{\xcy} \PP}(\PP_{\pnodeidx}, \PP_{(\xidx+\onehalf, \yidx+\threehalf)})  \right] \\
	  &+ \frac{\pardel \parnu \factorcol[\adv, \pcellidx]}{\pareps^2} \left[ \LCO[\adv]\left( \tilde\PM_{\pnodeidx} \FD_{\pnodeidx} + \pareps \PP_{\pnodeidx}  \right) \right] \\
      &+ \frac{\parthet \factorprol( \tilde \PM_{\pnodeidx} )}{\epsilon} \left[ \Source \left( \tilde \PM_{\pnodeidx} \FD_{\pnodeidx} + \pareps \PP_{\pnodeidx} \right)  -\tilde \PM_{\pnodeidx} \FD_{\pnodeidx} \right]
\end{align*}
with an average density $\tilde \PM_{\pnodeidx} = \frac{1}{4}(\PM_{\pnodell} + \PM_{\pnodelr} + \PM_{\pnodeur} + \PM_{\pnodeul} )$ over the primal cell $\mcell{\pnodeidx}$. 
The numerical flux function can be any of the usual methods, for example the upwind flux 
\begin{align*}
	\widehat{\vcoord_{\xcx} \PP}(\PP_{\pnodeidx}, \PP_{(\xidx+\threehalf, \yidx+\onehalf)}) = \max(\vcoord_{\xcx},0) \PP_{\pnodeidx} + \min(\vcoord_{\xcx}, 0) \PP_{(\xidx+\threehalf, \yidx+\onehalf)}.
\end{align*}

\subsection{Time discretization}
 \newcommand{\lcol}[1]{\begin{array}{l} #1  \end{array}}
 \newcommand{\rcol}[1]{\begin{array}{r} #1  \end{array}}
 
We use the IMEX schemes from \cite{Ascher1997}. The time-step size is denoted by $\timestep$. In the first-order scheme, the forward-backward Euler scheme is used. For the particular system \eqref{eq:semi-discrete-system}, this reads
{
\renewcommand{\arraycolsep}{1.5pt}
\renewcommand{\arraystretch}{1.2}
\begin{alignat*}{3}
	\rcol{\PMvec^{*} \\ \PPvec^{*}} &
	\lcol{= \\ =} &&
	\lcol{\PMvec^{\timeidx} + \timestep \macroexplvec(\PMvec^{\timeidx}, \PPvec^{\timeidx}) \\ \PPvec^{\timeidx} + \timestep \left(\microexplfibresvec(\PMvec^{\timeidx}, \PPvec^{\timeidx}) + \microexplvec(\PMvec^{\timeidx}, \PPvec^{\timeidx})  \right)} &
	\left. \lcol{~ \\ ~} \right\} &\, \text{explicit euler step} \\
	\rcol{\PMvec^{\timeidx+1} \\ \PPvec^{\timeidx+1} } &
	\lcol{= \\ =} &&
	\lcol{\PMvec^{*} + \timestep \macroimplvec(\PMvec^{\timeidx+1}, \PPvec^{*}) \\ \PPvec^{*} + \timestep \microimplvec(\PMvec^{*}, \PPvec^{\timeidx + 1}) } &
	\left. \lcol{~ \\ ~} \right\} &\, \lcol{\text{implicit solve} \\ \text{without coupling}}
\end{alignat*}
}
In the minimally implicit variant \mischeme[]{1} we have $\macroimplvec = 0$ and thus the implicit solve for density reduces to $\PMvec^{\timeidx + 1} = \PMvec^{*}$. 

Lemou and Mieussens proved that their scheme is stable under the time step restriction 
\begin{align}
	\label{eq:stability-time-step}
	\timestep \leq \frac{1}{2} \left( \timestep_{\text{micro}} + \timestep_{\text{macro} } \right). 
\end{align}
We do not try to prove a stability result, but all out computations indicate that this choice leads to a stable scheme. 
The microscopic time step restriction comes from the CFL condition in the discretization of the transport part and is given by
\begin{align*}
	\timestep_{\text{micro}} = \frac{1}{2} \frac{\gridsize}{\speed}. 
\end{align*}
On the macroscopic scale, the scheme must respect the stability condition of the diffusion approximation as well as  the CFL condition from advection:
\begin{align*}
	\timestep_{\text{macro}} = \max \left(\frac{\gridsize^2}{2 \norm{\Difftens}}, \frac{\gridsize}{2\norm{\drift}} \right). 
\end{align*}

\begin{remark}[Glioma equation]
	Considering the glioma equation \eqref{eq:lke-glioma-scaled}, the implicit part in the \mischeme[]{1} scheme can be solved analytically. 
	We have
	\begin{align*}
		\pPP^{\timeidx+1} &=  \pPP^{*} + \timestep \frac{\pardel \factorcol[\dif, \pcellidx]}{\pareps^2}\LCO[\dif] \pPP^{\timeidx+1}\\
		&= \pPP^{*} - \timestep \frac{\pardel \factorcol[\dif, \pcellidx]}{\pareps^2} \pPP^{\timeidx+1} \\
	\end{align*}
	which is easily solved for the update: 
	\begin{align*}
		\pPP^{\timeidx+1} &= \frac{1}{1 + \timestep \frac{\pardel \factorcol[\dif, \pcellidx]}{\pareps^2} } \pPP^{*}
	\end{align*}
	This is of course no longer possible for the schemes \ivscheme{1} and \ivscheme{2}  with implicitly discretized volume terms. 
\end{remark}

The second-order scheme has to be chosen carefully to keep the asymptotic preserving property. The subclass of stiffly accurate schemes in \cite{Ascher1997}, in which the updated solution is identical to the last stage in a time-step, seems to maintain the AP-property. 
The second-order time-stepping scheme for \eqref{eq:semi-discrete-system} is
 {
 \renewcommand{\arraycolsep}{1.5pt}
 \renewcommand{\arraystretch}{1.2}
\begin{alignat*}{3}
	\rcol{(\macroexplvec)^{(1)} \\ (\microexplfibresvec)^{(1)} \\ (\microexplvec)^{(1)}} & 
	\lcol{= \\ = \\ = } &&
	\lcol{\macroexplvec(\PMvec^{\timeidx}, \PPvec^{\timeidx}) \\ \microexplfibresvec(\PMvec^{\timeidx}, \PPvec^{\timeidx}) \\ \microexplvec(\PMvec^{\timeidx}, \PPvec^{\timeidx})  } & 
	\left.\lcol{~\\~\\~}\right\}& \, \lcol{\text{compute operators} \\ \text{at time $\tcoord$} }\\
	\rcol{\PMvec^{*} \\ \PPvec^{*}} &
	\lcol{= \\ =} &&
	\lcol{\PMvec^{\timeidx} + \timetwostepfraction \timestep (\macroexplvec)^{(1)} \\ \PPvec^{\timeidx} + \timetwostepfraction \timestep \left((\microexplfibresvec)^{(1)} + (\microexplvec)^{(1)} \right) } &
	\left.\lcol{~\\~}\right\}& \, \lcol{\text{intermediate explicit} \\  \text{step to } \tcoord+\timetwostepfraction \timestep } \\
	\rcol{\PMvec^{(\timeidx, 1)} \\ \PPvec^{(\timeidx,1)} } &
	\lcol{= \\ = } &&
	\lcol{\PMvec^{*} + \timetwostepfraction \timestep \macroimplvec(\PMvec^{(\timeidx,1)}, \PPvec^{*})\\ \PPvec^{*} + \timetwostepfraction \timestep \microimplvec(\PMvec^{*}, \PPvec^{(\timeidx, 1)})} &
	\left.\lcol{~ \\ ~}\right\}& \,  \lcol{\text{intermediate implicit} \\ \text{step}}\\
	\rcol{(\macroexplvec)^{(2)} \\  (\microexplfibresvec)^{(2)} \\ (\microexplvec)^{(2)}  \\ (\macroimplvec)^{(2)} \\ (\microimplvec)^{(2)} } &
	\lcol{= \\ = \\ = \\ = \\ = } && 
	\lcol{\macroexplvec(\PMvec^{(\timeidx,1)},\PPvec^{(\timeidx,1)}) \\ \microexplfibresvec(\PMvec^{(\timeidx,1)}, \PPvec^{(\timeidx,1)}) \\ \microexplvec(\PMvec^{(\timeidx,1)}, \PPvec^{(\timeidx,1)})  \\ \macroimplvec(\PMvec^{(\timeidx,1)}, \PPvec^{(\timeidx,1)}) \\ \microimplvec(\PMvec^{(\timeidx,1)}, \PPvec^{(\timeidx,1)}) } &
	\left.\lcol{~\\~\\~\\~\\~}\right\}&\,  \lcol{\text{compute operators} \\ \text{at time $\tcoord+\timetwostepfraction \timestep$}} \\
	\rcol{\PMvec^{**} \\ ~ \\ \PPvec^{**}\\ ~ } &
	\lcol{= \\ ~ \\ = \\ ~} &&
	\lcol{\PMvec^{\timeidx} +  (1-\timetwostepfraction)\timestep (\macroimplvec)^{(2)}  \\  +  \timestep (\timetwostepweight  (\macroexplvec)^{(1)} + (1-\timetwostepweight) (\macroexplvec)^{(2)} ) \\ \PPvec^{\timeidx} + ( 1 - \timetwostepfraction ) \timestep (\microimplvec)^{(2)} \\ + \timestep (\timetwostepweight  (\microexplfibresvec + \microexplvec)^{(1)} + (1-\timetwostepweight) (\microexplfibresvec + \microexplvec)^{(2)} )} &
	\left.\lcol{~\\~\\~\\~}\right\}&\, \text{explicit step to $\tcoord + \timestep$} \\
	\rcol{\PMvec^{\timeidx+1} \\ \PPvec^{\timeidx+1} } &
	\lcol{= \\ = } &&
	\lcol{\PMvec^{**} + \timetwostepfraction \timestep \macroimplvec(\PMvec^{\timeidx+1}, \PPvec^{**})\\ \PPvec^{**} + \timetwostepfraction \timestep \microimplvec(\PMvec^{**}, \PPvec^{\timeidx+1})} &
	\left.\lcol{~ \\ ~}\right\}& \,  \text{implicit step}
\end{alignat*}
with the constants 
\begin{align*}
	\timetwostepfraction &= \frac{2 - \sqrt{2}}{2} \\
	\timetwostepweight &= 1 - \frac{1}{2 \timetwostepfraction}. 
\end{align*}
}
Our numerical experiments indicate that the time step \eqref{eq:stability-time-step} needs to be restricted further by a factor of $0.2$ to achieve stability with this scheme. 

\subsection{The asymptotic limit of the scheme}
\label{sec:asymptotic-limit}
We consider the first-order minimally implicit variant which can, with some reordering of the steps, be written as 
\begin{align*}
\PPvec^{*} &= \PPvec^{\timeidx} + \timestep \left(\microexplfibresvec(\PMvec^{\timeidx}, \PPvec^{\timeidx}) + \microexplvec(\PMvec^{\timeidx}, \PPvec^{\timeidx}) \right) \\
\PPvec^{\timeidx + 1} &= \PPvec^{*} + \timestep \microimplvec(\PPvec^{\timeidx + 1})  \\
\PMvec^{\timeidx+1} &= \PMvec^{\timeidx} + \timestep \macroexplvec(\PMvec^{\timeidx}, \PPvec^{\timeidx+1}). 
\end{align*}
This looks already like a discrete version of the derivation of the diffusion limit \eqref{eq:diffusion-limit-general} where we first computed the perturbation and then inserted this into the density equation. 
In the diffusion limit, only those terms with an $\frac{1}{\pareps^2}$ in front remain. 
Thus the implicit perturbation update reduces to 
\begin{align*}
	\PP_{\pcellidx}^{\timeidx + 1} &= -\frac{\pareps^2}{\timestep \pardel \factorcol[\dif, \pcellidx]} \LCO[\dif]^{-1} \PP_{\pcellidx}^{*} 
\end{align*}
with 
\begin{align*}
	\PP_{\pcellidx}^{*} &= \timestep  \left( \microexplfibres_{\pcellidx}(\PMvec^{\timeidx}, \PPvec^{\timeidx}) + \microexpl_{\pcellidx}(\PMvec^{\timeidx}, \PPvec^{\timeidx})\right) \\
                     	&= \timestep \left(
	                     	 -\frac{\pardel}{\pareps^2} \frac{1}{\abs{\pcell}} \sum_{\pcellidxnb \in \mneighbors{\pcellidx}} \left( \sum_{\dcellidx  \in \madjacents{\pcellidx}{\pcellidxnb}} \abs{\pfacet} \vcoord \PM_{\dcellidx}^{\timeidx} \FD_{\pcellidx} \right) \cdot \outernormal_{\pcellidx, \pcellidxnb}
	                     	 + \frac{\pardel \parnu \factorcol[\adv]}{\pareps^2} \LCO[\adv]\left(\FD_{\pcellidx}\right) \tilde{\PM}_{\pcellidx}^{\timeidx}
                     	 \right).
\end{align*}
Combining these two expressions yields
\begin{align*}
	 \PP_{\pcellidx}^{\timeidx + 1} &= -\frac{1}{\factorcol[\dif, \pcellidx]} \left( 
	 	-\frac{1}{\abs{\pcell}} \sum_{\pcellidxnb \in \mneighbors{\pcellidx}} \left( \sum_{\dcellidx  \in \madjacents{\pcellidx}{\pcellidxnb}} \abs{\pfacet} \LCO[\dif]^{-1} (\vcoord \FD_{\pcellidx})  \PM_{\dcellidx}^{\timeidx}  \right) \cdot \outernormal_{\pcellidx, \pcellidxnb}
	 	+\parnu \factorcol[\adv, \pcellidx]  \LCO[\dif]^{-1} \LCO[\adv]\left(\FD_{\pcellidx} \right) \tilde{\PM}_{\pcellidx}^{\timeidx}
	 	\right).
\end{align*}
Finally, we get the limit of the scheme as $\pareps \rightarrow 0$, when we insert this expression into the update for the density: 
\begin{align*}
	\PM_{\dcellidx}^{\timeidx+1} &= \PM_{\dcellidx}^{\timeidx} + \timestep \left( -\pardel \frac{1}{\abs{\dcell}} \sum_{\dcellidxnb \in \mneighbors{\dcellidx}} \sum_{\pcellidx \in \madjacents{\dcellidx}{\dcellidxnb}} \abs{\dfacet} \ints{\vcoord \PP_{\pcellidx}^{\timeidx+1} } \cdot \outernormal_{\dcellidx, \dcellidxnb}^{\pcellidx} + \parthet \factorprol(\PM_{\dcellidx}^{\timeidx}) \PM_{\dcellidx}^{\timeidx} \right).
\end{align*}
This is an explicit scheme for the density $\PM_{\dcellidx}^{\timeidx + 1}$. 
The updated value $\PM_{\dcellidx}^{\timeidx + 1}$ only depends on the previous values on the same dual cell $\dcell$ and those cells $\mcell{\dcellidxnb'}$ which are connected to it with at least a vertex, i.e., 
$\dcell \isect \mcell{\dcellidxnb'} \neq \emptyset$ or $\exists \dvertidx : \dvert \in \mverticesof{\dcellidx} \wedge \mverticesof{\dcellidxnb'} $.  

On a square grid in two dimensions, this is equivalent to 
\begin{align*}
	\PP_{\pnodeidx}^{\timeidx+1} &= \frac{1}{\factorcol[\dif, \pcellidx]} \frac{1}{2\gridsize} \LCO[\dif]^{-1} \left(\FD_{\pnodeidx} \right.
		\left. \left[ -\vcoord_{\xcx} \left(\PM_{\pnodeul}^{\timeidx} + \PM_{\pnodell}^{\timeidx}\right) -\vcoord_{\xcy} \left(\PM_{\pnodell}^{\timeidx} + \PM_{\pnodelr}^{\timeidx} \right) \right. \right. \\
		& \qquad\qquad\qquad\qquad\qquad\qquad  \left. \left. + \vcoord_{\xcx} \left(\PM_{\pnodelr}^{\timeidx} + \PM_{\pnodeur}^{\timeidx}\right) + \vcoord_{\xcy} \left(\PM_{\pnodeur}^{\timeidx} + \PM_{\pnodeul}^{\timeidx} \right)\right]  \right) \\
		&\quad -\frac{\parnu \factorcol[\adv, \pcellidx]}{\factorcol[\dif]} \LCO[\dif]^{-1} \LCO[\adv](\FD_{\pnodeidx}) \frac{1}{4}\left(\PM_{\pnodell}^{\timeidx}  + \PM_{\pnodelr}^{\timeidx}  + \PM_{\pnodeur}^{\timeidx}  + \PM_{\pnodeul}^{\timeidx} \right)
\end{align*}
\begin{align*}
	\PM_{\dnodeidx}^{\timeidx + 1} = \PM_{\dnodeidx}^{\timeidx} - \frac{\timestep \pardel}{2\gridsize}  
	&\left\langle  -\vcoord_{\xcx}\left(\PP_{\dnodell}^{\timeidx+1}  + \PP_{\dnodeul}^{\timeidx+1}\right) -\vcoord_{\xcy} \left(\PP_{\dnodell}^{\timeidx+1} + \PP_{\dnodeur}^{\timeidx+1} \right) \right. \\
	& \left.  +\vcoord_{\xcx}\left(\PP_{\dnodelr}^{\timeidx+1} + \PP_{\dnodeur}^{\timeidx+1}\right) +\vcoord_{\xcy} \left(\PP_{\dnodeur}^{\timeidx+1} + \PP_{\dnodeul}^{\timeidx+1}\right) \right\rangle + \timestep \parthet \factorprol(\PM_{\dnodeidx}^{\timeidx})\PM_{\dnodeidx}^{\timeidx},
\end{align*}
For the special case that the equilibrium $\FD$ and the factors $\factorcol[\dif], \factorcol[\adv]$ are constant in space, we write the resulting scheme as one equation for the density by eliminating the perturbations. 
After tedious calculations, we arrive at 
\begin{align*}
	\PM_{\dnodeidx}^{\timeidx + 1} = \PM_{\dnodeidx}^{\timeidx} + \timestep\frac{\pardel}{\factorcol[\dif]}\left( \overline{\divx (\Difftens \gradx\PM)} \right) - \timestep\frac{\pardel \parnu \factorcol[\adv]}{\factorcol[\dif]} \left(\overline{\divx(\drift \PM)} \right) + \timestep \parthet \factorprol(\PM_{\dnodeidx}^{\timeidx})\PM_{\dnodeidx}^{\timeidx}
\end{align*}
with approximations to the diffusion 
\begin{align*}
	\overline{\divx (\Difftens \gradx\PM)} = \frac{1}{4 \gridsize^2}\Big(
	&\PM_{(\xidx  , \yidx  )}^{\timeidx} (-4\Difftens_{\xcx \xcx} -4 \Difftens_{\xcy \xcy}) \\
	&+\PM_{(\xidx-1, \yidx  )}^{\timeidx} (2\Difftens_{\xcx \xcx} - 2\Difftens_{\xcy \xcy}) \\	
	&+\PM_{(\xidx+1, \yidx  )}^{\timeidx} (2\Difftens_{\xcx \xcx} -2\Difftens_{\xcy \xcy}) \\
	&+\PM_{(\xidx  , \yidx-1)}^{\timeidx} (-2\Difftens_{\xcx \xcx} + 2\Difftens_{\xcy \xcy}) \\
	&+\PM_{(\xidx  , \yidx+1)}^{\timeidx} (-2\Difftens_{\xcx \xcx} + 2\Difftens_{\xcy \xcy}) \\
	&+\PM_{(\xidx-1, \yidx-1)}^{\timeidx} (\Difftens_{\xcx \xcx} + 2\Difftens_{\xcx \xcy} + \Difftens_{\xcy \xcy}) \\
	&+\PM_{(\xidx+1, \yidx-1)}^{\timeidx} (\Difftens_{\xcx \xcx} - 2\Difftens_{\xcx \xcy} + \Difftens_{\xcy \xcy}) \\
	&+\PM_{(\xidx-1, \yidx+1)}^{\timeidx} (\Difftens_{\xcx \xcx} - 2\Difftens_{\xcx \xcy} + \Difftens_{\xcy \xcy}) \\
	&+\PM_{(\xidx+1, \yidx+1)}^{\timeidx} (\Difftens_{\xcx \xcx} + 2\Difftens_{\xcx \xcy} + \Difftens_{\xcy \xcy}) \Big)
\end{align*}
and drift 
\begin{align*}
	\overline{\divx(\drift \PM)} = \frac{1}{8 \gridsize} \Big( 
	&\PM_{(\xidx-1, \yidx  )}^{\timeidx} (-2\drift_{\xcx}             ) \\
	&+\PM_{(\xidx+1, \yidx  )}^{\timeidx} ( 2\drift_{\xcx}             ) \\
	&+\PM_{(\xidx  , \yidx-1)}^{\timeidx} (             -2\drift_{\xcy}) \\
	&+\PM_{(\xidx  , \yidx+1)}^{\timeidx} (              2\drift_{\xcy}) \\
	&+\PM_{(\xidx-1, \yidx-1)}^{\timeidx} (- \drift_{\xcx} - \drift_{\xcy}) \\
	&+\PM_{(\xidx+1, \yidx-1)}^{\timeidx} (  \drift_{\xcx} - \drift_{\xcy}) \\	
	&+\PM_{(\xidx-1, \yidx+1)}^{\timeidx} (- \drift_{\xcx} + \drift_{\xcy}) \\
	&+\PM_{(\xidx+1, \yidx+1)}^{\timeidx} (  \drift_{\xcx} + \drift_{\xcy}) 	\Big)
\end{align*}
wherein $\Difftens$ is the diffusion tensor from \eqref{eq:general-diffusion-tensor} and $\drift$ is the drift vector from \eqref{eq:general-drift-vector}.
If the diffusion tensor is the identity $\Difftens = \Identity$, which is the case for example in the glioma equation with isotropic equilibrium $\FD(\vcoord) = 1$, then the discrete diffusion reduces to a diagonal five-point stencil: 
\begin{align*}
\overline{\divx (\Difftens \gradx\PM)} = \overline{\divx (\gradx\PM)} = \frac{1}{2 \gridsize^2} \left( - 4 \PM_{(\xidx  , \yidx  )}^{\timeidx} +
\PM_{(\xidx-1, \yidx-1)}^{\timeidx} + \PM_{(\xidx+1, \yidx-1)}^{\timeidx} + \PM_{(\xidx-1, \yidx+1)}^{\timeidx} + \PM_{(\xidx+1, \yidx+1)}^{\timeidx} \right).
\end{align*}
In this special case, the presented AP-method is identical to the nodal scheme proposed in \cite{buet2012design}. 
As already discussed therein, the scheme leads to a decoupling of meshes. 
If we start with a Dirac initial condition on cell $(\xidx, \yidx)$, only every other cell $(\xidx + \xidx', \yidx + \yidx')$ with $\xidx' + \yidx' = 2q$ will ever receive some mass. 
Computations of this linesource test show a strong checkerboard pattern, see \figref{fig:spmc-a}. 

The drift is approximated by a central scheme, which is also not ideal. 
For example, inserting the first unit vector $\drift = (1, 0)\trans$ for the drift, we get
\begin{align*}
	\overline{\divx (\drift \PM)} = \overline{\partial_{\xcx} \PM} = \frac{1}{4 \gridsize} \left( -2 \PM_{(\xidx-1, \yidx)}^{\timeidx}  +2\PM_{(\xidx+1, \yidx)}^{\timeidx} - \PM_{(\xidx-1, \yidx-1)}^{\timeidx} + \PM_{(\xidx+1, \yidx-1)}^{\timeidx} - \PM_{(\xidx-1, \yidx+1)}^{\timeidx} + \PM_{(\xidx+1, \yidx+1)}^{\timeidx}\right).
\end{align*}

In the next two subsections we show how to modify the AP-method in such a way that the diffusion and drift are better approximated in the limit. 
Particularly, on a tensor-product grid the diffusion will be approximated by a standard five-point stencil, and the drift by an upwind method. 

\subsection{An improved diffusion stencil in the limit}
\label{sec:spmc}
In the last section we have seen that the numerical diffusion approximation results from a concatenation of the macroscopic fluxes $\timestep \macroexpl_{\dcellidx}(\PMvec^{\timeidx}, \PPvec^{\timeidx+1})$ with $-\frac{\pareps^2}{\timestep \pardel \factorcol[\dif]} \LCO[\dif]^{-1} {\microexplfibres}_{\pcellidx}(\PMvec^{\timeidx}, \PPvec^{\timeidx})$ on overlapping primal cells $\pcellidx \in \mverticesof{\dcellidx}$.
The goal of this section is to modify $\macroexpl_{\dcellidx}$ and $\microexplfibres_{\pcellidx}$ such that---on a square grid in two dimensions---the resulting diffusion approximation becomes the standard five-point stencil. 
To simplify the following computations as much as possible, we set $\pardel = 1$,  $\factorcol[\dif] = 1$ and use a constant-in-space equilibrium $\FD(\xcoord, \vcoord) = \FD(\vcoord)$ such that the diffusion tensor is $\Difftens = \Identity$. 

Recall the flux over primal faces in the most general form:
\begin{align*}
	\fluxgr = \frac{\abs{\pface}}{\abs{\pcell}} \mapproxavg{\pface}{(\vcoord \sprec{\PM} \FD) \cdot \outernormal_{\pcellidx, \pcellidxnb}}. 
\end{align*}
Together with a piecewise constant reconstruction of the density $\evalat{\sprec{\PM}}{\dcell} = \dPM$ this results in the formulation 
\begin{align*}
	\fluxgr = \frac{1}{\abs{\pcell}}  \left( \sum_{\dcellidx \in \madjacents{\pcellidx}{\pcellidxnb}} \abs{\pfacet} \vcoord \dPM \FD_{\pcellidx} \right) \cdot \outernormal_{\pcellidx, \pcellidxnb}.
\end{align*}
This is a sum of constant fluxes over the facets $\pfacet$, weighted by the facet volumes $\abs{\pfacet}$. 
In the derivation of the AP scheme on square grids in \secref{sec:ap-scheme-square-grid} we used this method. 
Considering the primal face $(\xidx+\onehalf, \yidx+\onehalf), (\xidx+\threehalf, \yidx+\onehalf)$ in effect this method assigns equal weights $\frac{1}{2 \gridsize}$ to both overlapping dual cells $(\xidx+1, \yidx)$, $(\xidx+1, \yidx+1)$.
We get the same weights if we reconstruct $\sprec{\PM}$ as a globally continuous function from bilinear elements on each dual cell and use a midpoint quadrature rule on the faces. 
Starting from this interpretation, we define four variants of the microscopic flux 
\begin{align*}
\microexplfibressk{\pnodeidx}{\xcx,+}, \microexplfibressk{\pnodeidx}{\xcx,-}, \microexplfibressk{\pnodeidx}{\xcy,+}, \microexplfibressk{\pnodeidx}{\xcy,-}
\end{align*} 
that use different quadratures on different faces. 
In the $(\xcx,+)$-variant, the flux on $\xcx$-normal faces is evaluated at the upmost points, but for the $\xcy$-normal faces the midpoint rule is used:  
\begin{align*}
	\microexplfibressk{\pnodeidx}{\xcx,+} = -\frac{\pardel}{\pareps^2}\frac{1}{\gridsize} \FD_{\pnodeidx} \left[ -\vcoord_{\xcx} \PM_{(\xidx, \yidx+1)} - \onehalf\vcoord_{\xcy} (\PM_{(\xidx,\yidx)}+ \PM_{(\xidx+1, \yidx)}) + \vcoord_{\xcx} \PM_{(\xidx+1, \yidx+1)} + \onehalf\vcoord_{\xcy}(\PM_{(\xidx, \yidx+1)} + \PM_{(\xidx+1,\yidx+1)})\right]. 
\end{align*}
Similarly the $(\xcx,-)$-variant uses evaluations at the lowest points in $\xcx$-normal faces: 
\begin{align*}
\microexplfibressk{\pnodeidx}{\xcx,-} = -\frac{\pardel}{\pareps^2}\frac{1}{\gridsize} \FD_{\pnodeidx} \left[ -\vcoord_{\xcx} \PM_{(\xidx, \yidx)} - \onehalf\vcoord_{\xcy} (\PM_{(\xidx,\yidx)} +\PM_{(\xidx+1, \yidx)}) + \vcoord_{\xcx} \PM_{(\xidx+1, \yidx)} + \onehalf\vcoord_{\xcy}(\PM_{(\xidx, \yidx+1)} + \PM_{(\xidx+1,\yidx+1)})\right]. 
\end{align*}
The other variants are defined analogously for the $\xcy$-normal faces. 
The shifted evaluations are zeroth-order accurate Gauss-Radau quadrature rules, which is sufficient for a first-order scheme. 
In a second-order scheme, they have to be replaced by the correct first-order Gauss-Radau rules. 
We use each flux variant in the perturbation update $\PP_{\pnodeidx}^{\timeidx+1}$ in turn to compute the four modified perturbations 
\begin{align*}
	\PP_{\pnodeidx}^{\timeidx+1, (\xcx,+)}, \PP_{\pnodeidx}^{\timeidx+1, (\xcx,-)}, \PP_{\pnodeidx}^{\timeidx+1, (\xcy,-)}, \PP_{\pnodeidx}^{\timeidx+1, (\xcy,-)}.
\end{align*}

Now we modify the density flux $\macroexpl_{\dnodeidx}$. 
In each flux over a dual facet, the 'correct' variant of the perturbation is used:
\begin{align*}
\PM_{\dnodeidx}^{\timeidx + 1} = \PM_{\dnodeidx}^{\timeidx} - \frac{\timestep \pardel}{2\gridsize}  
&\left\langle  -\vcoord_{\xcx}\left(\PP_{\dnodell}^{\timeidx+1, (\xcx,+)}  + \PP_{\dnodeul}^{\timeidx+1, (\xcx,-)}\right) -\vcoord_{\xcy} \left(\PP_{\dnodell}^{\timeidx+1, (\xcy,-)} + \PP_{\dnodeur}^{\timeidx+1, (\xcy,-)} \right) \right. \\
& \left.  +\vcoord_{\xcx}\left(\PP_{\dnodelr}^{\timeidx+1, (\xcx,+)} + \PP_{\dnodeur}^{\timeidx+1, (\xcx,-)}\right) +\vcoord_{\xcy} \left(\PP_{\dnodeur}^{\timeidx+1, (\xcx,-)} + \PP_{\dnodeul}^{\timeidx+1, (\xcx,+)}\right) \right\rangle + \timestep \parthet \factorprol(\PM_{\dnodeidx}^{\timeidx})\PM_{\dnodeidx}^{\timeidx},
\end{align*}
The same tedious calculations as in the previous \secref{sec:asymptotic-limit} show that the diffusion is approximated by
\begin{align*}
\overline{\divx (\Difftens \gradx\PM)} = \frac{1}{4 \gridsize^2}\Big(
&\PM_{(\xidx  , \yidx  )}^{\timeidx} (-8\Difftens_{\xcx \xcx} -8 \Difftens_{\xcy \xcy}) \\
&+\PM_{(\xidx-1, \yidx  )}^{\timeidx} (4\Difftens_{\xcx \xcx}) 
+\PM_{(\xidx+1, \yidx  )}^{\timeidx} (4\Difftens_{\xcx \xcx}) 
+\PM_{(\xidx  , \yidx-1)}^{\timeidx} (4\Difftens_{\xcy \xcy}) 
+\PM_{(\xidx  , \yidx+1)}^{\timeidx} (4\Difftens_{\xcy \xcy}) \\
&+\PM_{(\xidx-1, \yidx-1)}^{\timeidx} (2\Difftens_{\xcx \xcy}) 
+\PM_{(\xidx+1, \yidx-1)}^{\timeidx} (-2\Difftens_{\xcx \xcy}) 
+\PM_{(\xidx-1, \yidx+1)}^{\timeidx} (-2\Difftens_{\xcx \xcy}) 
+\PM_{(\xidx+1, \yidx+1)}^{\timeidx} (2\Difftens_{\xcx \xcy}) \Big)
\end{align*}
in the limit, which is the classical five-point stencil
\begin{align*}
\overline{\divx (\Difftens \gradx\PM)} = \overline{\divx (\gradx\PM)} = \frac{1}{2 \gridsize^2} \left( - 4 \PM_{(\xidx  , \yidx  )}^{\timeidx} +
\PM_{(\xidx-1, \yidx)}^{\timeidx} + \PM_{(\xidx+1, \yidx)}^{\timeidx} + \PM_{(\xidx, \yidx-1)}^{\timeidx} + \PM_{(\xidx, \yidx+1)}^{\timeidx} \right).
\end{align*}
if the diffusion tensor is isotropic $\Difftens = \Identity$. 

\begin{remark}[Extension to three dimensions]
	In three space dimensions the procedure is structurally very similar but the notation becomes even more unwieldy. 
	The computational cost also increases, because we need twelve variants, four for each normal direction. 
	For example, in the variant $(\xcx,++)$, the fluxes over $\xcx$-normal faces are evaluated at the top right node. 
\end{remark}

\begin{figure}[h]
	\def\localpath{figures/SpuriousModesCorrector/}
	\centering
	\withfiguresize{\figurethreecol}{\figurethreecol}{\externaltikz{spurious_modes_corrector}{\begin{tikzpicture}

\def\dlplot#1#2#3#4#5{
	\nextgroupplot[
		xmin=-1,
		xmax=1,
		ymin=-1,
		ymax=1,
		xlabel={$\xcx / mm$},
		ylabel={$\xcy / mm$},
		title={\tikztitle{#1}},
		point meta min ={#2},
		point meta max ={#3},
		colorbar horizontal,
		colorbar style={,
			at={(0,1.02)},
			anchor=south west,
			height=0.02\textwidth,
			width=\figurewidth,
			xticklabel style={font=\figureticklabelfont{\figurewidth},anchor=south},
			ticklabel pos=right,
			max space between ticks=1000pt,
			try min ticks=3,
		},
		#4
	]
	\addplot graphics[xmin=-1.5,xmax=1.5,ymin=-1.5,ymax=1.5] {\localpath #5};
}

\def\errplot#1#2#3#4#5{
	\nextgroupplot[
		xmin=-1,
		xmax=1,
		ymin=-1,
		ymax=1,
		xlabel={$\xcx / mm$},
		ylabel={$\xcy / mm$},
		title={\tikztitle{#1}},
		colorbar horizontal,
		colorbar style={
			colormap name ={cmRdBu_r_minus},
			point meta min ={-#3}, 
			point meta max ={0},
			at={(0,1.02)},
			anchor=south west,
			height=0.02\textwidth,
			width=0.49\figurewidth,
			ticklabel pos=right,
			xticklabel style={font=\figureticklabelfont{\figurewidth},anchor=south, /pgf/number format/.cd,	/tikz/.cd},
			xticklabel={
			  	\pgfmathsetmacro\upBound{0.7*\pgfkeysvalueof{/pgfplots/xmax}+0.3*\pgfkeysvalueof{/pgfplots/xmin}}
			  	\ifdim \tick pt < \upBound pt 
					 $-10^{\pgfmathparse{-\tick+#2}\pgfmathprintnumber\pgfmathresult}$ 
				 \fi
			 },
			try min ticks=2,
		},
		colorbar/draw/.append code={
			\begin{axis}[
				every colorbar, 
				at={(1., 1.02)},
				anchor=south east,
				colormap name={cmRdBu_r_plus},
				colorbar horizontal,
				point meta min={0},
				point meta max={#3},
				colorbar=false,
				height=0.02\textwidth,
				width=0.49\figurewidth,
				ticklabel pos=right,
				xticklabel style={font=\figureticklabelfont{\figurewidth},anchor=south, /pgf/number format/.cd, /tikz/.cd},
				try min ticks=2, 
				xticklabel={
					\pgfmathsetmacro\loBound{0.3*\pgfkeysvalueof{/pgfplots/xmax}+0.7*\pgfkeysvalueof{/pgfplots/xmin}}
						\ifdim \tick pt > \loBound pt 
							$10^{\pgfmathparse{\tick+#2}\pgfmathprintnumber\pgfmathresult}$ 
						\fi
				},
				extra x ticks=0., 
				extra x tick labels=$\pm10^{#2}$,
			] 
			\pgfkeysvalueof{/pgfplots/colorbar addplot}
			\end{axis}
		},		
		#4
	] 
	\addplot graphics[xmin=-0.25,xmax=0.25,ymin=-0.25,ymax=0.25] {\localpath #5};	
}

\begin{groupplot}[
	group style={group size=3 by 1, horizontal sep = \figurehorizontalsep,  vertical sep = 3cm},
	axis on top,
	scale only axis,
	enlargelimits=false,
	width = \figurewidth,
	height =\figureheight,
	colormap name ={cmviridis},  
	scaled x ticks=false,
   	title style = {font=\figurelabelfont{\figurewidth}, yshift = 0.6cm},
   	xlabel style = {font=\figurelabelfont{\figurewidth}}, 
   	ylabel style = {font=\figurelabelfont{\figurewidth}}, 
   	xticklabel style={font=\figureticklabelfont{\figurewidth}},
   	yticklabel style={font=\figureticklabelfont{\figurewidth}},
]

	\dlplot{$\PM_{\times}$}            {-0.0025}{0.0025}{}{SpuriousModesCorrectorOff.png}
	\dlplot{$\PM_{+}$}                 {-0.0025}{0.0025}{}{SpuriousModesCorrectorOn.png}
	\errplot{$\reldiff{\PM_{\times}}{\PM_{+}}$}{-4}{5}{xmin=-0.25,xmax=0.25,ymin=-0.25,ymax=0.25}{SpuriousModesCorrectorDiff.png}
\end{groupplot}

 \end{tikzpicture}}}
	\settikzlabel{fig:spmc-a}
	\settikzlabel{fig:spmc-b}
	\settikzlabel{fig:spmc-c}
	\caption{Comparison between the direct application of the scheme $\mischeme{1}$ and the scheme with improved diffusion stencil \mischeme[+]{1} from \secref{sec:spmc} on the linesource benchmark. 
		 Plots of the density for \mischeme[\times]{1} (\ref{fig:spmc-a}) and \mischeme[+]{1} (\ref{fig:spmc-b}). In \ref{fig:spmc-c}, the relative difference  $\reldiff{\PM_{\times}}{\PM_{+}} = \frac{1}{\max\abs{\PM_{+}}} (\PM_{\times} - \PM_{+})$ is plotted on a signed truncated logarithmic scale.}
	\label{fig:spmc}
\end{figure}
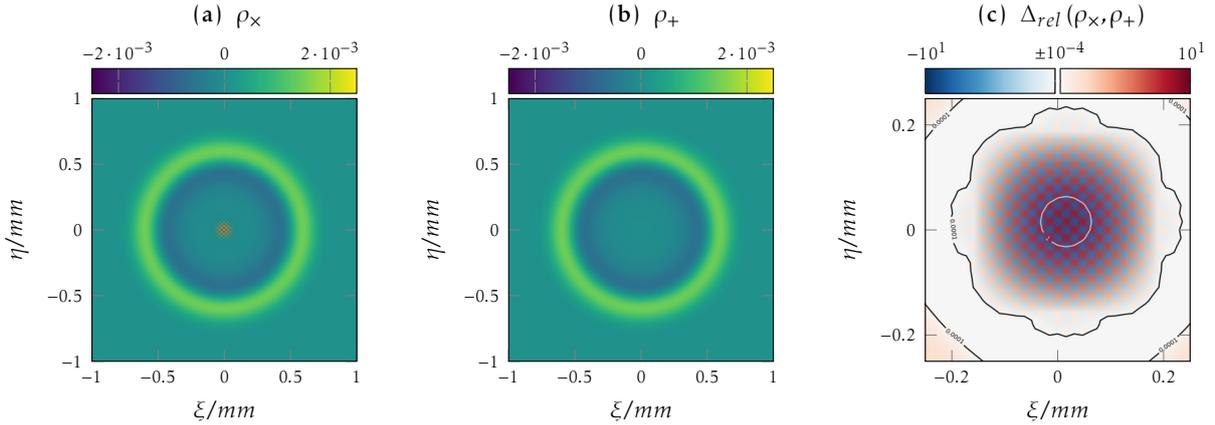

\subsection{Upwind discretization of the drift in the limit}
\label{sec:upwind-drift}
The limit drift approximation follows from a concatenation of the macroscopic flux $\timestep \macroexpl$ with $-\frac{\parnu \factorcol[\adv]}{\factorcol[\dif]} \LCO[\dif]^{-1} \LCO[\adv](\FD_{\pcellidx}) \tilde{\PM}_{\pcellidx}^{\timeidx}$. 
Using an average density $\tilde{\PM}_{\pcellidx}$ weighted by the subcell volumes as in \eqref{eq:density-on-primal-cell} leads to a central approximation of the drift. 
But we know the local drift direction 
\begin{align*}
	\drift_{\pcellidx} = \ints{\vcoord \LCO[\dif]^{-1}\LCO[\adv]\FD_{\pcellidx}}
\end{align*}
and want to assign more weight to those cells $\dcell$ that are upwind of the center $\dvert$. 
We write $\dvert^{*}$ for the intersection of the ray 
\begin{align*}
	\dvert - \tau \drift_{\pcellidx}, \quad  \tau \in \R^+
\end{align*}
with the cell boundary $\partial \pcell$.
Then we define 
\begin{align*}
	\tilde{\PM}_{\pcellidx} = \sprec{\PM}(\dvert^{*})
\end{align*}
with a continuous, piecewise linear reconstruction $\sprec{\PM}$ by hat-functions. 
This is of course only a first-order accurate approximation of the drift.

\subsection{Treatment of boundary conditions}
We consider only boundary conditions that preserve mass. 
On a macroscopic level this translates to a zero-flux Robin-type boundary condition for the density in \eqref{eq:diffusion-limit-general}: 
\begin{align}
	\label{eq:bc-limit-zero-flux}
	\evalat{-\divx (\density \Difftens) + \parnu \density \drift}{\partial \xdomain} = 0. 
\end{align}
This does not determine the boundary conditions on the microscopic level uniquely. 
All microscopic boundary conditions for $\PD$ that can be cast into the class of reflective boundary conditions preserve mass. 
At a reflective boundary, the values $\PD(\vcoord)$ are prescribed for incoming velocities $\vcoord \cdot \outernormal < 0$ and follow from the outgoing values via the reflection integral:
\begin{align}
	\label{eq:bc-full}
	\PD(\vcoord) &= \intBplus{\vcoord'}{\boundarykernel(\vcoord, \vcoord') \PD(\vcoord')} & \forall \Vminus{\vcoord}, 
\end{align}
Of course, the reflection kernel $\boundarykernel$ is defined such that the net mass flux across the boundary is zero, that is, it fulfills
\begin{equation}
	\label{eq:bc-zero-flux}
	\begin{aligned}
	0 = \intV{(\vcoord \cdot \outernormal) \PD(\vcoord)}
	&= \intBplus{\vcoord}{(\vcoord \cdot \outernormal) \PD(\vcoord)} + \intBminus{\vcoord}{(\vcoord \cdot \outernormal) \PD(\vcoord)}\\
	&= \intBplus{\vcoord}{(\vcoord \cdot \outernormal) \PD(\vcoord)} + \intBminus{\vcoord}{(\vcoord \cdot \outernormal) \intBplus{\vcoord'}{\boundarykernel(\vcoord, \vcoord') \PD(\vcoord')}}. 
	\end{aligned}
\end{equation} 
From the last line, we see that this is the case if 
\begin{align*}
		\intBminus{\vcoord}{(\vcoord \cdot \outernormal) \boundarykernel(\vcoord, \vcoord')} = -\vcoord' \cdot \outernormal 
\end{align*}
holds.
To see the boundary condition for $\PP$ that is equivalent to \eqref{eq:bc-full}, we insert the micro-macro decomposition \eqref{eq:APsplitPD} and obtain  
\begin{align*}
	\PP(\vcoord) &= \frac{\PM}{\pareps} \left[ \intBplus{\vcoord'}{\boundarykernel(\vcoord, \vcoord') \FD(\vcoord')} - \FD(\vcoord) \right] + \intBplus{\vcoord'}{\boundarykernel(\vcoord, \vcoord') \PP(\vcoord')}
\end{align*}
If the kernel is not compatible with the equilibrium state then in the limit when $\pareps$ tends to zero,  $\PP$ becomes unbounded at the boundary and we need to solve a half-space problem to compute the boundary condition. 
Here we do not want to consider boundary layers and therefore demand that \eqref{eq:bc-full} hold for the equilibrium state $\FD$. 
Then we have the condition 
\begin{align}
\label{eq:bc-micro}
\PP(\vcoord) &= \intBplus{\vcoord'}{\boundarykernel(\vcoord, \vcoord') \PP(\vcoord')}
\end{align}
for $\PP$. 
The value for $\PM$ is left unconstrained. 

For the kernel, we consider two options. 
The 'u-turn' kernel models that cells turn around 180 degrees when encountering a wall, independent of the angle of collision. 
It is given by
\begin{align*}
	\boundarykernel_{\text{uturn}}(\vcoord, \vcoord') &= \dirac{\vcoord}{-\vcoord'}. 
\end{align*}
Because the equilibrium fulfills $\FD(\vcoord) = \FD(-\vcoord)$, the reflection equation \eqref{eq:bc-full} holds for the equilibrium.
It is easy to check the zero-mass-flux condition \eqref{eq:bc-zero-flux} for this kernel. 

Another option is that after a collision with the wall, the incoming particles are in equilibrium
\begin{align*}
	\PD(\vcoord) &= \alpha \FD(\vcoord) & \forall \Vminus{\vcoord}. 
\end{align*}
This so-called thermal boundary condition can be achieved with the kernel 
\begin{align*}
	\boundarykernel_{\text{thermal}}(\vcoord, \vcoord') &= \frac{\alpha \FD(\vcoord)}{\intBplus{\vcoord'}{\PD(\vcoord')}}. 
\end{align*}
The parameter $\alpha$ is defined by
\begin{align*}
	\alpha &= -\frac{\intBplus{\vcoord}{(\vcoord \cdot \outernormal) \PD(\vcoord)}}{\intBminus{\vcoord}{(\vcoord \cdot \outernormal) \FD(\vcoord)}}
\end{align*}
to fulfill the zero-mass-flux condition \eqref{eq:bc-zero-flux}. 
For a symmetric equilibrium we have $\alpha = 1$ and thus the boundary condition is compatible with the equilibrium.

\begin{remark}[Specular reflection]
	The specular reflection kernel 
	\begin{align*}
		\boundarykernel_{\text{spec}}(\vcoord, \vcoord') &= \dirac{\vcoord}{\vcoord' - 2 (\vcoord' \cdot \outernormal) \outernormal}
	\end{align*}
	models hard-sphere collisions between particles and the wall. 
	It conserves mass, but is not compatible with the equilibrium in general, only if the equilibrium is mirror symmetric around the outer boundary
	\begin{align*}
		\FD(\vcoord) &= \FD(\vcoord - 2 (\vcoord \cdot \outernormal) \outernormal). 
	\end{align*}
\end{remark}

If we want to, we can additionally constrain the density
\begin{align*}
	\evalat{\density}{\partial \xdomain} = \density_b. 
\end{align*}
Then, together with \eqref{eq:bc-limit-zero-flux} this implies a condition for $\gradx \density$, which can always be fulfilled because $\Difftens$ is invertible. 
On the particle level, this means that we get the additional condition 
\begin{align*}
	\evalat{\PM}{\partial \xdomain} = \intV{\PD(\vcoord)} = \density_b. 
\end{align*}

\section{Discretization of the velocity space by a linear spectral method}
\label{sec:spectral-method}
The scheme that we derived in the previous sections is discrete in time and space. 
It remains to find a suitable discretization for the velocity. 
We use a linear spectral Galerkin method based on real-valued spherical harmonics, which is a slight modification of the well-known $P_N$ method \cite{brunner2005two,seibold2014starmap,garrett2013comparison}. 
First we define the spherical harmonics basis for the full space $\ltwospace(\US)$ of particle distributions $\PD$. 
A basis for the constrained space of perturbations $\PP$ from \lemref{lem:lcol-properties}
\begin{align}
	\label{eq:perturbed-space}
	\PP \in V := \Nsporth(\LCO[\dif]) := \left\{ \PP \in \weightedltwospace , (\PP,\FD)_{\FD} = \ints{\PP} = 0 \right\}, 
\end{align}  
is then obtained by removing the first element in the full basis. 

We collect the $2l+1$ harmonics of exactly order $l$ in the vector $\basisPDord{l}$. 
For example, there is one zeroth-order harmonic $\basisPDord{0} = \frac{1}{\sqrt{4\pi}}$, and there are three first-order harmonics $\basisPDord{1} = \frac{1}{\sqrt{4\pi}}(\sqrt{3} \vcoord_{\xcx}, \sqrt{3} \vcoord_{\xcy}, \sqrt{3} \vcoord_{\xcz})$.
For an exact definition of the real-valued spherical harmonics refer to \cite{seibold2014starmap}.
The $(N+1)^2$ spherical harmonics up to order $N$
\begin{align*}
\basisPD = \left(\basisPDord{0}, \basisPDord{1}, \dots, \basisPDord{N}\right) = \left(\basisPDcomp{0}, \dots, \basisPDcomp{n-1} \right), 
\end{align*}  
span a finite-dimensional subspace of $\ltwospace(\US)$---the space of polynomials up to order $N$.
An infinite-dimensional basis of the full space $\ltwospace(\US)$ is given by 
\begin{align*}
	\basisPD^{\infty} = \left(\basisPDord{0}, \basisPDord{1}, \dots \right). 
\end{align*}
One important property of the spherical harmonics is that they are orthonormal, that is
\begin{align*}
	\ints{\basisPDcomp{i} \basisPDcomp{j}} = \delta_{ij} 
\end{align*}
holds for any $i,j$. 
Thus all basis components except for $\basisPDcomp{0}$ fulfill the constraint $\ints{\PP} = 0$ in \eqref{eq:perturbed-space}: 
\begin{align*}
	\ints{\basisPDcomp{i}} = \frac{1}{\sqrt{4 \pi}} \ints{\basisPDcomp{i} \basisPDcomp{0}} = 0 \qquad i > 0. 
\end{align*}
We obtain bases for the constrained space $V$, and corresponding finite-dimensional subspaces by omitting the function $\basisPDcomp{0}$: 
\begin{align*}
	\basisPP^{\infty} &= \left(\basisPPord{1}, \basisPPord{2}, \dots\right) := \left(\basisPDord{1}, \basisPDord{2}, \dots\right), \\
	\basisPP &= \left(\basisPPord{1}, \dots, \basisPPord{N} \right).
\end{align*}
A perturbation $\PP$ has a unique basis representation
\begin{align*}
	\PP(\vcoord) &= \momPP^{\infty} \cdot \basisPP^{\infty}(\vcoord), 
\end{align*}
wherein the coefficients $\momPPcomp{i}^{\infty}$ are equal to the moments 
\begin{align*}
	\ints{\PP \basisPPcomp{i}^{\infty}} = \ints{\sum_{j} \momPPcomp{j}^{\infty} \basisPPcomp{j}^{\infty} \basisPPcomp{i}^{\infty} } = \momPPcomp{i}^{\infty}. 
\end{align*}
because of the orthonormal property. 
The orthogonal projection of $\PP$ onto the finite-dimensional subspace $V_h$ is 
\begin{align*}
	\recPP(\vcoord) &= \momPP \cdot \basisPP(\vcoord), 
\end{align*}
with moments
\begin{align*}
	\momPPcomp{i} = \ints{\recPP \basisPP} = \begin{cases}
	\momPPcomp{i}^{\infty} = \ints{\PP \basisPP} &, i < n \\
	0 &, i \geq n.
	\end{cases}
\end{align*}
The discrete-in-velocity approximation of problem \eqref{eq:continuous-system} is to find $(\PM, \recPP)$ that solve 
\begin{equation}
\label{eq:velocity-discrete-system}
\begin{alignedat}{2}
\dt \PM &= \macroexpl(\PM, \recPP) &+& \macroimpl(\PM,\recPP), \\
\dt \ints{\recPP \basisPP} = \dt \momPP  &= \ints{\microexplfibres(\PM) \basisPP} + \ints{\microexpl(\PM, \recPP) \basisPP} &+& \ints{\microimpl(\PM,\recPP) \basisPP}. 
\end{alignedat}
\end{equation}
This is a set of $n+1 = (N+1)^2$ equations for the $n+1$ unknowns $(\PM, \momPP)$. 
The individual terms therein are 
\begin{align*}
\macroexpl(\PM,\recPP)   &= -\pardel \divx \ints{\vcoord \recPP} + \parthet \factorprol \PM,\\
\ints{\microexplfibres(\PM) \basisPP} &= -\frac{\pardel}{\pareps^2} \divx \left(\PM \ints{\vcoord\FD\basisPP}\right), \\
\ints{\microexpl(\PM, \recPP) \basisPP}  &=  -\frac{\pardel}{\pareps} \left[\divx \ints{\vcoord \recPP \basisPP} - \divx \ints{\vcoord \recPP} \ints{\FD \basisPP} \right] + \frac{\pardel \parnu \factorcol[\adv]}{\pareps^2} \ints{\LCO[\adv] (\PM \FD + \pareps \recPP) \basisPP} + \frac{\parthet \factorprol}{\pareps} \ints{ \Ptbproj \Source (\PM \FD + \pareps \recPP) \basisPP}, 
\end{align*}
and
\begin{align*}
\macroimpl(\PM,\recPP)   &= 0, \\
\ints{ \microimpl(\PM, \recPP) \basisPP}  &= \frac{\pardel \factorcol[\dif]}{\pareps^2} \ints{ \LCO[\dif](\recPP) \basisPP}.
\end{align*}
The equations are coupled through the flux moments $\ints{\vcoord \recPP} \in \R^{\spacedim}$, $\ints{\vcoord \recPP \basisPP} \in \R^{n \times \spacedim}$ and moments of the collision term and source on the right hand side. 
The macro equation is coupled with the micro equations through the moments
\begin{align*}
	\ints{\vcoord \recPP} = \frac{\sqrt{4\pi}}{\sqrt{3}}\ints{\basisPPord{1} \recPP} = \frac{\sqrt{4\pi}}{\sqrt{3}} \momPPord{1}. 
\end{align*}
In general, $i$-th order flux moments $\ints{\vcoord \recPP \basisPPord{i}}$ can be written as a combination of the moments $\ints{\recPP \basisPP^{(i+1)}} = \momPP^{(i+1)}$ of order $i+1$.
Usually this relation is written in matrix form.
For instance for the $\xcx$-component of the velocity, we write
\begin{align*}
	\ints{\vcoord_{\xcx} \recPP \basisPP} &= M_{\xcx} \momPP \\
	&:= \ints{\vcoord_{\xcx} \basisPP \basisPP\trans} \momPP
\end{align*}
For details on how to compute these matrices for the full basis $\basisPD$, see for example \cite{seibold2014starmap}. 
Due to orthogonality of the basis, we can simply remove the first row and column of the matrix $\ints{\vcoord_{\xcx} \basisPD \basisPD\trans}$ to get the matrices for the restricted basis $\basisPP$. 
Because the turning operators are linear, we can also write their contribution to the moment system in matrix form: 
\begin{align*}
	\ints{\LCO[\dif](\recPP) \basisPP} &= C_{\dif} \momPP, \\
	\ints{\LCO[\adv](\recPP) \basisPP} &= C_{\adv} \momPP.  
\end{align*} 

\begin{remark}[Turning operators in the glioma equation]
	From equation \eqref{eq:glioma-lco-inverse} we have 
	\begin{align*}
		\LCO[\dif](\recPP) = -\recPP, 
	\end{align*}
	thus 
	\begin{align*}
		\ints{\LCO[\dif](\recPP) \basisPP} = -\ints{\recPP \basisPP} = -\momPP, 
	\end{align*}
	and $C_{\dif} = -\Identity$. 	
	The turning perturbation is given by 
	\begin{align*}
		\LCO[\adv](\recPP) &= \hat\lambda_H \gradx \VF \cdot (\FD \ints{\vcoord \recPP} - \vcoord \recPP). 
	\end{align*}
	Its moments are 
	\begin{align*}
		\ints{\LCO[\adv](\recPP) \basisPP} &= \hat\lambda_H \gradx\VF \cdot \left(\ints{\FD \basisPP} \ints{\vcoord \recPP} - \ints{\vcoord \recPP \basisPP} \right) 
	\end{align*}
	The dot product is between components of the gradient $\gradx\VF$ and components of the velocity $\vcoord$. 
	The moments appearing in this expression have been calculated before. 
	With some abuse of vector notation, we have 
	\begin{align*}
	\ints{\LCO[\adv](\recPP) \basisPP} = \hat\lambda_H \gradx\VF \cdot \left(\ints{\FD \basisPP} \frac{\sqrt{4\pi}}{\sqrt{3}} \momPPord{1} - (M_{\xcx} \momPP, M_{\xcy} \momPP, M_{\xcz} \momPP)  \right). 
	\end{align*}
	Because the source is just the identity $\Source \PD = \PD$, the source moments can be simplified to 
	\begin{align*}
		\ints{ \Ptbproj \Source (\PM \FD + \pareps \recPP) \basisPP} = \pareps \momPP. 
	\end{align*}
\end{remark}

\begin{remark}
	Equation \eqref{eq:velocity-discrete-system} is equivalent to the moment system 
	\begin{align*}
		\dt \momPD &= -\frac{\pardel}{\pareps} \divx \ints{\vcoord \recPD \momPD}  + \frac{\pardel}{\pareps^2}\factorcol[\dif] \ints{\LCO[\dif] \recPD \momPD} + \frac{\pardel \parnu}{\pareps} \factorcol[\adv] \ints{\LCO[\adv] \recPD \momPD} + \parthet \factorprol \ints{\Source \recPD \momPD}
	\end{align*}
	for the original equation \eqref{eq:lke-scaled} with the approximation $\recPD$ and moments $\momPD$ of the particle distribution $\PD$ given by
	\begin{align*}
	\recPD &= \momPD \cdot \basisPD = \PM \FD + \pareps \recPP , \\
	\momPDcomp{i} &= \ints{\recPD \basisPDcomp{i}} = \begin{cases}
	\frac{1}{\sqrt{4\pi}} \PM  & i = 0\\
	\PM \ints{\FD \basisPPcomp{i}} + \pareps \momPPcomp{i} & i > 0
	\end{cases}
	\end{align*}
\end{remark}

The space and time discretization can be carried over to the moment system without change.

\section{Results}
\label{sec:results}
Whenever we know the analytical solution to a problem, we use it to numerically evaluate the convergence of our code with respect to grid refinement.
One such convergence test consists of $L+1$ runs with identical parameters but increasing grid refinement, starting with $M_{0}$ grid points per space direction and increasing by a constant factor $r$ in each step. 
In run $l$, the number of grid points per dimension is then
\begin{align*}
	M_l &= \lfloor M_0 r^l \rfloor \qquad l = 0,\dots, L, 
\end{align*} 
and the size of each grid cell 
\begin{align*}
\gridsize_l = \frac{1}{M_l} = \frac{1}{\lfloor M_0 r^l \rfloor } \qquad l = 0,\dots,L. 
\end{align*}
The error $e_l$ in each run is defined as the $\ltwospace$-difference between the computed density $\PM_l$ and the exact solution $\PM_{ex}$, evaluated at the final time $T$
\begin{align*}
e_l = \Vert \PM_l(T, \xcoord) - \PM_{ex}(T, \xcoord) \Vert_{2} = \left( \int_{\xdomain} (\PM_l - \PM_{ex})^2 d\xcoord \right)^{\frac{1}{2}}. 
\end{align*}
The integral is computed by a quadrature of appropriate order. 
Convergence rates between successive refinement steps are computed with the formula
\begin{align*}
\frac{\log(e_l) - \log(e_{l+1})}{\log(\gridsize_l) - \log(\gridsize_{l+1})}. 
\end{align*}
In the presentation and discussion of results, we will make use of the pointwise relative difference 
\begin{align*}
\reldiff{f}{g}(\xcoord) = \frac{1}{\underset{\xcoord \in \xdomain}{\max}\abs{g}} (f(\xcoord) - g(\xcoord))
\end{align*}  between two functions $f(\xcoord), g(\xcoord)$. 
In error plots, a signed truncated logarithmic scale 
\begin{align*}
\text{sign}(f) \big(\log(\max(\abs{f}, f_L))  - \log(f_L)\big)
\end{align*}
is useful to show a wide range of absolute values as well as their signs. 

All computations are performed on the glioma model from \secref{sec:glioma-model} with the peanut distribution \eqref{eq:peanut}. 
When not otherwise mentioned, we use the minimally implicit scheme with the stencil improvements from \secref{sec:spmc} and \secref{sec:upwind-drift}. 
For the computations in \secref{sec:diffusion-limit-convergence} and \secref{sec:discontinuities} we need to prescribe the macroscopic diffusion tensor $\DT$. 
We achieve this by constructing artificial values for the water diffusion tensor
\begin{align*}
\DW = \frac{1}{2} \left(5 \DT - \Identity \right),
\end{align*}
according to the inverse of \eqref{eq:DT-glioma}. 
Whenever we prescribe the macroscopic drift $\driftT$, we define the volume fraction $\VF$ according to the inverse of \eqref{eq:drift-glioma}:
\begin{align*}
	\gradx \VF &= \frac{1}{\hat \lambda_H} \driftT\trans \DT^{-1}
\end{align*} 
When the physical values of $\DT, \driftT$ are given together with $\xref$ we can compute corresponding parameters $\speed, \lambda_0, \lambda_1$ for the microscopic glioma equation using the scaling relations in \secref{sec:parabolic-scaling}:
\begin{align*}
	\speed = \frac{\Difftens_0}{\xref \pareps}, \quad \lambda_0 = \frac{\Difftens_0}{\xref^2 \pareps^2}, \quad	\lambda_1 = \frac{\drift_0}{\xref \pareps^2}.
\end{align*}

\subsection{Fundamental solution of the limit equation}
\label{sec:diffusion-limit-convergence}
When the diffusion tensor $\Difftens$ and drift $\drift$ are constant and the growth factor $\parthet$ is zero, the limiting advection-diffusion equation in physical coordinates \eqref{eq:diffusion-limit-physical} has the fundamental solution 
\begin{align}
	\label{eq:fundamental-solution}
	\density_{,f} &= \left((4\pi)^{\spacedim} \det \Difftens\right)^{-\onehalf} \tcoord^{-\frac{\spacedim}{2}} \exp\left(-\frac{1}{4\tcoord} (\xcoord - \drift \tcoord)\trans \Difftens^{-1} (\xcoord - \drift \tcoord)\right). 
\end{align} 
Our scheme should reproduce this solution when $\pareps$ is small. 
For the test we choose
\begin{align*}
	\DT &= \Difftens_0 \frac{1}{4.5} R \begin{pmatrix}
	2.5 & 0 & 0 \\ 0 & 1 & 0 \\ 0 & 0 & 1
	\end{pmatrix} R\trans, \\
	\driftT &= \drift_0 \frac{1}{\sqrt{10}} \begin{pmatrix}
	3 \\ 1 \\ 0
	\end{pmatrix}.
\end{align*}
Herein the matrix $R$ rotates $e_1$ onto the main diffusion direction $(-1, 2, 0)\trans$. 
We choose a characteristic diffusion speed $\Difftens_0 = \frac{1}{100}$. 
We perform two tests, one without drift, i.e., $\drift_0 = 0$,  and one with drift speed $\drift_0 = 0.1$. 

To smoothen the initial Dirac-delta distribution, we choose the initial condition $\PM(0, \xcoord) = \density_{,f}(\tcoord_O, \xcoord)$ with the time offset $\tcoord_O = 0.2$. 
Then the solution at time $\tcoord$ is given by $\density_{,f}(\tcoord + \tcoord_O, \xcoord)$. 

First we test convergence of the first and second order schemes with respect to grid refinement, starting at a $40\times40$ grid and refining by factor $1.5$ five times. 
The analytical solution is of course only valid in the diffusion limit, therefore we choose $\pareps = 10^{-5}$. 
The $\ltwospace$ error over the number of grid points is plotted in \figref{fig:fundamental-convergence-grid}. 
Without the drift term, both schemes converge with second order accuracy to the analytic solution, as is to be expected for a discretization of the pure diffusion equation. 
With the drift, the order of both schemes is reduced to about $0.9$ and absolute errors are also much greater. 

We are also interested in convergence as $\pareps$ tends to zero. 
From the grid refinement study, we see that at about $200 \times 200$ grid points, the error is roughly $\scinum{2}{-5}$ without drift and $\scinum{4}{-4}$ with the drift term. 
As $\pareps$ approaches zero, we expect the total error to be dominated by this discretization error. 
In \figref{fig:fundamental-convergence-epsilon}, the $\ltwospace$ error of the first order scheme at $200 \times 200$ grid points is plotted, over values of $\pareps$ from one to $\scinum{1}{-9}$. 
We observe that the error levels out at the expected discretization error below a threshold value of $\pareps$---roughly $\scinum{1}{-4}$ without drift and $\scinum{1}{-3}$ with drift. 
Note that for certain intermediate values of $\pareps$, the error reaches a local minimum slightly below the limit discretization error because kinetic effects cancel out some of the numerical diffusion of the scheme. 
Numerical solutions in the kinetic to intermediate regime ($\pareps \in [0.1, 0.01]$) are shown in \figref{fig:fundamental-limit-kinetic}. 
In the kinetic regime, the problem is similar to the linesource problem \cite{garrett2013comparison}; only for anisotropic scattering. 
Indeed the $P_1$ solutions feature a single ellipsoid wave, which travels at speed $\frac{1}{\sqrt{3}}\speed$ in the main diffusion direction and is biased towards the drift direction. 
With decreasing $\pareps$ the diffusion dominates and the wave maximum is smeared out into a Gaussian. 
Below $\pareps \approx \scinum{1}{-2}$ the solutions are too similar for direct visual comparisons. 
Therefore, in \figref{fig:fundamental-limit-diffusive} we show relative differences on a signed logarithmic scale instead. 
\figref{fig:fundamental-limit-diffusive-a} to \figref{fig:fundamental-limit-diffusive-f} show relative differences between the numerical solution and the fundamental solution to the diffusion equation \eqref{eq:fundamental-solution} . 
Although not visible from a plot of the solution, at $\pareps = \scinum{1}{-2}$ still has some small kinetic effects(see \figref{fig:fundamental-limit-diffusive-a}) of relative magnitude $\scinum{1}{-2}$. 
In \figref{fig:fundamental-limit-diffusive-f} the relative difference between the numerical solutions at $\pareps = \scinum{1}{-3}$ and $\pareps = \scinum{1}{-9}$ is plotted.  
We see that already at $\pareps = \scinum{1}{-3}$ the discretization error dominates the kinetic effects. 
From \figref{fig:fundamental-limit-diffusive-b} we see how the kinetic effects cancel some of the numerical diffusion. 
The numerical diffusion from the drift discretization becomes apparent from \figref{fig:fundamental-limit-diffusive-e}: Looking in drift direction the solution at $\pareps = \scinum{1}{-9}$ overestimates the fundamental solution before and after the peak and underestimates at the peak. 

\begin{figure}[h]
	\def\localpath{figures/DiffusionLimit/}
	\centering
	\withfiguresize{\figurethreecol}{\figurethreecol}{\externaltikz{numerical_diffusion}{	
\begin{tikzpicture}
\def\signedlogoffset{-7}
\def\signedlog#1#2{%
	sign(#1) * (ln(abs(#1))/ln(10) - #2)%
}
\begin{groupplot}[
    group style={group size=3 by 1, horizontal sep = \figurehorizontalsep,  vertical sep = 2.5cm},
    width =\figurewidth,
    height =\figureheight,
   	title style = {font=\figurelabelfont{\figurewidth}},
   	xlabel style = {font=\figurelabelfont{\figurewidth}}, 
   	ylabel style = {font=\figurelabelfont{\figurewidth}}, 
   	xticklabel style={font=\figureticklabelfont{\figurewidth}},
   	yticklabel style={font=\figureticklabelfont{\figurewidth}},
   	scale only axis,
    xmode=log,
    xtick={40, 60, 90, 135, 203, 304},
    xlabel={\# gridpoints},
    xticklabel=\pgfmathparse{(1.5)^\tick)}\pgfmathprintnumber{\pgfmathresult},
    x tick label style={/pgf/number format/.cd, fixed, precision=0},
    log basis x=1.5, 
    ylabel style={yshift=0.5em},
   	yticklabel={
   		\ifdim \tick pt > 0 pt 
   		$+10^{\pgfmathparse{\tick+\signedlogoffset}\pgfmathprintnumber\pgfmathresult}$ 
   		\else
   		\ifdim \tick pt < 0 pt
   		$-10^{\pgfmathparse{-\tick+\signedlogoffset}\pgfmathprintnumber\pgfmathresult}$
   		\else
   		$\pm10^{\pgfmathparse{-\tick+\signedlogoffset}\pgfmathprintnumber\pgfmathresult}$
   		\fi
   		\fi
   	},
    extra y ticks={0},
    extra tick style={grid=major},
]
\nextgroupplot[
	title=\tikztitle{Larger eigenvalue},
	ylabel={$\lambda_1$},
	legend columns=4,
	legend entries={\mischeme[]{1} $\driftT = 0$, \mischeme[]{1} $\driftT = 0.1$, \mischeme[]{2} $\driftT = 0$, \mischeme[]{2} $\driftT = 0.1$, $\mathcal{O}(\Delta x)$, $\mathcal{O}(\Delta x^{2})$},
	legend to name={CommonLegend},
	legend style={font=\figurelabelfont{\figurewidth}}, 
	extra y tick labels={},
]
\addplot table [x=nx, y expr=(\signedlog{\thisrowno{1}}{\signedlogoffset}), col sep=comma] {\localpath numerical_diffusion_STO11_nodrift.csv};
\addplot table [x=nx, y expr=(\signedlog{\thisrowno{1}}{\signedlogoffset}), col sep=comma] {\localpath numerical_diffusion_STO11.csv};
\addplot table [x=nx, y expr=(\signedlog{\thisrowno{1}}{\signedlogoffset}), col sep=comma] {\localpath numerical_diffusion_STO22_nodrift.csv};
\addplot table [x=nx, y expr=(\signedlog{\thisrowno{1}}{\signedlogoffset}), col sep=comma] {\localpath numerical_diffusion_STO22.csv};
\addplot[domain=40:304, black, thick, mark=] {\signedlog{40 * 0.001*x^(-1)}{\signedlogoffset}}; 
\addplot[domain=40:304, black, thick, dashed, mark=] {\signedlog{40 * 40 * 0.0001*x^(-2)}{\signedlogoffset}}; 
\addplot[domain=40:304, black, thick, dashed, mark=] {\signedlog{-40 * 40 * 0.0001*x^(-2)}{\signedlogoffset}}; 

\nextgroupplot[
title=\tikztitle{Smaller eigenvalue},
ylabel={$\lambda_2$},
extra y tick labels={},
]
\addplot table [x=nx, y expr=(\signedlog{\thisrowno{2}}{\signedlogoffset}), col sep=comma] {\localpath numerical_diffusion_STO11_nodrift.csv};
\addplot table [x=nx, y expr=(\signedlog{\thisrowno{2}}{\signedlogoffset}), col sep=comma] {\localpath numerical_diffusion_STO11.csv};
\addplot table [x=nx, y expr=(\signedlog{\thisrowno{2}}{\signedlogoffset}), col sep=comma] {\localpath numerical_diffusion_STO22_nodrift.csv};
\addplot table [x=nx, y expr=(\signedlog{\thisrowno{2}}{\signedlogoffset}), col sep=comma] {\localpath numerical_diffusion_STO22.csv};
\addplot[domain=40:304, black, thick, dashed, mark=] {\signedlog{40 * 40 * 0.0001*x^(-2)}{\signedlogoffset}}; 
\addplot[domain=40:304, black, thick, dashed, mark=] {\signedlog{-40 * 40 * 0.0001*x^(-2)}{\signedlogoffset}}; 
\nextgroupplot[
title=\tikztitle{Diffusion angle},
ylabel={$\theta$ / deg},
ymode=linear,
ytick={-45, 0},
yticklabel=\pgfmathparse{\tick}\pgfmathprintnumber{\pgfmathresult},
y tick label style={/pgf/number format/.cd, fixed, precision=2},
extra y ticks={18.43494882292201, -63.43494882292201},
]
\addplot table [x=nx, y expr=(\thisrowno{3} - 90), col sep=comma] {\localpath numerical_diffusion_STO11_nodrift.csv};
\addplot table [x=nx, y=theta, col sep=comma] {\localpath numerical_diffusion_STO11.csv};
\addplot table [x=nx, y=theta, col sep=comma] {\localpath numerical_diffusion_STO22_nodrift.csv};
\addplot table [x=nx, y=theta, col sep=comma] {\localpath numerical_diffusion_STO22.csv};
\end{groupplot}

\path[] (group c1r1.north) -- node[anchor=south, yshift=1cm]{\ref{CommonLegend}} (group c3r1.north); 
\end{tikzpicture}}}
	\settikzlabel{fig:fundamental-numerical-diffusion-a}
	\settikzlabel{fig:fundamental-numerical-diffusion-b}
	\settikzlabel{fig:fundamental-numerical-diffusion-c}
	\caption{Estimates of the numerical diffusion of the \mischeme[]{1} and \mischeme[]{2} schemes using the fundamental solution. Shown are the larger eigenvalue of the numerical diffusion tensor in \ref{fig:fundamental-numerical-diffusion-a}, the smaller eigenvalue in \ref{fig:fundamental-numerical-diffusion-b} and the direction of the main eigenvector in \ref{fig:fundamental-numerical-diffusion-c} for each scheme without and with the drift term.}
	\label{fig:fundamental-numerical-diffusion}
\end{figure}
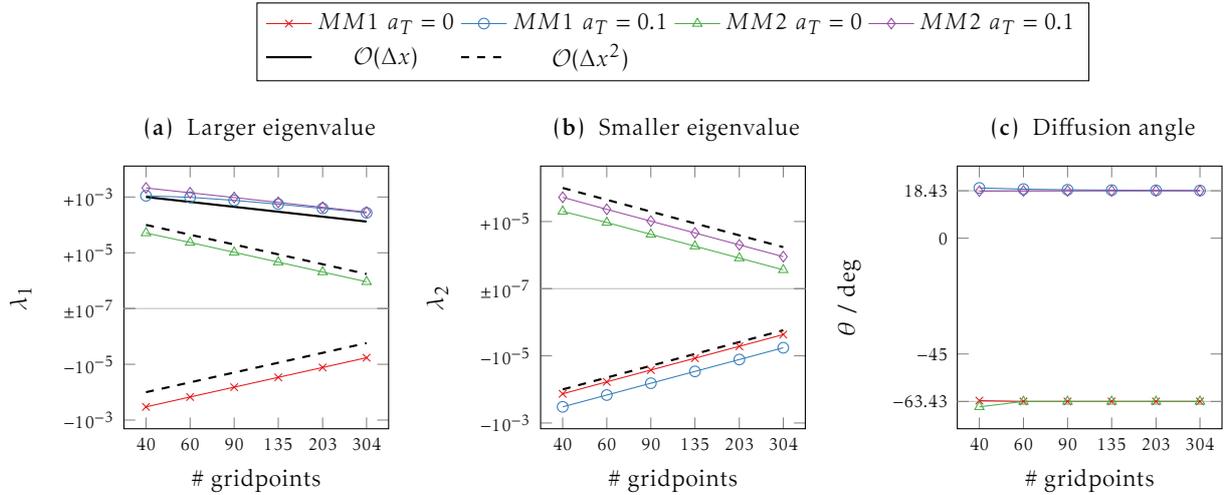

With the fundamental solution we can also quantify the numerical diffusion of the scheme. 
We fit a multivariate Gaussian to the numerical result and view the corresponding estimated diffusion tensor as the sum of the exact diffusion tensor and a contribution from the numerical scheme. 
In \figref{fig:fundamental-numerical-diffusion}, the two eigenvalues and the main direction of this estimated numerical diffusion are plotted. 
We observe that numerical diffusion converges at the same rate as the $\ltwospace$ error. 
When the drift term is active, it dominates the overall numerical diffusion by two orders of magnitude and the main axis of the numerical diffusion is parallel to the drift direction.
Without the drift, we observe an interesting difference between the \mischeme[]{1} scheme and the \mischeme[]{2} scheme. 
For the \mischeme[]{2} scheme, both eigenvalues are positive and their ratio is close to the anisotropy factor $2.5$.  
Additionally, the main axes of physical and numerical diffusion are aligned. 
Thus, the numerical diffusion is proportional to the physical diffusion. 
In the \mischeme[]{1} scheme the ratio of eigenvalues and main axis is the same.
However, both eigenvalues are negative, which indicates that the leading numerical error is dispersive rather than diffusive. 

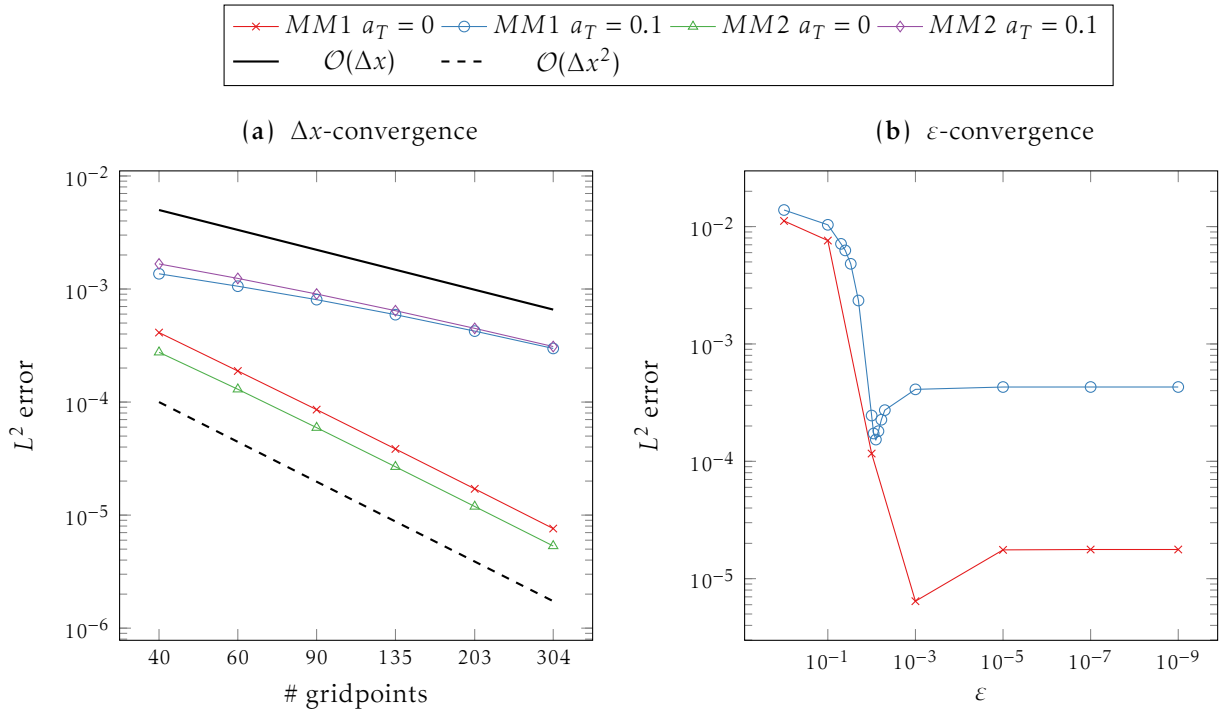
\begin{figure}[h]
	\def\localpath{figures/DiffusionLimit/}
	\centering
	\withfiguresize{\figuretwocol}{\figuretwocol}{\externaltikz{diffusion_limit_convergence}{\begin{tikzpicture}
\begin{groupplot}[
    group style={group size=2 by 2, horizontal sep = \figurehorizontalsep,  vertical sep = 2.5cm},
    scale only axis,
    width = \figurewidth,
    height = \figureheight,
   	title style = {font=\figurelabelfont{\figurewidth}},
   	xlabel style = {font=\figurelabelfont{\figurewidth}}, 
   	ylabel style = {font=\figurelabelfont{\figurewidth}}, 
	xticklabel style={font=\figureticklabelfont{\figurewidth}},
	yticklabel style={font=\figureticklabelfont{\figurewidth}},
    xmode=log,
    ymode=log,
    ylabel style={yshift=0.5em},
    ylabel={$\ltwospace$ error},
    log basis y=10,
]
\nextgroupplot[
	title=\tikztitle{$\Delta x$-convergence},
	log basis x=1.5, 
	xlabel={\# gridpoints},
	xtick={40, 60, 90, 135, 203, 304},
	xticklabel=\pgfmathparse{(1.5)^\tick)}\pgfmathprintnumber{\pgfmathresult},
	x tick label style={/pgf/number format/.cd, fixed, precision=0},
	legend columns=4,
	legend entries={\mischeme[]{1} $\driftT = 0$, \mischeme[]{1} $\driftT = 0.1$, \mischeme[]{2} $\driftT = 0$, \mischeme[]{2} $\driftT = 0.1$, $\mathcal{O}(\Delta x)$, $\mathcal{O}(\Delta x^{2})$},
	legend to name={DiffusionLimitLegend}, 
	legend style={font=\figurelabelfont{\figurewidth}},
]
\addplot table [x=nx, y=sto11_nodrift, col sep=comma] {\localpath convergence_space.csv};
\addplot table [x=nx, y=sto11,         col sep=comma] {\localpath convergence_space.csv};
\addplot table [x=nx, y=sto22_nodrift, col sep=comma] {\localpath convergence_space.csv};
\addplot table [x=nx, y=sto22,         col sep=comma] {\localpath convergence_space.csv};
\addplot[domain=40:304, black, thick, mark=] {40*0.005*x^(-1)}; 
\addplot[domain=40:304, black, thick, dashed, mark=] {40*40*0.0001*x^(-2)}; 

\nextgroupplot[
	title=\tikztitle{$\varepsilon$-convergence},
	log basis x=10, 
	xlabel={$\varepsilon$},
	x dir = reverse,
	]
\addplot table [x=epsilon, y=sto11_nodrift, col sep=comma] {\localpath convergence_epsilon.csv};
\addplot table [x=epsilon, y=sto11        , col sep=comma] {\localpath convergence_epsilon.csv};
\end{groupplot}
\path[] (group c1r1.north) -- node[anchor=south, yshift=1cm]{\ref{DiffusionLimitLegend}} (group c2r1.north); 
\end{tikzpicture}}}
	\settikzlabel{fig:fundamental-convergence-grid}
	\settikzlabel{fig:fundamental-convergence-epsilon}
	\caption{Convergence study for the fundamental solution test from \secref{sec:diffusion-limit-convergence}. \ref{fig:fundamental-convergence-grid}: $\ltwospace$ errors over number of gridpoints on each axis. Shown are the errors for both \mischeme[]{1}, and \mischeme[]{2}, each without drift $\driftT = 0$ and with some drift $\driftT = 0.1$.  \ref{fig:fundamental-convergence-epsilon}: $\ltwospace$ errors over parabolic scaling parameter $\pareps$, for a fixed grid with $200 \times 200$ cells. Errors for the \mischeme[]{1} are shown both without and with drift.}
	\label{fig:fundamental-convergence}
\end{figure}

\begin{figure}[h]
	\def\localpath{figures/DiffusionLimit/}
	\centering
	\withfiguresize{\figurethreecol}{\figurethreecol}{\externaltikz{diffusion_limit_kinetic}{\begin{tikzpicture}
\def\dlplot#1#2#3#4{
	\nextgroupplot[
		xmin=-0.5,
		xmax=0.5,
		ymin=-0.5,
		ymax=0.5,
		xlabel={$\xcx / mm$},
		ylabel={$\xcy / mm$},
		title={\tikztitle{#1}},
		point meta min ={#2},
		point meta max ={#3},
		colorbar horizontal,
		colorbar style={,
			at={(0,1.02)},
			anchor=south west,
			height=0.02\textwidth,
			width=\figurewidth,
			xticklabel style={font=\figureticklabelfont{\figurewidth},anchor=south},
			ticklabel pos=right,
			max space between ticks=1000pt,
			try min ticks=3,
		},
	]
	\addplot graphics[xmin=-1,xmax=1,ymin=-1,ymax=1] {\localpath #4};
}

\begin{groupplot}[
	group style={group size=3 by 3, horizontal sep = \figurehorizontalsep,  vertical sep = 3.cm},
	axis on top,
	scale only axis,
	enlargelimits=false,
	width = \figurewidth,
	height = \figureheight,
	colormap name ={cmviridis},  
	scaled x ticks=false,
   	title style = {font=\figurelabelfont{\figurewidth}, yshift = 0.6cm},
   	xlabel style = {font=\figurelabelfont{\figurewidth}}, 
   	ylabel style = {font=\figurelabelfont{\figurewidth}}, 
   	xticklabel style={font=\figureticklabelfont{\figurewidth}},
   	yticklabel style={font=\figureticklabelfont{\figurewidth}},	
]
	\dlplot{$\varepsilon = 10^{-1}$}        {-0.011}{ 0.1}{Fundamental_STO11_Eps1e-1.png}
	\dlplot{$\varepsilon = 5\times10^{-2}$} {-0.011}{ 0.1}{Fundamental_STO11_Eps5e-2.png}
	\dlplot{$\varepsilon = 4\times10^{-2}$} {-0.011}{ 0.1}{Fundamental_STO11_Eps4e-2.png}
	\dlplot{$\varepsilon = 3\times10^{-2}$} {-0.011}{ 0.1}{Fundamental_STO11_Eps3e-2.png}
	\dlplot{$\varepsilon = 2\times10^{-2}$} {-0.011}{ 0.1}{Fundamental_STO11_Eps2e-2.png}
	\dlplot{$\varepsilon = 10^{-2}$}        {-0.011}{ 0.1}{Fundamental_STO11_Eps1e-2.png}
\end{groupplot}

 \end{tikzpicture}}}
	\settikzlabel{fig:fundamental-limit-kinetic-a}
	\settikzlabel{fig:fundamental-limit-kinetic-b}
	\settikzlabel{fig:fundamental-limit-kinetic-c}
	\settikzlabel{fig:fundamental-limit-kinetic-d}
	\settikzlabel{fig:fundamental-limit-kinetic-e}
	\settikzlabel{fig:fundamental-limit-kinetic-f}
	\caption{The numerical solution to the fundamental solution test of \secref{sec:diffusion-limit-convergence}, using the \mischeme{1}-$P_1$ scheme. The density $\PM$ is depicted for solutions with various values of $\pareps$, ranging from the kinetic regime $\pareps = 0.1$ in \ref{fig:fundamental-limit-kinetic-a} to the intermediate regime in \ref{fig:fundamental-limit-kinetic-f} with $\pareps = 10^{-2}$.}
	\label{fig:fundamental-limit-kinetic}
\end{figure}
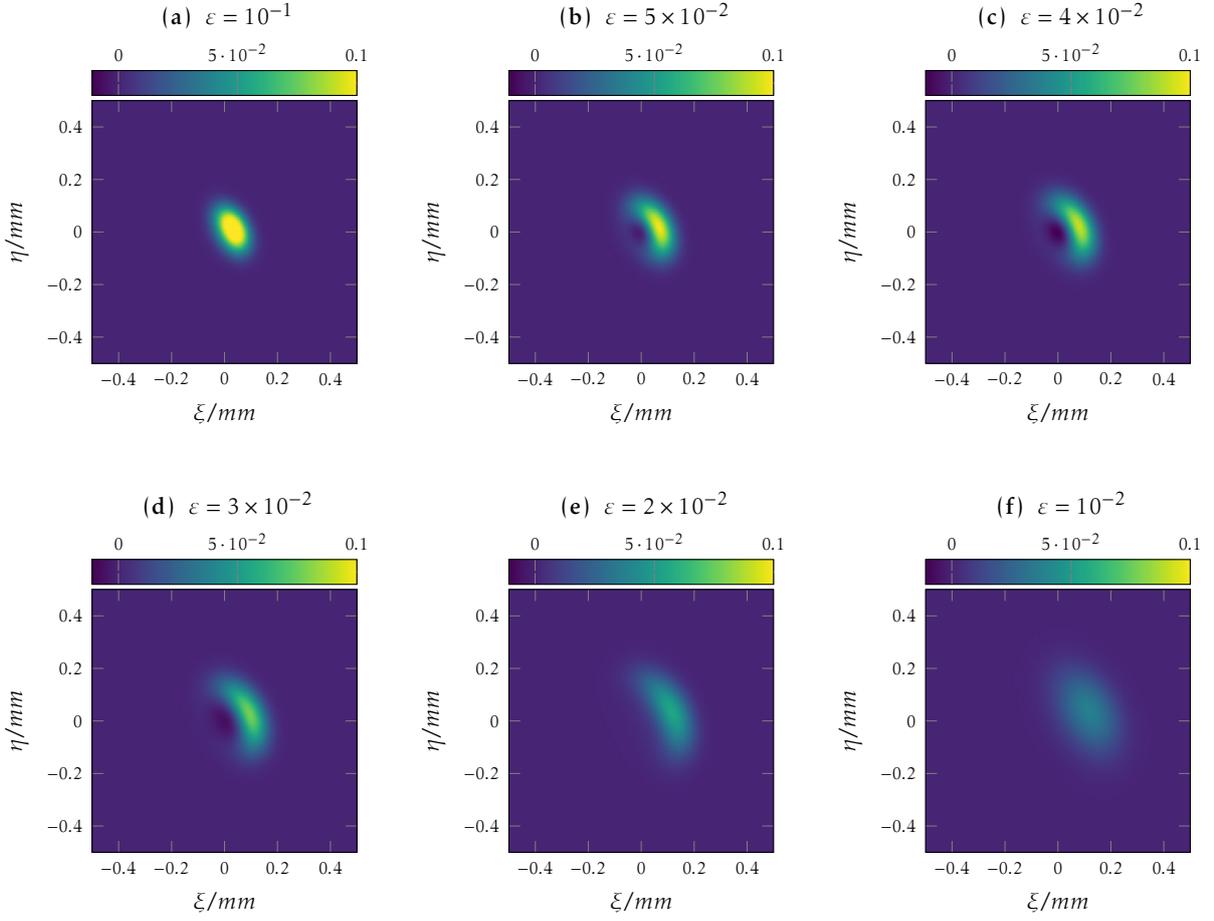

\begin{figure}[h]
	\def\localpath{figures/DiffusionLimit/}
	\centering
	\withfiguresize{\figurethreecol}{\figurethreecol}{\externaltikz{diffusion_limit_diffusive}{\begin{tikzpicture}

\def\dlplot#1#2#3#4#5{
	\nextgroupplot[
		xmin=-0.5,
		xmax=0.5,
		ymin=-0.5,
		ymax=0.5,
		xlabel={$\xcx / mm$},
		ylabel={$\xcy / mm$},
		title={\tikztitle{#1}},
		point meta min ={#2},
		point meta max ={#3},
		colorbar horizontal,
		colorbar style={,
			at={(0,1.02)},
			anchor=south west,
			height=0.02\textwidth,
			width=0.65\figurewidth,
			xticklabel style={font=\figureticklabelfont{\figurewidth},anchor=south},
			ticklabel pos=right,
			max space between ticks=1000pt,
			try min ticks=3,
		},
		#4
	]
	\addplot graphics[xmin=-1,xmax=1,ymin=-1,ymax=1] {\localpath #5};
}

\def\errplot#1#2#3#4#5{
	\nextgroupplot[
		xmin=-0.5,
		xmax=0.5,
		ymin=-0.5,
		ymax=0.5,
		xlabel={$\xcx / mm$},
		ylabel={$\xcy / mm$},
		title={\tikztitle{#1}},
		colorbar horizontal,
		colorbar style={
			colormap name ={cmRdBu_r_minus},
			point meta min ={-#3}, 
			point meta max ={0},
			at={(0,1.02)},
			anchor=south west,
			height=0.02\textwidth,
			width=0.49\figurewidth,
			ticklabel pos=right,
			xticklabel style={font=\figureticklabelfont{\figurewidth},anchor=south, /pgf/number format/.cd,	/tikz/.cd},
			xticklabel={
			  	\pgfmathsetmacro\upBound{0.7*\pgfkeysvalueof{/pgfplots/xmax}+0.3*\pgfkeysvalueof{/pgfplots/xmin}}
			  	\ifdim \tick pt < \upBound pt 
					 $-10^{\pgfmathparse{-\tick+#2}\pgfmathprintnumber\pgfmathresult}$ 
				 \fi
			 },
			try min ticks=2,
		},
		colorbar/draw/.append code={
			\begin{axis}[
				every colorbar, 
				at={(1., 1.02)},
				anchor=south east,
				colormap name={cmRdBu_r_plus},
				colorbar horizontal,
				point meta min={0},
				point meta max={#3},
				colorbar=false,
				height=0.02\textwidth,
				width=0.49\figurewidth,
				ticklabel pos=right,
				xticklabel style={font=\figureticklabelfont{\figurewidth},anchor=south, /pgf/number format/.cd, /tikz/.cd},
				try min ticks=2, 
				xticklabel={
					\pgfmathsetmacro\loBound{0.3*\pgfkeysvalueof{/pgfplots/xmax}+0.7*\pgfkeysvalueof{/pgfplots/xmin}}
						\ifdim \tick pt > \loBound pt 
							$10^{\pgfmathparse{\tick+#2}\pgfmathprintnumber\pgfmathresult}$ 
						\fi
				},
				extra x ticks=0., 
				extra x tick labels=$\pm10^{#2}$,
			] 
			\pgfkeysvalueof{/pgfplots/colorbar addplot}
			\end{axis}
		},		
		#4
	] 
	\addplot graphics[xmin=-1,xmax=1,ymin=-1,ymax=1] {\localpath #5};	
}

\begin{groupplot}[
	group style={group size=3 by 2, horizontal sep = \figurehorizontalsep,  vertical sep = 3cm},
	axis on top,
	scale only axis,
	enlargelimits=false,
	width = \figurewidth,
	height =\figureheight,
	colormap name ={cmviridis},  
	scaled x ticks=false,
   	title style = {font=\figurelabelfont{\figurewidth}, yshift = 0.6cm},
   	xlabel style = {font=\figurelabelfont{\figurewidth}}, 
   	ylabel style = {font=\figurelabelfont{\figurewidth}}, 
   	xticklabel style={font=\figureticklabelfont{\figurewidth}},
   	yticklabel style={font=\figureticklabelfont{\figurewidth}},	
]

	\errplot{$\reldiff{\PM_{\scinum{1}{-2}}}{\density}$}{-4}{3}{}{Fundamental_STO11_Eps1e-2-Eps0.png}
	\errplot{$\reldiff{\PM_{\scinum{8}{-3}}}{\density}$}{-4}{3}{}{Fundamental_STO11_Eps8e-3-Eps0.png}
	\errplot{$\reldiff{\PM_{\scinum{5}{-3}}}{\density}$}{-4}{3}{}{Fundamental_STO11_Eps5e-3-Eps0.png}
	\errplot{$\reldiff{\PM_{\scinum{1}{-3}}}{\density}$}{-4}{3}{}{Fundamental_STO11_Eps1e-3-Eps0.png}
	\errplot{$\reldiff{\PM_{\scinum{1}{-9}}}{\density}$}{-4}{3}{}{Fundamental_STO11_Eps1e-9-Eps0.png} 
	\errplot{$\reldiff{\PM_{\scinum{1}{-3}}}{\PM_{\scinum{1}{-9}}}$} {-4}{3}{}{Fundamental_STO11_Eps1e-3-Eps1e-9.png}
\end{groupplot}

 \end{tikzpicture}}}
	\settikzlabel{fig:fundamental-limit-diffusive-a}
	\settikzlabel{fig:fundamental-limit-diffusive-b}
	\settikzlabel{fig:fundamental-limit-diffusive-c}
	\settikzlabel{fig:fundamental-limit-diffusive-d}
	\settikzlabel{fig:fundamental-limit-diffusive-e}
	\settikzlabel{fig:fundamental-limit-diffusive-f}
	\caption{The numerical solution to the fundamental solution test of \secref{sec:diffusion-limit-convergence}, using the \mischeme{1}-$P_1$ scheme. Each plot shows the relative difference in density $\PM$ between two solutions on a signed truncated logarithmic scale. \figref{fig:fundamental-limit-diffusive-a} - \ref{fig:fundamental-limit-diffusive-e} show the relative difference between the numerical solution at various $\pareps$, and the exact solution \eqref{eq:fundamental-solution}. \figref{fig:fundamental-limit-diffusive-f} shows the difference between the numerical solutions at $\pareps = 10^{-3}$ and $\pareps = 10^{-9}$.}
	\label{fig:fundamental-limit-diffusive}
\end{figure}

\subsection{Convergence analysis with manufactured solutions}
\label{sec:ms-convergence}
Convergence tests with manufactured solutions are useful to detect errors in the scheme and bugs in its implementation. 
If we achieve the expected convergence order we can be more confident that we actually solve the correct problem. 

We only consider the two-dimensional setting.
On the domain 
\begin{align*}
	\domain &= [0, \frac{1}{4}] \times [0,1]^2 \times \US
\end{align*}
we prescribe the solution 
\begin{align*}
	\PD_{ex}(\tcoord, \xcoord, \vcoord) &= \FD(\xcoord, \vcoord) \left( \cos(2 \pi \tcoord) (p_6(\xcx) + p_6(\xcy)) + 2 \right). 
\end{align*}
In terms of the density and perturbation, this is expressed as 
\begin{equation}
	\label{eq:ms-micro-macro}
	\begin{aligned}
	\PM_{ex}(\tcoord, \xcoord) &=  \cos(2 \pi \tcoord) (p_6(\xcx) p_6(\xcy)) + 2 \\
	\PP_{ex}(\tcoord, \xcoord, \vcoord) &= 0. 
	\end{aligned}
\end{equation}
The analytic solution at final time is simply $\evalat{\PM_{ex}}{t = \frac{1}{4}} \equiv 2$,  $\evalat{\PP_{ex}}{t = \frac{1}{4}} \equiv 0$. 
We choose a solution with zero micro part, because this makes the expression for the source easier. 
Nevertheless, due to the coupling of the micro and macro parts, errors in the $\PP$ equation can still be detected with this method. 
The sixth-order polynomial 
\begin{align*}
	p_6(\xi) &= 32 \left( -\xi^6 + 3 \xi^5 - 3 \xi^4 + \xi^3 \right) 
\end{align*}
is carefully chosen such that its value, and its first and second derivative are zero at the boundary:
\begin{align*}
	0 = p_6(0) = p_6(1) = p_6'(0) = p_6'(1) = p_6''(0) = p_6''(1). 
\end{align*}
We add artificial source terms $\hat S_{\PM}, \hat S_{\PP}$ to the right hand side of \eqref{eq:APrho}, \eqref{eq:APremainder} and insert the solution \eqref{eq:ms-micro-macro} to obtain
\begin{align*}
	\hat S_{\PM} &= \dt \PM_{ex} = -2 \sin(2 \pi \tcoord)  (p_6(\xcx) + p_6(\xcy)) \\
	\hat S_{\PP} &= \frac{\pardel}{\pareps^2} \divx (\vcoord \PM_{ex} \FD) -\frac{\pardel}{\pareps^2} \PM_{ex} \lambda_H \gradx \VF \cdot (\vcoord \FD)
\end{align*}
that will produce the desired solution. 

To see the correct order, we need of course a smoothly varying fiber distribution. 
Here we use a distribution with increasing anisotropy along the $\xcx$-axis:
\begin{align*}
	\DW(\xcoord) &= \begin{pmatrix}
	1 + \xcx & 0 & 0 \\ 
	0 & 1 & 0 \\
	0 & 0 & 1
	\end{pmatrix}
\end{align*}

In each convergence test, we refine the grid five times, starting at $20$ grid points and increasing by a factor of 1.5 in each step. 
\begin{figure}[h]
	\def\localpath{figures/convergence/}
	\centering
	\withfiguresize{\figuretwocol}{\figuretwocol}{\externaltikz{convergence_all}{\begin{tikzpicture}
\begin{groupplot}[
    group style={group size=2 by 2, horizontal sep = \figurehorizontalsep,   vertical sep = 2.5cm},
    width = \figurewidth,
    height = \figureheight,
   	title style = {font=\figurelabelfont{\figurewidth}},
   	xlabel style = {font=\figurelabelfont{\figurewidth}}, 
   	ylabel style = {font=\figurelabelfont{\figurewidth}}, 
   	xticklabel style={font=\figureticklabelfont{\figurewidth}},
   	yticklabel style={font=\figureticklabelfont{\figurewidth}},
    scale only axis, 
    xmode=log,
    ymode=log,
    xlabel={\# gridpoints},
    ylabel style={yshift=0.5em},
    ylabel={$\ltwospace$ error},
    log basis x=1.5, 
    log basis y=10,
    xtick={20, 30, 45, 68, 102, 152},
    xticklabel=\pgfmathparse{(1.5)^\tick)}\pgfmathprintnumber{\pgfmathresult},
    x tick label style={/pgf/number format/.cd, fixed, precision=0},
]
\nextgroupplot[title=\tikztitle{\mischeme{1}, $\parnu = 0$}, legend columns=5, legend entries={$\varepsilon = 1$, $\varepsilon = 10^{-1}$, $\varepsilon = 10^{-2}$, $\varepsilon = 10^{-3}$, $\varepsilon = 10^{-5}$, $\mathcal{O}(\Delta x)$, $\mathcal{O}(\Delta x^{2})$}, legend to name={MSConvergenceLegend}, ]
\addplot table [x=nx, y=epsem0, col sep=comma] {\localpath convergence_sto11.csv};
\addplot table [x=nx, y=epsem1, col sep=comma] {\localpath convergence_sto11.csv};
\addplot table [x=nx, y=epsem2, col sep=comma] {\localpath convergence_sto11.csv};
\addplot table [x=nx, y=epsem3, col sep=comma] {\localpath convergence_sto11.csv};
\addplot table [x=nx, y=epsem5, col sep=comma] {\localpath convergence_sto11.csv};
\addplot[domain=20:152, black, thick, mark=] {20*0.1*x^(-1)}; 
\addplot[domain=20:152, black, thick, dashed, mark=] {20*20*0.001*x^(-2)}; 
\nextgroupplot[title=\tikztitle{\mischeme{2}, $\parnu = 0$}]
\addplot table [x=nx, y=epsem0, col sep=comma] {\localpath convergence_sto22_nodrift.csv};
\addplot table [x=nx, y=epsem1, col sep=comma] {\localpath convergence_sto22_nodrift.csv};
\addplot table [x=nx, y=epsem2, col sep=comma] {\localpath convergence_sto22_nodrift.csv};
\addplot table [x=nx, y=epsem3, col sep=comma] {\localpath convergence_sto22_nodrift.csv};
\addplot table [x=nx, y=epsem5, col sep=comma] {\localpath convergence_sto22_nodrift.csv};
\addplot[domain=20:152, black, thick, mark=] {20*0.001*x^(-1)}; 
\addplot[domain=20:152, black, thick, dashed, mark=] {20*20*0.0001*x^(-2)}; 
\nextgroupplot[title=\tikztitle{\mischeme[]{1}, $\parnu = 10$}]
\addplot table [x=nx, y=epsem0, col sep=comma] {\localpath convergence_sto11_eta10.csv};
\addplot table [x=nx, y=epsem1, col sep=comma] {\localpath convergence_sto11_eta10.csv};
\addplot table [x=nx, y=epsem2, col sep=comma] {\localpath convergence_sto11_eta10.csv};
\addplot table [x=nx, y=epsem3, col sep=comma] {\localpath convergence_sto11_eta10.csv};
\addplot table [x=nx, y=epsem5, col sep=comma] {\localpath convergence_sto11_eta10.csv};
\addplot[domain=20:152, black, thick, mark=] {20*0.1*x^(-1)}; 
\addplot[domain=20:152, black, thick, dashed, mark=] {20*20*0.001*x^(-2)}; 
\nextgroupplot[title=\tikztitle{\mischeme[]{2}, $\parnu = 10$}]
\addplot table [x=nx, y=epsem0, col sep=comma] {\localpath convergence_sto22_eta10.csv};
\addplot table [x=nx, y=epsem1, col sep=comma] {\localpath convergence_sto22_eta10.csv};
\addplot table [x=nx, y=epsem2, col sep=comma] {\localpath convergence_sto22_eta10.csv};
\addplot table [x=nx, y=epsem3, col sep=comma] {\localpath convergence_sto22_eta10.csv};
\addplot table [x=nx, y=epsem5, col sep=comma] {\localpath convergence_sto22_eta10.csv};
\addplot[domain=20:152, black, thick, mark=] {20*0.001*x^(-1)}; 
\addplot[domain=20:152, black, thick, dashed, mark=] {20*20*0.0001*x^(-2)}; 
\end{groupplot}
\path[] (group c1r1.north) -- node[anchor=south, yshift=1cm]{\ref{MSConvergenceLegend}} (group c2r1.north); 
\end{tikzpicture}}}
	\settikzlabel{fig:ms-convergence-sto11}
	\settikzlabel{fig:ms-convergence-sto22}
	\settikzlabel{fig:ms-convergence-sto11-drift}
	\settikzlabel{fig:ms-convergence-sto22-drift}
	\caption{$\ltwospace$ errors to manufactured solution at the final time for various values of $\pareps$}
	\label{fig:ms-convergence}
\end{figure}
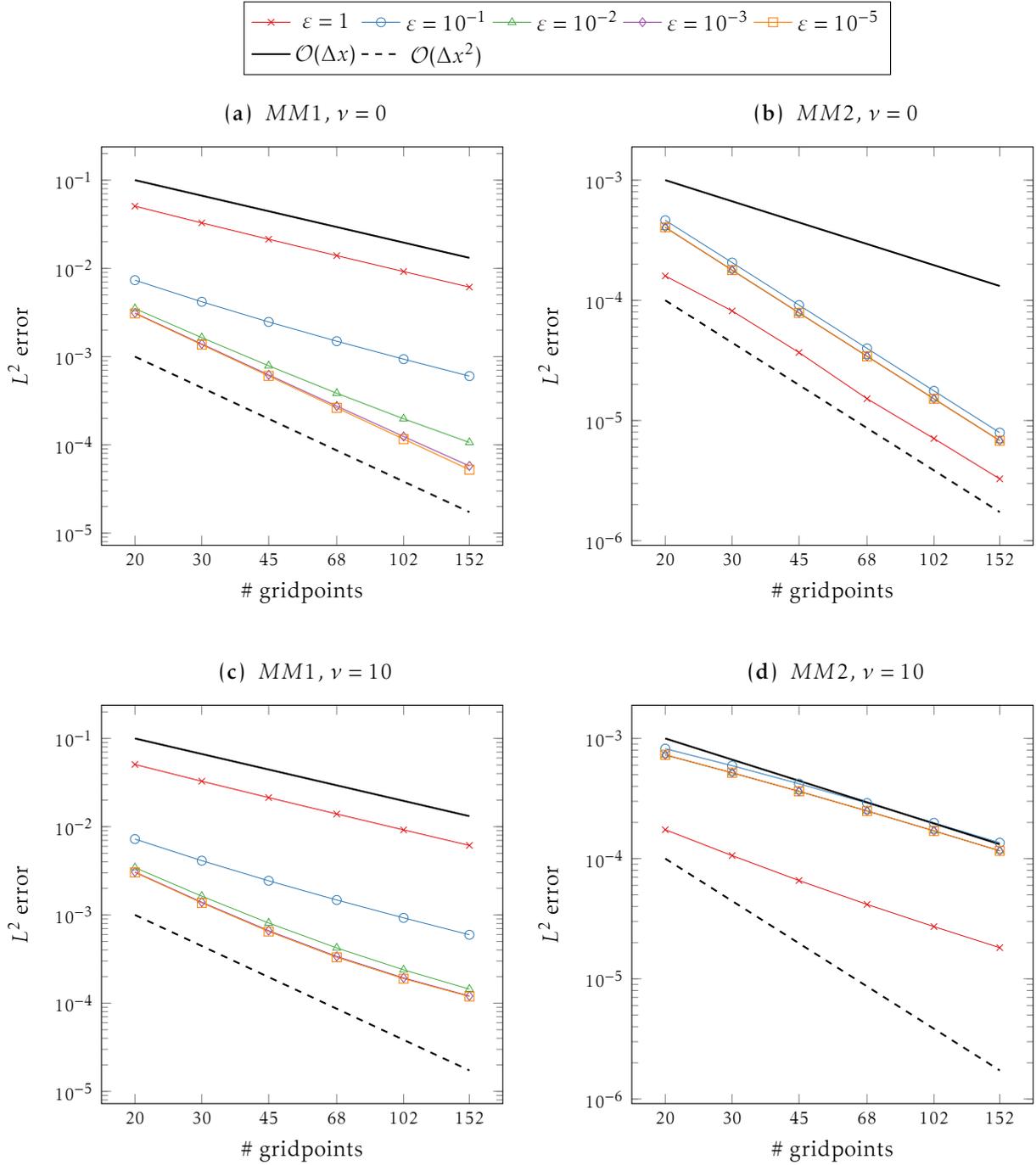

We set $\pardel = 0.1$ and ignore natural growth, i.e., set $\parthet = 0$. 
Convergence tests were run with the first and the second order code, each with advection $\parnu =  10$ and without advection $\parnu = 0$. 
Each of these tests was repeated for different values of the scaling parameter $\pareps$ ranging from one to $10^{-5}$. 
The results are plotted in \figref{fig:ms-convergence}. 

Without the drift $\parnu = 0$, the first order code (see \figref{fig:ms-convergence-sto11}) shows the expected first order of convergence in the kinetic regime $\pareps = 1$ and second order of convergence in the diffusive regime $\pareps = 10^{-5}$. 
In the transition between the regimes, the convergence order increases from one to two. 
As expected, this increase in order is lost when the drift term is active (see \figref{fig:ms-convergence-sto11-drift}) and the convergence order is one for all considered values of $\pareps$. 
We observe second order convergence for the second order code without drift, independently of the flow regime(see \figref{fig:ms-convergence-sto22}). 
However, presence of the drift term reduces the order to one(see \figref{fig:ms-convergence-sto22-drift}). 
This is due to the first order approximation of the drift term. 
The second order code still produces smaller absolute errors than the first order code. 
Interestingly, absolute errors for the second order code are much smaller with $\pareps = 1$ compared to all other values of $\pareps$. 

\subsection{Strong discontinuities in the diffusion coefficients}
\label{sec:discontinuities}
The coefficients in the glioma model from \secref{sec:glioma-model} are estimated from DTI measurements of the brain, which give a water diffusion tensor $\DW$ per voxel. 
Voxels typically have a length of a few millimeters. 
On each voxel, the tensor is assumed constant and as such the resulting coefficients jump across the voxel boundaries. 
Apart from these artifacts, there are genuine jumps in the data when the underlying tissue orientation changes rapidly. 
Thus we are interested in the behavior of our scheme in the presence of discontinuous coefficients, especially if $\pareps$ is small. 

In the context of flow through porous media, a number of benchmarks with strong jumps in the diffusion coefficient have been developed\cite{eigestad2005convergence, rivie2000part}. 
We adapt two benchmarks with an analytical solution for our scheme. 
The first is a special case of a benchmark with discontinuities in permeability at quadrant boundaries from Eigestad and Klausen \cite{eigestad2005convergence} which we call isotropic quadrants test.
The domain is divided into four quadrants of which each is assigned a constant and isotropic permeability. 
The other test is similar to the 'piecewise uniform flow' in \cite{eigestad2005convergence}. 
It features two domains of constant diffusion tensor with a single discontinuity. 
But here we align the discontinuity with the $\xcoord_2$-axis and choose constant anisotropic diffusion tensors whose main axes meet at an angle at the interface. 

Note that the benchmarks are designed for the stationary porous media equation
\begin{align*}
	\divx (\Difftens \gradx \density ) &= 0
\end{align*}
Our code is neither stationary nor does it solve the porous media equation. 
If growth and drift are neglected, the code should approximately solve 
\begin{align}
	\label{eq:diffusion-benchmark-eq}
	\dt \density - \pardel \divx (\divx (\Difftens \density)) &= 0. 
\end{align}
for small $\pareps$. 
However, we can run the simulations for a long enough time $T^*$, until a steady state is reached and choose a very small $\pareps$, e.g., $10^{-5}$. 
In the steady state, the choice of $\pardel$ does not play a role. 
Effectively, this is a very inefficient iterative solver of the stationary equation. 
As a convergence criterion we use the relative $\ltwospace$-difference between successive time steps, i.e, we abort the simulation if 
\begin{align*}
	\frac{\Vert \PM(\tcoord_{i-1}) - \PM(\tcoord_i)\Vert_2}{\Vert \PM(\tcoord_i) \Vert_2 \timestep_i} &< tol. 
\end{align*}

In the benchmarks, we prescribe Dirichlet boundary conditions for $\PM$ according to the exact solution and Maxwellian boundary conditions \eqref{eq:bc-micro} for the micro equation $\PP$.

\subsubsection{Quadrants with jump in permeability}
\label{sec:quadrants}
First, we switch to polar coordinates
\begin{align*}
	\begin{pmatrix}
	\xcx \\ \xcy 
	\end{pmatrix}
	&= 
	\polradius \begin{pmatrix}
	\cos(\polangle) \\
	\sin(\polangle)
	\end{pmatrix}.
\end{align*}
The $i$-th quadrant is then $\quadrant{i} = (\polradius, \polangle) \in [0, \infty) \times [\frac{i \pi}{2}, \frac{(i+1) \pi}{2})$, for $i = 0,\dots 3$. 
On each quadrant, we have a constant isotropic diffusion tensor $\Difftens_i = \kappa_i \Identity$. 
\begin{figure}
	\centering
	\def\localpath{figures/DiffusionBenchmark/}
	\withfiguresize{\figuretwocollegend}{\figuretwocollegend}{\externaltikz{quadrants_solution}{\begin{tikzpicture}[
	mainline/.style={draw, thick},
	auxline/.style ={draw, gray},
]

\begin{groupplot}[
	group style={group size=2 by 1, horizontal sep = \figurehorizontalsep,  vertical sep = 2.5cm},
	scale only axis,
	enlargelimits=false,
	width = \figurewidth,
	height = \figureheight,
   	title style = {yshift = 0.6cm, font=\figurelabelfont{\figurewidth}},
   	xlabel style = {font=\figurelabelfont{\figurewidth}}, 
   	ylabel style = {font=\figurelabelfont{\figurewidth}}, 
   	xticklabel style={font=\figureticklabelfont{\figurewidth}},
   	yticklabel style={font=\figureticklabelfont{\figurewidth}},
	colormap name ={cmviridis},  
]
	
	\nextgroupplot[
		axis on top,
		scaled x ticks=false,
		xmin=-1,
		xmax=1,
		ymin=-1,
		ymax=1,
		xticklabel shift = 2pt,
		yticklabel shift = 2pt,
		title={\tikztitle{Analytic solution}},
		xlabel={$\xcx/$mm},
		ylabel={$\xcy/$mm},
		colorbar horizontal,
		colorbar style={,
			at={(0,1.02)},
			anchor=south west,
			height=0.02\textwidth,
			width=\figurewidth,
			xticklabel style={font=\figureticklabelfont{\figurewidth},anchor=south},
			ticklabel pos=right,
			max space between ticks=1000pt,
			try min ticks=3,
		},
		point meta min ={-1.050174217},
		point meta max ={1.050174217}, 
	]
	\addplot graphics[xmin=-1,xmax=1,ymin=-1,ymax=1] {\localpath quadrants_solution_plot.png};
	
	\nextgroupplot[
		xmode=log,
		ymode=log,
		xlabel={\# gridpoints},
		ylabel style={yshift=0.5em},
		ylabel={$L_2$ error},
		log basis x=1.5, 
		log basis y=10,
		xtick={20, 30, 45, 68, 102, 152},
		xticklabel=\pgfmathparse{(1.5)^\tick)}\pgfmathprintnumber{\pgfmathresult},
		x tick label style={/pgf/number format/.cd, fixed, precision=0},
		legend style={at={(1.02,1.)},anchor=north west, font=\figureticklabelfont{\figurewidth}},
		title=\tikztitle{Convergence},
	]
	\addplot table [x=nx, y=epsem5, col sep=comma] {\localpath quadrants_sto11.csv};
	\addplot[domain=20:152, black, thick, mark=] {20^(0.253804139)*0.1602*x^(-0.253804139)}; 
	\addplot[domain=20:152, black, thick, dashed, mark=] {20^(0.4)*0.1602*x^(-0.4)}; 
	\addlegendentry{$e_n$};
	\addlegendentry{$\mathcal{O}(\Delta x^{2\alpha})$}; 
	\addlegendentry{$\mathcal{O}(\Delta x^{0.4})$}; 
\end{groupplot}

  
\end{tikzpicture}}}
	\settikzlabel{fig:quadrants-solution}
	\settikzlabel{fig:quadrants-convergence}
	\caption{The benchmark described in \secref{sec:quadrants}, with discontinuities in permeability at the quadrant boundaries.
 \ref{fig:quadrants-solution}: Analytic solution \eqref{eq:quadrants-solution} for the permeability values in \tabref{tab:quadrants-coeffs}. \ref{fig:quadrants-convergence}: Convergence of $\ltwospace$-error with respect to grid refinement.}
	\label{fig:quadrants}
\end{figure}
The stationary solution to \eqref{eq:diffusion-benchmark-eq} has the form 
\begin{align}
	\label{eq:quadrants-solution}
	\density_{,ex}(\polradius, \polangle) &= \polradius^{\alpha} \left( a_i \cos(\alpha \polangle) + b_i \sin(\alpha \polangle) \right)  & (\polradius, \polangle) \in \quadrant{i}, 
\end{align}
with coefficients $\alpha, a_i, b_i$ determined by the continuity of the density and the flux at the interfaces. 
Continuity of the density gives the four conditions
\begin{align*}
	\density_{,ex}(\polradius, \polangle_i^-) = \density_{,ex}(\polradius, \polangle_i^+), 
\end{align*}
wherein $\polangle_i^\pm$ mean that the interface at $\frac{i \pi}{2}$ is approached from the left or the right. 
Continuity of the fluxes translates into the conditions
\begin{align*}
	\pfrac{}{\outernormal} \Difftens \density_{,ex}(\polradius, \polangle_i^-) = \pfrac{}{\outernormal} \Difftens \density_{,ex}(\polradius, \polangle_i^+), 
\end{align*}
with 
\begin{align*}
	\pfrac{}{\outernormal} \Difftens \density_{,ex}(\polradius, \polangle) = \kappa \pfrac{}{\outernormal} \density_{,ex} = \alpha \polradius^{\alpha-1}(-a_i \sin(\alpha \polangle) + b_i \cos(\alpha \polangle)).
\end{align*}
Here we used that on each quadrant the coefficients are constant. 
Altogether we have eight conditions for nine coefficients. 
We arbitrarily set $a_0 = 1$ and solve for the remaining coefficients numerically. 

Similar to \cite{eigestad2005convergence}, we take the permeability $\kappa$ equal at diagonally opposite quadrants, and set 
\begin{align*}
	\kappa_0 &= \kappa_2 = 100, \\
	\kappa_1 &= \kappa_3 = 1.
\end{align*}
In the code this is achieved by prescribing the turning rate 
\begin{align*}
\lambda_0(\xcoord) = \frac{3}{\kappa}
\end{align*}
and an isotropic water diffusion tensor $\DW = \Identity$. 

The coefficients that belong to this choice are listed in \tabref{tab:quadrants-coeffs}. 
They are identical to the values reported in \cite{eigestad2005convergence}. 
A plot of the analytic solution \eqref{eq:quadrants-solution} corresponding to these coefficients is shown in \figref{fig:quadrants-solution}.
Due to the discontinuous permeability, the solution to the diffusion equation only belongs to the fractional Sobolev space $H^{1+\alpha-\nu}, \forall \nu > 0$, i.e., it is at most $1+\alpha$ times differentiable. 
Therefore the maximum order of convergence we can expect with respect to grid refinement is $2\alpha$. 
We performed a grid refinement study with five refinements, a refinement factor of $1.5$ and 20 grid points on the coarsest grid. 
Surprisingly the observed order of convergence(see \figref{fig:quadrants-convergence}) is about $0.4$---significantly greater than the theoretical order $2\alpha \approx 0.25$. 
The error at $45$ grid points is exceptionally large because for an odd number of grid points, the quadrant boundary does not coincide with primal cell edges. 

\begin{table}
	\centering
	\begin{tabular}{lllll}
		$i$ & 0 & 1 & 2 & 3 \\
		\hline
		$\kappa_i $ &  100. & 1. & 100. &1. \\
		$a_i$ & 1.  & 2.96039604 &-0.88275659 &-6.45646175 \\
		$b_i$ & 0.1 & -9.6039604 &-0.48035487 & 7.70156488 \\
		\hline
		$\alpha$ & \multicolumn{4}{l}{0.126902069721}	
	\end{tabular}
	\caption{Coefficients for the exact solution \eqref{eq:quadrants-solution} of the quadrants test described in \secref{sec:quadrants}.}
	\label{tab:quadrants-coeffs}
\end{table}

\subsubsection{Interface with change in diffusion tensor axis}
\label{sec:halfplane}
In this test, the diffusion tensor is constant but anisotropic on the left and right half-planes. 
At the interface--the $\xcy$-axis--there is an abrupt change in the main direction of diffusion. 
Let $R(\theta) \in SO(3)$ a rotation around the $\xcz$-axis with angle $\theta$. 
The diffusion tensor field is parametrized by left and right anisotropies $a^L, a^R$ and left and right angles of main diffusion $\theta^L, \theta^R$: 
\begin{align*}
	\Difftens(\xcoord) &= \begin{cases}
		\Difftens^{L} = \frac{1}{a^L + 2} R\trans(\theta^L) \begin{pmatrix}
		a^L & 0  & 0 \\
		0 & 1 & 0 \\
		0 & 0 & 1 
		\end{pmatrix} 
		R(\theta^L) & \xcx < 0\\
		\Difftens^{R} = \frac{1}{a^R + 2} R\trans(\theta^R) \begin{pmatrix}
		a^R & 0  & 0 \\
		0 & 1 & 0 \\
		0 & 0 & 1 
		\end{pmatrix} 
		R(\theta^R) & \xcx >0
	\end{cases}
\end{align*}
The piecewise linear function 
\begin{align*}
	\density(\xcoord) = \begin{cases}
	\density^{L} = s^L \cdot \xcoord & \xcx < 0 \\
	\density^{R} = s^R \cdot \xcoord & \xcx > 0
	\end{cases} 
\end{align*}
is a stationary solution to the diffusion equation \eqref{eq:diffusion-benchmark-eq} on each half-plane. 
For a given left slope $s^L$, we use the continuity of the solution and normal fluxes at the interface to compute the right slope $s^R$. 
Continuity of the solution gives us 
\begin{align*}
	s^{R}_{\xcy} = s^{L}_{\xcy},
\end{align*}
and continuity of the normal fluxes translates to 
\begin{align*}
	\divx( \Difftens^L \density^{L}(0^-, \xcy)) \cdot e_1 &= \divx( \Difftens^R \density^{R}(0^+, \xcy) ) \cdot e_2\\
	\Difftens^{L}_{\xcx \xcx} s^{L}_{\xcx} + \Difftens^{L}_{\xcy \xcx} s^{L}_{\xcy} &= \Difftens^{R}_{\xcx \xcx} s^{R}_{\xcx} + \Difftens^{R}_{\xcy \xcx} s^{R}_{\xcy} \\
	s^{R}_{\xcx} &= \frac{1}{\Difftens^{R}_{\xcx \xcx}} \left( -\Difftens^{R}_{\xcy \xcx} s^{R}_{\xcy} + \Difftens^{L}_{\xcx \xcx} s^{L}_{\xcx} + \Difftens^{L}_{\xcy \xcx} s^{L}_{\xcy}\right). 
\end{align*}
We compute two different situations whose parameters are summarized in \tabref{tab:diffusion-benchmark-halfplane-coeffs} and that only differ in the tangential flux at the interface, which is determined by $s^{L}_{\xcy}$.  
In the first test there is no tangential flux at the interface. 
In this case the numeric solution is identical to the analytic solution. 
However, in the second test in which the tangential flux component is not zero, the numeric solution differs significantly from the analytic solution. 
Relative differences in density and flux between the computation results on a $50\times50$ grid and the analytic solution are plotted in \figref{fig:tensorjump}.
The errors are largest at the interface, especially at the lower and upper boundary. 
In density the error is about $10\%$, but in the fluxes it reaches $300\%$. 

\begin{table}
	\centering
	\renewcommand{\arraystretch}{1.2}	
	\begin{tabular}{llll}
		&  & L & R\\
		\hline
		&$\theta$ & $80$\textdegree & $20$\textdegree \\
		&$a$ & 2.5 & 2.5 \\ 
		\hline
		Test 1&$s$ & (1,0) & (0.44965177,  0. ) \\
		Test 2&$s$ & (1,1) & (0.35261053,  1. ) 
	\end{tabular}
	\caption{Coefficients for the anisotropic half-plane test described in \secref{sec:halfplane}. }
	\label{tab:diffusion-benchmark-halfplane-coeffs}
\end{table}

\begin{figure}
	\centering
	\def\localpath{figures/DiffusionBenchmark/}
	\withfiguresize{\figurethreecol}{\figurethreecol}{\externaltikz{tensorjump}{\begin{tikzpicture}

\def\dlplot#1#2#3#4#5{
	\nextgroupplot[
		xmin=-1,
		xmax=1,
		ymin=-1,
		ymax=1,
		xlabel={$\xcx / mm$},
		ylabel={$\xcx / mm$},
		title={\tikztitle{#1}},
		point meta min ={#2},
		point meta max ={#3},
		colorbar horizontal,
		colorbar style={,
			at={(0,1.02)},
			anchor=south west,
			height=0.02\textwidth,
			width=0.65\figurewidth,
			xticklabel style={font=\figureticklabelfont{\figurewidth},anchor=south},
			ticklabel pos=right,
			max space between ticks=1000pt,
			try min ticks=3,
		},
		#4
	]
	\addplot graphics[xmin=-1,xmax=1,ymin=-1,ymax=1] {\localpath #5};
}

\def\errplot#1#2#3#4#5{
	\nextgroupplot[
		xmin=-1,
		xmax=1,
		ymin=-1,
		ymax=1,
		xlabel={$\xcx / mm$},
		ylabel={$\xcy / mm$},
		title={\tikztitle{#1}},
		colorbar horizontal,
		colorbar style={
			colormap name ={cmRdBu_r_minus},
			point meta min ={-#3}, 
			point meta max ={0},
			at={(0,1.02)},
			anchor=south west,
			height=0.02\textwidth,
			width=0.49\figurewidth,
			ticklabel pos=right,
			xticklabel style={font=\figureticklabelfont{\figurewidth},anchor=south, /pgf/number format/.cd,	/tikz/.cd},
			xticklabel={
			  	\pgfmathsetmacro\upBound{0.7*\pgfkeysvalueof{/pgfplots/xmax}+0.3*\pgfkeysvalueof{/pgfplots/xmin}}
			  	\ifdim \tick pt < \upBound pt 
					 $-10^{\pgfmathparse{-\tick+#2}\pgfmathprintnumber\pgfmathresult}$ 
				 \fi
			 },
			try min ticks=2,
		},
		colorbar/draw/.append code={
			\begin{axis}[
				every colorbar, 
				at={(1., 1.02)},
				anchor=south east,
				colormap name={cmRdBu_r_plus},
				colorbar horizontal,
				point meta min={0},
				point meta max={#3},
				colorbar=false,
				height=0.02\textwidth,
				width=0.49\figurewidth,
				ticklabel pos=right,
				xticklabel style={font=\figureticklabelfont{\figurewidth},anchor=south, /pgf/number format/.cd, /tikz/.cd},
				try min ticks=2, 
				xticklabel={
					\pgfmathsetmacro\loBound{0.3*\pgfkeysvalueof{/pgfplots/xmax}+0.7*\pgfkeysvalueof{/pgfplots/xmin}}
						\ifdim \tick pt > \loBound pt 
							$10^{\pgfmathparse{\tick+#2}\pgfmathprintnumber\pgfmathresult}$ 
						\fi
				},
				extra x ticks=0., 
				extra x tick labels=$\pm10^{#2}$,
			] 
			\pgfkeysvalueof{/pgfplots/colorbar addplot}
			\end{axis}
		},		
		#4
	] 
	\addplot graphics[xmin=-1,xmax=1,ymin=-1,ymax=1] {\localpath #5};	
}

\begin{groupplot}[
	group style={group size=3 by 1, horizontal sep = \figurehorizontalsep,  vertical sep = 3cm},
	axis on top,
	scale only axis,
	enlargelimits=false,
	width = \figurewidth,
	height =\figureheight,
	colormap name ={cmviridis},  
	scaled x ticks=false,
   	title style = {yshift = 0.6cm, font=\figurelabelfont{\figurewidth}},
   	xlabel style = {font=\figurelabelfont{\figurewidth}}, 
   	ylabel style = {font=\figurelabelfont{\figurewidth}}, 
   	xticklabel style={font=\figureticklabelfont{\figurewidth}},
   	yticklabel style={font=\figureticklabelfont{\figurewidth}},
]
%
	\errplot{$\reldiff{\PM}{\density}$}{-4}{5}{}{Tensorjump_diagonal_rho.png}
	\errplot{$\reldiff{\ints{\vcoord_{\xcx} \PP}}{\ints{\vcoord_{\xcx} \PP}_0}$}{-4}{5}{}{Tensorjump_diagonal_flux0.png}
	\errplot{$\reldiff{\ints{\vcoord_{\xcx} \PP}}{\ints{\vcoord_{\xcx} \PP}_0}$}{-4}{5}{}{Tensorjump_diagonal_flux1.png}
\end{groupplot}

 \end{tikzpicture}}}
	\settikzlabel{fig:tensorjump-a}
	\settikzlabel{fig:tensorjump-b}
	\settikzlabel{fig:tensorjump-c}
	\caption{Numerical solution to the half-plane test in \secref{sec:halfplane}. Shown are the relative errors in density $\PM$ (\ref{fig:tensorjump-a}) and flux components $\ints{\vcoord_{\xcx} \PP}, \ints{\vcoord_{\xcy} \PP}$ (\ref{fig:tensorjump-b}, \ref{fig:tensorjump-c}) on a signed truncated logarithmic scale. }
	\label{fig:tensorjump}
\end{figure}

\subsection{Computation using DTI data from human brain}
\label{sec:brain}
To demonstrate the full capabilities of the scheme we compute the model of glioma invasion in the human brain (see \secref{sec:glioma-model}) with the parameters in \tabref{tab:brain-parameters}. 
We do not claim that these parameters, which are similar to those in \cite{EHS}, are accurate at all, not even to an order of magnitude. 
But the results are qualitatively similar to clinical observations (see e.g. \cite{swan2018patient}) and therefore serve as a starting point to test the scheme under more realistic conditions. 
The diffusion tensor field $\DW$ is the same as in \cite{EHS, corbin2018higher}. 
It remains to estimate the volume fraction $\VF[\DW](\xcoord)$ and the function $\hat \lambda_H$. 
We use the same estimates as in \cite{EHKS14,corbin2018higher}, namely  
\begin{align*}
\nonumber\VF[\DW](\xcoord) &= CL(\DW(\xcoord)) := 1 - \left( \frac{\trace (\DW)}{4 \max\Eigenvalue(\DW)} \right)^{\frac{3}{2}}.
\end{align*}
for the volume fraction and 
\begin{align*}
\hat \lambda_H[\VF](\xcoord) &= \frac{1}{1 + \frac{\alpha(\VF)}{\lambda_0}}h'(\VF),  \\
\alpha &= k^+ \VF + k^-, \\
h  &= \frac{k^+ \VF}{\alpha},
\end{align*}
for the activation function, with positive constants $k^+, k^-$. 
We are not interested in absolute values of $\PM$ but rather in the ratio $\frac{\PM}{\carryingcapacity}$ and therefore set the carrying capacity $\carryingcapacity = 1$ in the computations. 
A two dimensional slice through the three dimensional data set is visualized in \figref{fig:brain-setup}. 
The two dimensional computations are performed on a $40mm \times 40mm$ square subset of that slice.
We simulate the tumor growth over a time span of two years starting from its original appearance at an isolated site. 
Therefore, initially we set $\PM = 1$ on the grid cell at the center of the computation domain and $\PM = 0$ everywhere else. 
It is reasonable to assume that the tumor starts in equilibrium, i.e., $\PP(0, \xcoord) = 0$ everywhere. 
In \figref{fig:brain-evolution} snapshots of the simulated density $\PM$ at half-year intervals are shown. 
The tumor evolves basically like a traveling wave in the Fisher equation with heterogeneous wave speed due to the heterogeneous diffusion tensor and drift. 
We observe an increased speed of the invasion front along the white matter tracts. 
The solution inside the tumor is almost stationary and fluctuates around the carrying capacity of the growth. 
Note that due to the drift, the model allows migration into regions that are already full and thus the density can become larger than the carrying capacity.
This can be seen in \figref{fig:brain-evolution-4c}, wherein we show contours of $\rho$ at selected percentages of the carrying capacity. 

Next we compare solutions of the model with different settings. 
Therefore, in \figref{fig:brain-comparisons} we plot the $10\%$ contour lines of the solution at the final time of two years. 

In \figref{fig:brain-comparisons-a}, we compare solutions for various values of $\pareps$, between $\pareps = \scinum{1}{-3}$ and $\pareps = \scinum{1}{-10}$.  
The original parameters describe a situation very close to the diffusion limit, with $\pareps \approx \scinum{3.3}{-6}$.  
As we can expect from the results in \secref{sec:diffusion-limit-convergence} there is no difference between the original model and the model with $\pareps = \scinum{1}{-10}$. 
However, we start to see differences when we artificially choose a greater $\pareps$. 
Generally the invasion front is faster for greater $\epsilon$ because due to a reduced turning rate individual cells have a higher chance of overtaking the diffusive invasion front. 
At $\pareps = \scinum{1}{-4}$ we observe a distance of contours comparable to the $2mm$ resolution of the DTI data set. 
This value of $\pareps$ corresponds to the cell speed $\speed \approx \scinum{6.9}{-6} \frac{mm}{s}$, which is about a hundredth of the original value, and turning rates $\lambda_0 \approx \scinum{8.6}{-4}, \lambda_1 \approx \scinum{1.1}{-1}$ approximately one thousandth of the original rates. 
Thus the kinetic model could be relevant for cell species that migrate very slowly and change their orientation very rarely (in this example once every twenty minutes).  

We also investigate the influence of the spatial and temporal discretization scheme on the solution and compare the \mischeme[]{1}, \ivscheme[]{1}, \mischeme[]{2}, \ivscheme{2} schemes (see \figref{fig:brain-comparisons-b}).
The second-order variants \mischeme[]{2} and \ivscheme[]{2} agree very well in most of the domain. 
The contours of both second-order schemes lie between the contour for the \mischeme[]{1} scheme on the inside and the contour for the \ivscheme[]{1} scheme on the outside everywhere. 
Hence the \mischeme[]{1} scheme seems to underestimate the invasion front, whereas the \ivscheme[]{1} scheme overestimates it. 
Considering the explicit and implicit distretizations of $\dot x = x$, this behavior is to be expected. 

Because the situation is very close to the diffusion limit, higher moment orders in the velocity discretization make no difference to the solution. 
In \figref{fig:brain-comparisons-c} the contours for the $P_1$ and the $P_3$ solutions are plotted and are visually identical. 

Finally, we compare the solution of the two-dimensional model with a slice of the three dimensional model (see \figref{fig:brain-comparisons-d}). 

\begin{table}
	\centering
	\renewcommand{\arraystretch}{1.2}
	\begin{tabular}{l@{\hspace{5pt}}lr@{.\hspace{0pt}}l@{$\times$\hspace{0pt}}lcl}
		\multicolumn{2}{l}{Parameter}   &\multicolumn{3}{l}{Value}&                  & Description \\
		\hline 
		T           && $6$&$31$&  $10^{ 7}$  & $s$              & time span = one year \\
		c           && $2$&$1$   &  $10^{-4}$  & $ \frac{mm}{s}$  & cell speed \\
		$\lambda_0$ && $8$&$0$   &  $10^{-1}$  & $\frac{1}{s}$    & cell-state independent part of turning rate\\
		$\lambda_1$ && $1$&$5$   &  $10^{ 2}$  & $\frac{1}{s}$    & cell-state dependent part of turning rate\\
		$k^+$       && $1$&$0$   &  $10^{-1}$  & $\frac{1}{s} $   & attachment rate of cells to ECM\\
		$k^-$       && $1$&$0$   &  $10^{-1}$  & $\frac{1}{s} $   & detachment rate of cells to ECM\\	
		$\factorprolref$ && $8$&$44$ &$10^{-7}$& $\frac{1}{s} $   & growth rate\\
		\hline		
		$\St$       &&$1$&$21$  &  $10^{-2}$  &                  & Strouhal number\\
		$\Kn$       &&$3$&$96$  &  $10^{-8}$  &                  & Knudsen number\\
		$\pareps$   &&$3$&$28$  &  $10^{-6}$  &                  & parabolic scaling number\\
		$\pardel$   &&$2$&$72$  &  $10^{-4}$  &                  & parabolic scaling number\\
		$\parnu$   &&$1$&$25$  &  $10^{2}$   &                  & Ratio of turning rate coefficients\\ 
	\end{tabular}
	\caption{The reference parameters and the resulting characteristic numbers used in the simulations of glioma invasion in the human brain.}
	\label{tab:brain-parameters}
\end{table}

\begin{figure}
	\centering
	\def\localpath{figures/Brain/}
	\withfiguresize{\figurethreecol}{\figurethreecol}{\externaltikz{brain_time_evolution}{\begin{tikzpicture}

\def\colorplot#1#2#3#4#5{
	\nextgroupplot[
		xmin=32,
		xmax=112,
		ymin=32,
		ymax=112,
		xlabel={$\xcx / mm$},
		ylabel={$\xcy / mm$},
		title={\tikztitle{#1}},
		point meta min ={#2},
		point meta max ={#3},
		colorbar horizontal,
		colorbar style={,
			at={(0,1.02)},
			anchor=south west,
			height=0.02\textwidth,
			width=\figurewidth,
			xticklabel style={font=\figureticklabelfont{\figurewidth},anchor=south},
			ticklabel pos=right,
			max space between ticks=1000pt,
			try min ticks=3,
		},
		#4
	]
	\addplot graphics[xmin=32,xmax=112,ymin=32,ymax=112] {\localpath #5};
}

\begin{groupplot}[
	group style={group size=3 by 2, horizontal sep = \figurehorizontalsep,  vertical sep = 3cm},
	axis on top,
	scale only axis,
	enlargelimits=false,
	width = \figurewidth,
	height =\figureheight,
	colormap name ={cmviridis},  
	scaled x ticks=false,
   	title style = {yshift = 0.6cm, font=\figurelabelfont{\figurewidth}},
   	xlabel style = {font=\figurelabelfont{\figurewidth}}, 
   	ylabel style = {font=\figurelabelfont{\figurewidth}}, 
   	xticklabel style={font=\figureticklabelfont{\figurewidth}},
   	yticklabel style={font=\figureticklabelfont{\figurewidth}},
]
\nextgroupplot[
width = \figurewidth/ 284 * 234,
height = \figureheight,
xmin=0,
xmax=234,
ymin=0,
ymax=284,
xlabel={$\xcx / mm$},
ylabel={$\xcy / mm$},
title={\tikztitle{Initial}},
point meta min ={0},
point meta max ={1.05},
]
\addplot graphics[xmin=0,xmax=234,ymin=0,ymax=284] {\localpath basic2d_initial.png};
\draw[white!90!black, very thick] (32, 32) rectangle (112, 112); 
\draw[white!90!black, very thick, ->] (102, 72) -- (72, 72); 

\colorplot{$\tcoord = 0.5y$}{0}{1.05}{}{basic2d_1.png}
\colorplot{$\tcoord = 1y$}{0}{1.05}{}{basic2d_2.png}
\colorplot{$\tcoord = 1.5y$}{0}{1.05}{}{basic2d_3.png}
\colorplot{$\tcoord = 2y$}{0}{1.05}{}{basic2d_4.png}
\colorplot{$\tcoord = 2y$}{-3}{0}{
	colorbar style={
		xticklabel style={font=\footnotesize,anchor=south, /pgf/number format/.cd,	/tikz/.cd},
		xticklabel={$10^{\pgfmathparse{\tick}\pgfmathprintnumber\pgfmathresult}$},
	}
}{basic2d_4_contours.png}
\end{groupplot}

 \end{tikzpicture}}}
	\settikzlabel{fig:brain-setup}
	\settikzlabel{fig:brain-evolution-1}
	\settikzlabel{fig:brain-evolution-2}
	\settikzlabel{fig:brain-evolution-3}
	\settikzlabel{fig:brain-evolution-4}
	\settikzlabel{fig:brain-evolution-4c}
	\caption{\ref{fig:brain-setup}: A two dimensional slice through the DTI data set. The white box indicates the computational domain and the white arrow the initial tumor location. \ref{fig:brain-evolution-1} - \ref{fig:brain-evolution-4}: Plots of the glioma simulation in six-month intervals. \ref{fig:brain-evolution-4c}: Contours at $100\%, 10\%, 1\%$ and $0.1\%$ of the carrying capacity at the final time. Tumor density $\PM$ is shown in color. The volume fraction $\VF$ is encoded in the grayscale background image; brighter color means greater $\VF$. The black arrows show the limit drift vector $\driftT$. }
	\label{fig:brain-evolution}
\end{figure}
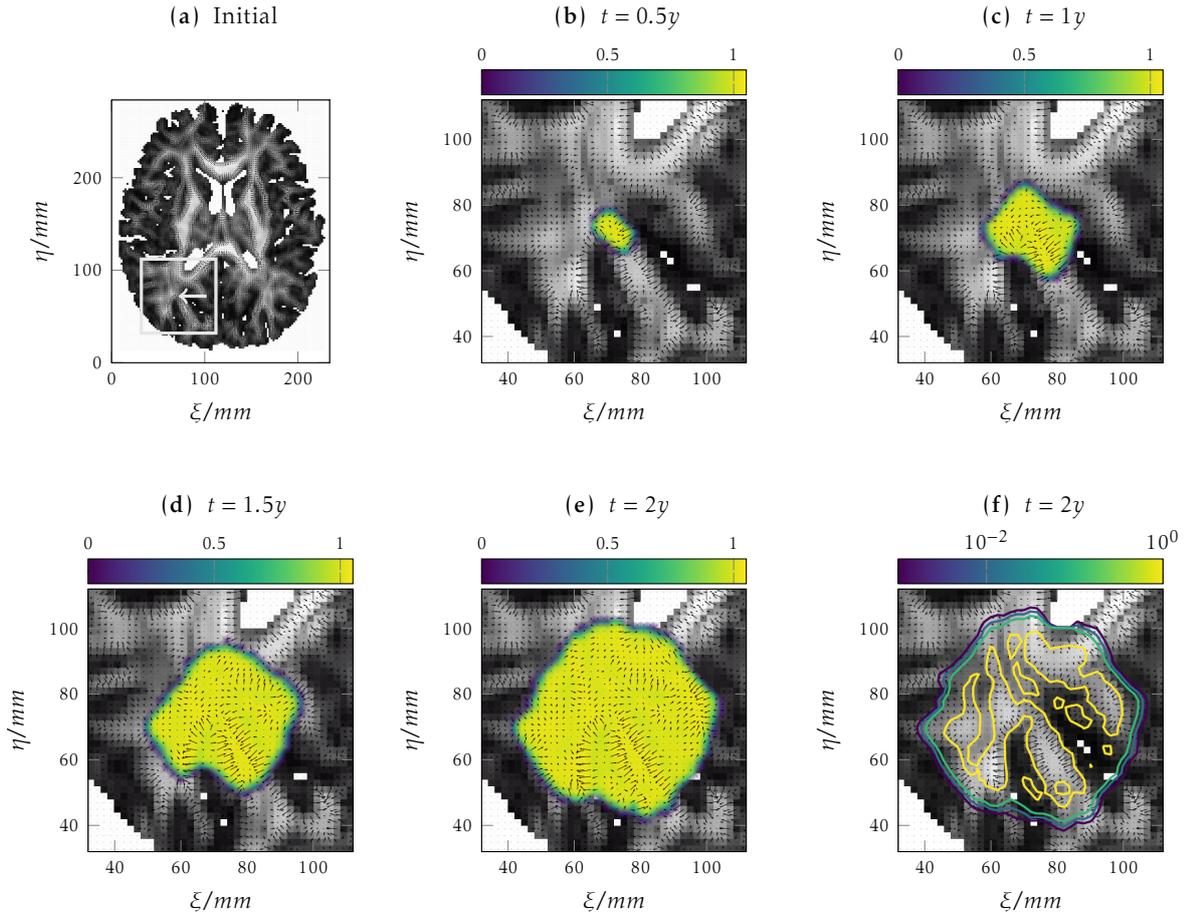

\begin{figure}
	\centering
	\def\localpath{figures/Brain/}
	\withfiguresize{\figuretwocollegend}{\figuretwocollegend}{\externaltikz{brain_comparisons}{\begin{tikzpicture}
\pgfplotscreateplotcyclelist{marklist}{%
	mark=\\
}
\pgfplotsset{
	cycle list/Set1-5,
	cycle multiindex* list={
		marklist\nextlist
		Set1-5\nextlist},
}

\def\colorplot#1#2#3#4#5{
	\nextgroupplot[
		xmin=32,
		xmax=112,
		ymin=32,
		ymax=112,
		xlabel={$\xcx / mm$},
		ylabel={$\xcy / mm$},
		title={\tikztitle{#1}},
		point meta min ={#2},
		point meta max ={#3},
		#4
	]
	\addplot graphics[xmin=32,xmax=112,ymin=32,ymax=112] {\localpath #5};
}

\begin{groupplot}[
	group style={group size=2 by 2, horizontal sep = \figurelegendhorizontalsep,  vertical sep = 2cm},
	axis on top,
	scale only axis,
	enlargelimits=false,
	width = \figurewidth,
	height =\figureheight,
	colormap name ={cmviridis},  
	scaled x ticks=false,
   	title style = {font=\figurelabelfont{\figurewidth}},
   	xlabel style = {font=\figurelabelfont{\figurewidth}}, 
   	ylabel style = {font=\figurelabelfont{\figurewidth}}, 
   	xticklabel style={font=\figureticklabelfont{\figurewidth}},
   	yticklabel style={font=\figureticklabelfont{\figurewidth}},
	legend style={at={(1.02,1.)},anchor=north west, font=\figurelabelfont{\figurewidth}},
]
\colorplot{$\epsilon$}{0}{1.05}{legend entries={$\pareps=\scinum{1}{-3}$, $\pareps=\scinum{1}{-4}$, $\pareps=\scinum{1}{-5}$,  $\pareps=\scinum{3.3}{-6}$, $\pareps=\scinum{1}{-10}$},}{compare_eps_4.png}
\addplot{0};
\addplot{0};
\addplot{0};
\addplot{0};
\colorplot{Discretization order}{0}{1.05}{legend entries={\mischeme[]{1}, \ivscheme[]{1}, \mischeme[]{2}, \ivscheme[]{2}},}{compare_implvol_4.png}
\addplot{0};
\addplot{0};
\addplot{0};
\colorplot{Moment Order}{0}{1.05}{legend entries={$P_1$, $P_3$},}{compare_pn_4.png}
\addplot{0};
\colorplot{Dimension}{0}{1.05}{legend entries={2D,3D},}{compare_dimension_4.png}
\addplot{0};
\end{groupplot}

 \end{tikzpicture}}}
	\settikzlabel{fig:brain-comparisons-a}
	\settikzlabel{fig:brain-comparisons-b}
	\settikzlabel{fig:brain-comparisons-c}
	\settikzlabel{fig:brain-comparisons-d}
	\caption{Results of glioma simulations(\secref{sec:brain}) with varied parameters and schemes. Shown are always the $10\%$ contours. For the interpretation of the background image, refer to \figref{fig:brain-evolution}. \ref{fig:brain-comparisons-a}: Solutions for various $\pareps$ in the intermediate to diffusive regime. \ref{fig:brain-comparisons-b}: Comparison between the numerical schemes. \ref{fig:brain-comparisons-c}: Comparison between moment orders. \ref{fig:brain-comparisons-d}: Comparison between the two-dimensional and three-dimensional models.}
	\label{fig:brain-comparisons}
\end{figure}
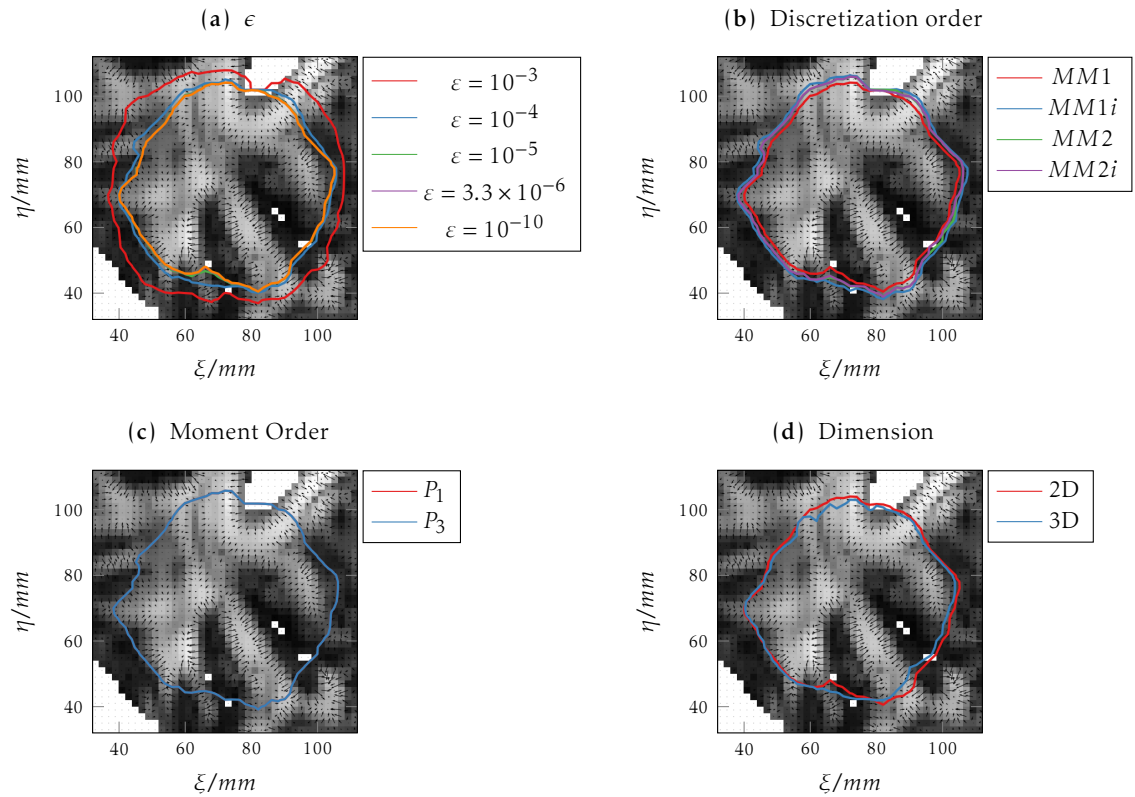

\section{Discussion}
\label{sec:discussion}
The goal of this work was to develop a numerical tool for a special class of transport equations that lead to an advection-diffusion-reaction equation in the parabolic limit. 
This method should be applicable to a wide range of scaling regimes, from almost free transport to very close to the diffusion limit. 
One example of an application that is very close to the parabolic limit is a model of glioma invasion in the human brain. 
The method was developed mainly with this model and the corresponding data in mind. 
This means that in the implementation, we could take advantage of the simplifications it offers; for example that the turning operator is explicitly invertible or that the equilibrium distribution is of a quadratic form. 
But probably the most significant influence on the method development came from the associated data. 
DTI data are measured and delivered on regular grids with fixed spatial resolution. 
On each grid cell, the water diffusion tensor is assumed constant, because there is no natural way to interpolate between those tensors. 
To avoid interpolation artifacts in the solution, the space discretization has to use the same grid as the original data. 
As a consequence, the method was implemented only for tensor-product grids and not tested for more general grids. 
However, the method has to address the strong heterogeneities and discontinuities of the DTI data. 

As a starting point for our scheme, we used the method developed by Lemou and Mieussens \cite{lemou2008ap}. 
This scheme employs a micro-macro decomposition and discretizes the microscopic and macroscopic components on different parts of a staggered grid.   
In this work, we generalized the method to an asymptotic preserving finite-volume formulation on primal-dual mesh pairs that works in two and three space dimensions.
In the description of the method, we used a mostly mesh-agnostic notation because we are confident that it also is applicable on unstructured meshes. 
Most parts of the implementation in DUNE \cite{dune-web-page} are already written mesh-independently, but a complete implementation is still only available for tensor-product grids. 
Development and testing of the unstructured implementation are left for the future. 

To discretize the velocity space in the micro equation, we employ the method of moments. 
More specifically, we use spherical harmonic basis functions and a linear reconstruction ansatz. 
In the diffusive regime, first-order basis polynomials are accurate enough, which means that only one degree of freedom per space dimension is needed.  
Compare this to the discrete ordinates method, that needs at least two degrees of freedom per space dimension to maintain symmetry. 
For successively less diffusive regimes, higher moment orders can be added as needed. 
Of course, in the kinetic regime the linear moment method has the usual drawback of producing unphysical Gibb's phenomena.
But this is not a problem in the diffusive regime. 

For asymptotic preserving methods, one special point of interest is resulting discretization in the parabolic limit. 
We show the limit diffusion and drift approximations only for a very simplified setting---regular grid with constant isotropic coefficients---but this is enough to identify two drawbacks of the basic method. 
First, the limit diffusion approximation is a five-point diagonal stencil that leads to a decoupling of grids and spurious oscillations. 
The same effect was also described in \cite{buet2012design} and seems to be a general problem for primal-dual discretizations. 
We propose alterations of the basic method that effectively allows us to modify the limiting discretization of the diffusion and drift terms. 
In effect, this leads to the classical five-point stencil for the diffusion and an upwind approximation of the drift. 
However, the drift discretization comes at the price of being inherently first-order accurate. 

We perform a wide range of benchmarks to numerically test some of the method's properties. 
The fundamental solution test demonstrates that the method indeed is asymptotic preserving and in the limit converges with the correct order to the fundamental solution. 
Moreover, we use this benchmark to estimate properties of the modified equation of the scheme. 
Of special interest is the behavior of the method in presence of strong discontinuities as encountered in the DTI data. 
For this, we adapt two stationary benchmark tests from the porous media community. 
The scheme deals well with strong jumps in permeability and has surprisingly a higher rate of convergence than could be expected from the regularity of the solution. 
Also, jumps in diffusion direction across an interface are resolved well, as long as the flux is only normal to the interface. 
Any tangential flux drastically reduces the approximation quality. 
Last but not least we demonstrate the capabilities of the method in the glioma invasion model. 
Although the parameters are only very rough estimates, the overall situation is similar to the application. 
The method performs well even on the coarse and heterogeneous real-world DTI data.

\newpage

\addreferencestotoc

\bibliographystyle{plain}
\bibliography{bibliography}

\end{document}